\documentclass{article}
\usepackage[T1]{fontenc}
\usepackage[utf8]{inputenc}
\usepackage{lmodern}
\usepackage{stmaryrd}
\usepackage[spanish,english]{babel}
\usepackage{geometry} 
\usepackage{amsmath,amssymb,amsfonts,amsthm,mathtools}
\usepackage{mathrsfs}
\usepackage{standalone}
\usepackage{graphicx}
\usepackage{caption}
\usepackage{subcaption}
\usepackage{array}
\usepackage{fancyhdr}
\usepackage{textcomp} 
\usepackage{xcolor} 
\usepackage{epigraph} 
\usepackage{comment}
\usepackage[colorlinks=true, urlcolor=purple, linkcolor=purple!87!black, citecolor=green!60!black, pdfborder={0 0 0} ]{hyperref}
\usepackage{bookmark}
\usepackage{url}
\usepackage{float}
\usepackage[style=numeric,sorting=nyt,backref=true,maxbibnames=99]{biblatex}
\usepackage{tikz}
\usetikzlibrary{calc,fadings,arrows,decorations.pathmorphing,backgrounds,positioning,fit,shapes,arrows,automata,matrix,fit,patterns,external,knots}
\usepackage{cleveref}

\theoremstyle{plain}
\newtheorem{theorem}{Theorem}[section] 
\newtheorem{prop}[theorem]{Proposition}
\newtheorem{lemma}[theorem]{Lemma}
\newtheorem{corollary}[theorem]{Corollary}
\theoremstyle{definition}
\newtheorem{remark}[theorem]{Remark}
\newtheorem{defn}[theorem]{Definition} 
\newtheorem{example}[theorem]{Example} 
\newtheorem{question}[theorem]{Question}

\theoremstyle{plain}
\newtheorem{thmx}{Theorem}

\newcommand{\N}{\mathbb{N}}
\newcommand{\Z}{\mathbb{Z}}

\newcommand{\NN}{\mathbb{N}}
\newcommand{\HH}{\mathbb{H}}
\newcommand{\ZZ}{\mathbb{Z}}
\newcommand{\FF}{\mathcal{F}}
\newcommand{\G}{\mathcal{G}}
\newcommand{\X}{\mathcal{X}}
\newcommand{\RR}{\mathbb{R}}
\newcommand{\define}[1]{\textbf{#1}}
\newcommand{\act}{\curvearrowright}
\newcommand{\orb}{\textnormal{orb}}
\newcommand{\DP}{\textnormal{DP}}
\newcommand{\SDSP}{\textnormal{SDSP}}

\makeatletter
\def\moverlay{\mathpalette\mov@rlay}
\def\mov@rlay#1#2{\leavevmode\vtop{%
   \baselineskip\z@skip \lineskiplimit-\maxdimen
   \ialign{\hfil$\m@th#1##$\hfil\cr#2\crcr}}}
\newcommand{\charfusion}[3][\mathord]{
    #1{\ifx#1\mathop\vphantom{#2}\fi
        \mathpalette\mov@rlay{#2\cr#3}
      }
    \ifx#1\mathop\expandafter\displaylimits\fi}
\makeatother

\newcommand{\fadeimage}[3][15pt]{%
  \begin{tikzpicture}
    \node[inner sep=0pt, outer sep=0pt] (img) {\includegraphics[#2]{#3}};
    \fill[white, path fading=south] (img.north west) rectangle ($(img.north east)+(0,-#1)$);
    \fill[white, path fading=north] (img.south west) rectangle ($(img.south east)+(0,#1)$);
    \fill[white, path fading=east]  (img.north west) rectangle ($(img.south west)+(#1,0)$);
    \fill[white, path fading=west]  (img.north east) rectangle ($(img.south east)+(-#1,0)$);
  \end{tikzpicture}
}

\newcommand{\Addresses}{{
		\bigskip
		\hskip-\parindent   S.~Barbieri, \textsc{Departamento de Matem\'{a}ticas, Universidad de Chile, Santiago, Chile.}\par\nopagebreak
		\textit{E-mail address}: \texttt{sbarbieri@uchile.cl}
		
		\medskip
		
		\hskip-\parindent   N.~Bitar, \textsc{
        Institute de Recherche en Informatique Fondamentale (UMR CNRS 8243), Universit\'{e} de Paris Cité, Paris, France.}\par\nopagebreak
		\textit{E-mail address}: \texttt{nicolas.bitar@irif.fr}

}}

\usepackage{morefloats}

\usepackage{scalerel}
\usepackage[ruled]{algorithm2e}

\usepackage{tikz}
\usetikzlibrary{arrows,decorations.pathmorphing,backgrounds,positioning,fit,petri,knots}

\usepackage{imakeidx}
\makeindex

\begin{document}

\author{Sebasti\'an Barbieri and Nicol\'as Bitar}

\date{}

\title{Space-filling groups}
\maketitle
\begin{abstract}
\hskip -.2in
\noindent

We introduce the notion of a space-filling group. These are finitely generated groups which admit an action on a compact space that carries a continuous cocycle which covers arbitrarily large portions of the space in a uniform manner. This notion is inspired by J. Kari's construction of a $\ZZ^2$-tileset that simulates Hilbert curves and its application to reversibility problems on cellular automata.

We show that the property of being space-filling is always witnessed by a subshift, and consider computable, finitary and local variants. We show that the property of being space-filling is an invariant of quasi-isometry, and that the computable and finitary variants are also invariant under natural hypotheses. Furthermore, we provide interesting examples within the class of groups which are quasi-isometric to the hyperbolic space $\mathbb{H}^d$ for some $d\geq 1$ along with several stability properties of the class and applications to cellular automata and snake tiling problems on groups.

\medskip
        
        
        \textit{Secondary:} 20F67, 
                20F10. 
\end{abstract}


\section{Introduction}

A type of construction which is found nowadays quite often in the study of symbolic dynamics on groups is a shift space $X$ on a group $G$, in which every configuration $x\in X$ encodes the structure of a group $H$ in a geometric manner while respecting a given set of constraints. This is particularly fruitful in the case where $H$ is not a subgroup of $G$.

A well-known construction that follows this general idea is the encoding of a group $H$ through a subshift $X$ on $G$ in a quasi-isometric way. More precisely, if $H$ and $G$ are finitely generated quasi-isometric groups, Cohen~\cite{Cohen_2017} constructs a $G$-subshift $X$ in which every configuration encodes a quasi-isometry $f\colon G \to H$ that can be decoded locally. This is done in such a way that if both groups are finitely presented, the subshift $X$ can be described by local rules and thus is a subshift of finite type (SFT). This construction provides a useful tool to encode $H$-actions into a $G$-subshift, and it has been used to prove the quasi-isometry invariance (for finitely presented groups) of the domino problem~\cite{Cohen_2017}, the existence of strongly aperiodic SFTs~\cite{Cohen_2017}, the property of self-simulation~\cite{BarSabSal2026}, and the existence of nontrivial Medvedev degrees~\cite{BarCar2026_medvedev}. A variation on more general structures has also been introduced by the authors to show the undecidability of the domino problem for geometric tilings with finite local complexity~\cite{BarBit2026}.

A second important variant of these constructions was introduced by Jeandel~\cite{jeandel2015translationlikeactionsaperiodicsubshifts}. Here, we consider two finitely generated groups $H,G$ and suppose that $H$ acts translation-like on $G$, that is, freely and such that every generator of $H$ moves elements a bounded distance (see~\cite{Whyte1999}). In this construction, every configuration in the $G$-subshift $X$ must encode a translation-like action locally. Hence we may think of this subshift as a space of covers of a Cayley graph of $G$ by copies of a Cayley graph of $H$. This technique has been used to construct examples of weakly aperiodic subshifts~\cite{jeandel2015translationlikeactionsaperiodicsubshifts}, simulation theorems on direct products~\cite{Barbieri2019}, to study topological entropies of SFTs on general groups~\cite{Bar2021_entropy_fails} and to study Medvedev degrees on effective subshifts~\cite{Carrasco2026} and SFTs~\cite{BarCar2026_medvedev} among other results.

In the case of $H=\ZZ$, we can think of Jeandel's construction as a $G$-subshift in which every configuration encodes disjoint copies of bi-infinite paths on a Cayley graph of $G$. In this setting, it would be desirable to have a subsystem in which a single copy of $\ZZ$ covers the whole space, that is, a subshift which encodes the space of Hamiltonian paths on the Cayley graph of $G$ with respect to a set of generators (these Hamiltonian paths exist on some Cayley graph of any finitely generated one-ended group $G$ by a result of Seward~\cite{Seward2014}). However, the space of configurations which encode  Hamiltonian paths is not closed unless the group $G$ is virtually cyclic (\Cref{prop:transitive_implies_virt_cyclic}) and thus this property is too strong.

Instead of pursuing the property above, we settle for a slightly different property that encodes the idea that every path must cover translated copies of every ball of the group $G$ in a uniform manner. This turns out to be a shift-invariant and closed property, and thus one that can be encoded by subshifts. 

Let us proceed more generally. Let $G$ be a discrete group and $X$ be a nonempty compact metrizable space. Consider an action $G\curvearrowright X$ by homeomorphisms. Given a continuous map $R\colon X\to G$ and $x\in X$, we denote by $R_x\colon G\to G$ the induced map given by $R_x(g)=gR(g^{-1} x)$. For $n \geq 0$ we inductively define $R^0_x(g) = g$ and $R^n_x(g) = R_x(R_x^{n-1}(g))$. 

\begin{defn}\label{defn:space_filling_intro}
    An action $G\curvearrowright X$ by homeomorphisms is \define{space-filling} if it admits a continuous map $R\colon X \to G$ such that for every finite $F\subset G$ there exists $n \in \NN$ such that for all $x\in X$, there is $g\in G$ such that
  \[gF \subseteq\left\{R_x^i(1_G) : 0\leq i\leq n\right\}.\]
  We say that a group $G$ is \define{space-filling} if it admits some space-filling action $G\curvearrowright X$.
\end{defn}

Notice first that as $R$ is continuous and $G$ is discrete, its image is a finite subset of $G$. In the case where $X$ is a $G$-subshift and $G\curvearrowright X$ is the left shift action, then one can interpret $R$ as a map which depends upon finitely many coordinates and provides a direction. The condition above says that starting from any position, following the directions dictated by $R$ will cover translates of any finite subset of the group. In particular, a space-filling group must be finitely generated (see~\Cref{prop:spacefilling_isFG}).

Our definition is inspired by a construction of Jarkko Kari~\cite{Kari_reversibility_defofpeanocurve}, in which a $\ZZ^2$-tileset is constructed with the property that any tiling by it encodes (through a sliding block code) arbitrarily large copies of the Hilbert curve as shown in~\Cref{fig:el_jilbert}. This construction was used by Kari to provide a proof of the undecidability of the reversibility problem for cellular automata on $\ZZ^2$. In particular, Kari's construction provides an example of a space-filling $\ZZ^2$-SFT.

\newcommand{\hilbertstep}[2]{
    \draw[thick, blue, -stealth] (\hx,\hy) -- ++(0.5*#1,0.5*#2);
    \draw[thick, blue] (\hx+0.5*#1,\hy+0.5*#2) -- ++(0.5*#1,0.5*#2);
    \pgfmathparse{\hx+#1}\xdef\hx{\pgfmathresult}
    \pgfmathparse{\hy+#2}\xdef\hy{\pgfmathresult}}
\expandafter\def\csname hilbert@R\endcsname{\hilbertstep{1}{0}}
\expandafter\def\csname hilbert@L\endcsname{\hilbertstep{-1}{0}}
\expandafter\def\csname hilbert@U\endcsname{\hilbertstep{0}{1}}
\expandafter\def\csname hilbert@D\endcsname{\hilbertstep{0}{-1}}

\begin{figure}[ht!]
    \centering
    \begin{tikzpicture}[scale = 0.4]
        \draw[opacity = 0.5] (0,0) grid (16,16);
        \clip (0,0) rectangle (16,16);
        \draw[thick, blue] (0.5,0) -- (0.5,0.5);
        \xdef\hx{0.5}\xdef\hy{0.5}
        \foreach \m in {%
            R,U,L,U,U,R,D,R,U,R,D,D,L,D,R,R,U,R,D,R,R,U,L,U,R,U,L,L,D,L,U,U,%
            U,R,D,R,R,U,L,U,R,U,L,L,D,L,U,L,L,D,R,D,D,L,U,L,D,L,U,U,R,U,L,U,%
            U,R,D,R,R,U,L,U,R,U,L,L,D,L,U,U,R,U,L,U,U,R,D,R,U,R,D,D,L,D,R,R,%
            R,U,L,U,U,R,D,R,U,R,D,D,L,D,R,D,D,L,U,L,L,D,R,D,L,D,R,R,U,R,D,R,%
            U,R,D,R,R,U,L,U,R,U,L,L,D,L,U,U,R,U,L,U,U,R,D,R,U,R,D,D,L,D,R,R,%
            R,U,L,U,U,R,D,R,U,R,D,D,L,D,R,D,D,L,U,L,L,D,R,D,L,D,R,R,U,R,D,D,%
            L,D,R,D,D,L,U,L,D,L,U,U,R,U,L,L,D,L,U,L,L,D,R,D,L,D,R,R,U,R,D,D,%
            D,L,U,L,L,D,R,D,L,D,R,R,U,R,D,R,R,U,L,U,U,R,D,R,U,R,D,D,L,D,R,R}{\csname hilbert@\m\endcsname}
    \end{tikzpicture}
    \caption{The Hilbert curve encoded by Kari's tileset.}
    \label{fig:el_jilbert}
\end{figure}
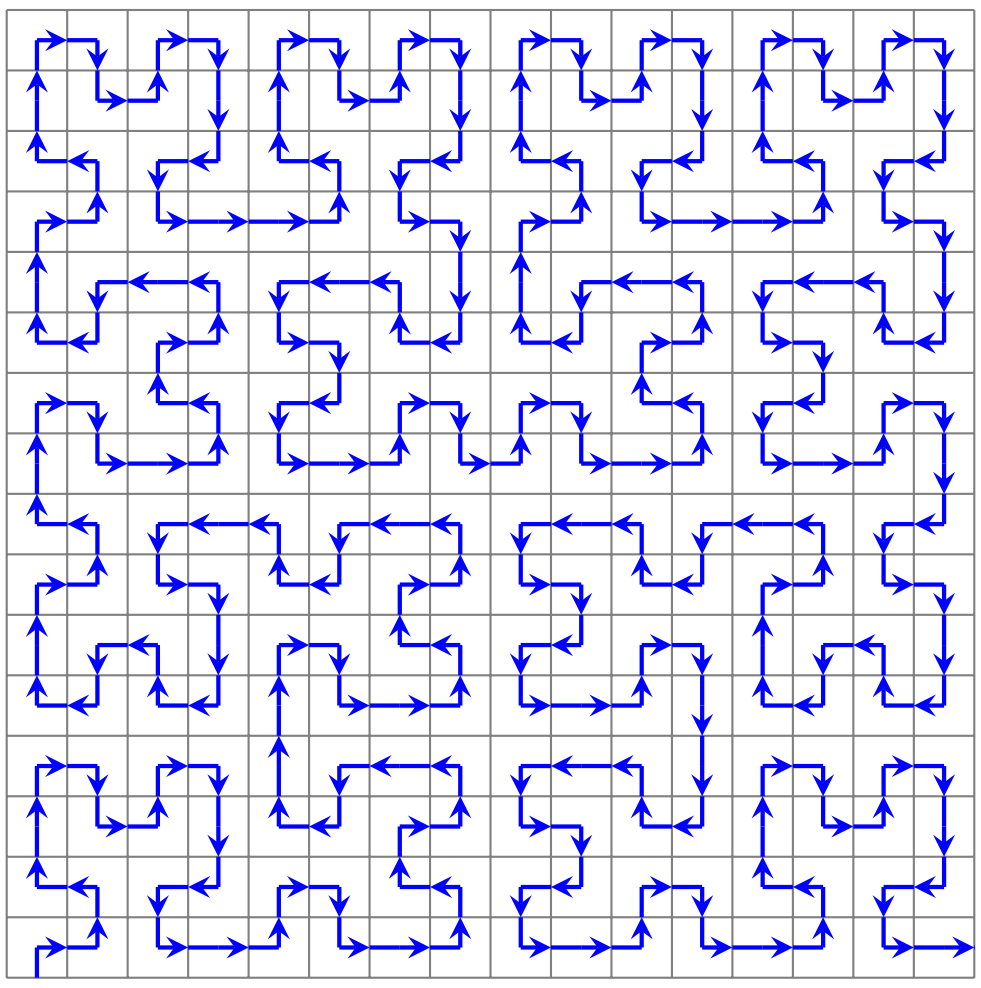

The first fundamental observation on space-filling groups is that the property is always witnessed by subshifts. 

\begin{thmx}\label{mainthm:los_subshifts_son_magnificos}
    A group is space-filling if and only if it admits a space-filling subshift.
\end{thmx}

This allows us to restrict our view to this class and think of these objects symbolically: up to topological conjugacy the continuous map $R$ from~\Cref{defn:space_filling_intro} can be thought of as a map whose domain is the alphabet of the subshift, and each symbol provides a direction. This recovers the behavior of the Hilbert curve, where following local directions one is able to cover arbitrarily large portions of the plane.

Within this framework, we consider three variants of space-filling groups which are given by considering increasingly harsh constraints on the underlying subshift. These variants are collected in the following definition.

\begin{defn}
    Let $G$ be a finitely generated group. 
    \begin{enumerate}
        \item $G$ is \define{computably space-filling} if it admits an effective space-filling subshift (\Cref{def:comp_SF}).
        \item $G$ is \define{finitary space-filling} if it admits a space-filling subshift of finite type (\Cref{def:finitary_SF}).
        \item $G$ is \define{locally space-filling} if there exists a finite symmetric set of generators $S$ of $G$ for which there exists a cycle-free strongly space-filling directed tileset over $(G,S)$ (\Cref{def:locally_SF}). 
    \end{enumerate}
\end{defn}

For these classes we study natural obstructions, such as having infinitely many ends (\Cref{prop:inf_ends}) and provide a few stability and realization results. Most notably, we show that the properties of being space-filling, computably space-filling and finitary space-filling are stable under direct products of infinite groups (\Cref{cor:productos}). We also study variants of the space-filling property where we ask that the maps $R_x$ are invertible and that they induce transitive actions, and state some results based on simulation theory that give tools to turn computably space-filling groups into finitary space-filling groups.

Our first major theorem is that these three properties are invariant under quasi-isometries under natural conditions.

\begin{thmx}\label{mainthm:QI}
    Let $G,H$ be finitely generated groups which are quasi-isometric. The following hold.
    \begin{enumerate}
        \item $G$ is space-filling if and only if $H$ is space-filling (\Cref{cor:QI_groups}).
        \item If both $G,H$ have decidable word problem, then $G$ is computably space-filling if and only if $H$ is computably space-filling (\Cref{cor:QI_effective}).
        \item If both $G,H$ are finitely presented, then $G$ is finitary space-filling if and only if $H$ is finitary space-filling (\Cref{thm:QI_local}).
    \end{enumerate}
\end{thmx}

A direct consequence of these results and the product stability theorem is that every finitely generated virtually abelian group is finitary space-filling (\Cref{cor:virtab_is_finitary_SF}).

Beyond virtually abelian groups, our first examples of computably space-filling groups are those which are quasi-isometric to a hyperbolic space.

\begin{thmx}\label{mainthm:hyperbolic}
    Every finitely generated group which is quasi-isometric to the hyperbolic space $\HH^d$ for $d\geq 1$ is computably space-filling.
\end{thmx}

Of course, the interesting case is $d\geq 2$, as for $d=1$ we get precisely the infinite virtually cyclic groups which were already covered by our previous results. In particular, our theorem applies to every group $G$ that acts properly and cocompactly by isometries on $\HH^{d}$, such as cocompact Fuchsian groups and the fundamental groups of closed hyperbolic $d$-manifolds. We also remark that in an upcoming article of the first author with Salo and Sablik, it will be shown that groups which are quasi-isometric to $\HH^d$ for $d\geq 3$ are self-simulable in the sense of~\cite{BarSabSal2026}, and thus this would imply that they are finitary space-filling as well. In the meantime, we obtain new examples of finitary space-filling groups through simulation theorems; for instance, if $G$ is quasi-isometric to $\HH^d$ for some $d\geq 2$, then $\ZZ^n \times G$ is finitary space-filling for any $n\geq 1$ (\Cref{ex:productitos}).

Finally, we turn towards applications of this theory to cellular automata and snake tiling problems in groups. For this purpose, we introduce the strongest notion of space-filling of the article. We say a directed tileset with respect to a generating set $S$ (that is, a tileset where every letter additionally contains an element of $S$ pointing towards the next tile) is \define{strongly space-filling} if it admits a valid configuration whose path is injective, and the path of every infinite snake of the tileset covers translates of every finite subset of the group (\Cref{def:SSF}). A group is \define{locally space-filling} if, for some generating set, it admits such a directed tileset which is moreover \define{cycle-free}, meaning that no configuration has a valid directed cycle. For these groups and notions, we can establish the following reductions.

\begin{thmx}\label{mainthm:reduction}
     Let $G$ be a finitely generated group and $S$ a finite symmetric set of generators of $G$ such that there is a cycle-free strongly space-filling directed tileset over $(G,S)$. The following hold.
    \begin{itemize}
        \item The domino problem many-one reduces to the strong directed snake problem with respect to $S$ (\Cref{thm:reduction}).
        \item The complement of the domino problem many-one reduces to the injectivity and reversibility problems of cellular automata over $G$ (\Cref{thm:CA}).
    \end{itemize}
\end{thmx}

We note that \Cref{thm:reduction} in fact holds without the assumption that the underlying strongly space-filling tileset is cycle-free and that~\Cref{thm:CA} is the natural generalization of a proof of J. Kari to this general setting~\cite{Kari_reversibility_defofpeanocurve}.

\subsubsection*{Structure of the article}
In~\Cref{sec:space-filling} we introduce basic notation, provide the new definitions and study the general properties of space-filling actions. This section begins with a study of the associated movement functions and how the property is preserved under topological extensions and subsystems. Then we turn to space-filling subshifts and show that they are the quintessential examples of space-filling actions, justified by~\Cref{mainthm:los_subshifts_son_magnificos}. Next, we study two-sided space-filling actions, followed by group theoretical restrictions on space-filling groups. The section then ends with the theorem which establishes stability under direct products, and the study of computably space-filling groups.

\Cref{sec:QI} is devoted to the quasi-isometry invariance of space-filling (and its computable and finitary variants) under suitable conditions (\Cref{mainthm:QI}). We do this by enlarging the scope of space-filling actions to countable graphs of bounded degree, and providing a combinatorial characterization that is preserved under quasi-isometries. The section ends by showing that the finitary space-filling property is a quasi-isometry invariant for finitely presented groups using techniques of Cohen~\cite{Cohen_2017}. 

In~\Cref{sec:hyp-group} we construct the first examples of computably space-filling groups beyond virtually abelian groups by proving \Cref{mainthm:hyperbolic}. More precisely, we show that a specific graph $\mathcal{B}_{d+1}$ which is quasi-isometric to $\HH^{d+1}$ satisfies the combinatorial characterization of space-filling given in~\Cref{sec:QI}. The cases $d=1$ and $d\geq 2$ are treated separately as the constructions are slightly different. 

\Cref{sec:applications} is independent of the two previous sections. Here we introduce the strongest notion of space-filling, locally space-filling groups. With this notion, we study how different decision problems relating to tilings are related to the domino problem of the group and prove~\Cref{mainthm:reduction}.

Finally, we include~\Cref{sec:questions} where we present a list of questions that arose from our study of space-filling actions.

\subsubsection*{Acknowledgements} S. Barbieri was supported by ANID Fondecyt regular 1240085. N.~Bitar was supported by the ANR project IZES ANR-22-CE40-0011 and by the ERC grant DynAMiCs (101167561) of the European Research Council.

\subsubsection*{AI declaration} The final version of this article was completely written by the authors. Claude was used to draft initial versions of the quasi-isometry invariance and hyperbolic spaces proofs using our preliminary proof sketches of the results as input. In particular, it allowed us to simplify and improve the proofs in~\Cref{sec:hyp-group} considerably. These sections were then rewritten by the authors. AI was also used to produce a review of a preliminary version of the article, which we used to improve this final version and correct several mistakes. As it should always be, the authors are entirely responsible for the contents of this article.

\section{Space-filling actions}
\label{sec:space-filling}

From our analysis of the Hilbert-curve subshift, we propose the following elements for the definition of a space-filling subshift, and more generally a space-filling action. Consider a discrete group $G$ and a nonempty compact metrizable space $X$, such that the group acts on $X$ by homeomorphisms. In this general context, the directions that must be followed at each step in the Hilbert-curve subshift are understood as a continuous map $R\colon X \to G$.

Given such a map $R\colon X \to G$, we denote by $R_{x}\colon G \to G$ the induced map given by $R_x(g) = gR(g^{-1}x)$. For $n \geq 0$ we inductively define $R^0_x(g) = g$ and $R^n_x(g) = R_x(R_x^{n-1}(g))$. In the case where $R_x$ is invertible we extend the above definition naturally to negative integers using the inverse and note that it defines a $\ZZ$-cocycle. Finally, given a subset $I\subseteq\ZZ$ we use the notation $R^I_x(g) = \{R^k_x(g) \mid k\in I\}$. For $a \leq b \in \ZZ$ we will write $\llbracket a,b\rrbracket = [a,b]\cap \ZZ$ for the interval of integers between $a$ and $b$.

For a group action $G\act X$ and $x\in X$ we write $\orb(x) = \{g\cdot x\mid g\in G\}$ and $\operatorname{stab}(x) = \{g \in G : gx =x \}$ for the orbit and stabilizer of $x$ respectively. Below we note a few identities that will be used repeatedly throughout the text. We omit their proof as they follow directly from the definitions above.

\begin{prop}
\label{prop:faciles}
    Let $G\curvearrowright X$ be an action and let $R\colon X \to G$. For every $x\in X$, $g,h \in G$ and $k \geq 0$ the following hold:
    \begin{enumerate}
        \item $R^k_x(gh) = gR^k_{g^{-1}x}(h)$. Moreover, if $R_y$ is invertible for every $y\in\orb(x)$, this holds for every $k\in\Z$.
        \item If $R_x$ is invertible, $m\geq 0$, and $g=R^k_x(1_G)$ for some $k\in\llbracket 0,m\rrbracket$, then $R^{\llbracket 0,m\rrbracket}_x(1_G)\subseteq R^{\llbracket -m,m\rrbracket}_x(g)$. 
    \end{enumerate}
\end{prop}

We can now proceed to the main definition of the article.

\begin{defn}
    An action $G\curvearrowright X$ by homeomorphisms is \define{space-filling} if it admits a continuous map $R\colon X \to G$ such that for every $F\Subset G$ there exists $n \in \NN$ such that for all $x\in X$, there is $g\in G$ such that
  \[gF \subseteq R_x^{\llbracket 0,n\rrbracket}(1_G).\]
  We say that $R$ is a \define{movement function} associated to $G\curvearrowright X$. A group $G$ is \define{space-filling} if it admits a space-filling action.
\end{defn}

Intuitively, a space-filling action is one for which the movement function covers translates of every finite set in a uniform manner, and thus, up to translation, it fills the group. Our first fundamental observation is that only finitely generated groups may admit space-filling actions.

\begin{prop}\label{prop:spacefilling_isFG}
    If a group $G$ is space-filling, then $G$ is finitely generated. More precisely, for every movement function $R$ the finite set $R(X)$ generates $G$.
\end{prop}

\begin{proof}
    Suppose $G$ admits a space-filling action $G\curvearrowright X$ and let $R\colon X \to G$ be a movement function. As $R$ is continuous and $X$ is compact, it follows that $R(X)$ is a compact subset of the discrete group $G$ and thus $R(X)$ is finite. Therefore for any $x\in X$ we have that \[ \{R_x^k(1_G) : k \in \NN\} \subset \langle R(X)\rangle.  \]
    If $G$ is not finitely generated, it follows that $G \setminus \langle R(X)\rangle\neq \varnothing$. Let $g \in G \setminus \langle R(X)\rangle$ and take $F=\{1_G,g\}$. It is immediate that no translate of $F$ may lie in $\langle R(X)\rangle$, which contradicts the assumption that $G\curvearrowright X$ is space-filling. 
\end{proof}

From the previous proposition, we may restrict our study to the case of finitely generated groups. Given a finitely generated group $G$ and a finite symmetric set of generators $S$, we will denote the associated word metric by $|\cdot|_S$ and the ball of radius $n$ centered in the origin by $B_S(n)=\{g \in G: |g|_S \leq n\}$. We will often let the underlying set of generators be implicit and just write $B(n)$.

\begin{remark}\label{rem:covering_function}
    Let $G$ be a finitely generated group endowed with a word metric and let $G\curvearrowright X$ be a space-filling action. From the definition it follows that there is a map $f\colon \NN\to \NN$ such that for all $n\in\NN$ and $x \in X$ there is $g \in G$ such that \[ gB(n) \subset \{R_x^i(1_G) \mid 0 \leq i \leq f(n)\}.  \] We call such a map $f$ a \define{covering function} for $R$. Notice that any function $f'\colon\N\to\N$ that upper-bounds $f$ is also a covering function. Moreover, for every $x\in X$ and $i\in\NN$, if we let $w = R^i_x(1_G)$ then $R^{i+j}_x(1_G) = wR^j_{w^{-1}x}(1_G)$ for every $j\geq 0$ by~\Cref{prop:faciles}. Applying the above to $w^{-1}x$ we obtain that $R_x^{\llbracket j,j+f(n)\rrbracket}(1_G)$ contains a translate of $B(n)$ for every $x\in X$ and $j,n\in\NN$.
\end{remark}

\begin{remark}
\label{rem:detalpalotalastilla}
If $G\curvearrowright X$ is a space-filling action and $Y\subset X$ is a nonempty closed $G$-invariant set, then the induced action $G\curvearrowright Y$ is also space-filling with the same movement function. In particular if a space-filling action exists for $G$, then we may assume it is minimal.
\end{remark}

A natural question that arises from our definition is if it carries meaningful dynamical information, that is, if the 
property of an action being space-filling is invariant under isomorphisms of dynamical systems. 

Given continuous actions $G\curvearrowright X$ and $G\curvearrowright Y$ on compact spaces $X,Y$, a continuous map $\phi\colon X \to Y$ is a \define{topological morphism} if it is $G$-equivariant, that is, if $\phi(g\cdot x)=g\cdot\phi(x)$ for all $x\in X$ and $g\in G$. A topological morphism $\phi$ is called a \define{topological factor map} if it is surjective, and a \define{topological conjugacy} if it is a homeomorphism. Finally, if $\phi\colon X \to Y$ is a topological factor map, we say that $G\curvearrowright X$ is a \define{topological extension} of $G\curvearrowright Y$, and that $G\curvearrowright Y$ is a \define{topological factor} of $G\curvearrowright X$.
 
\begin{prop}\label{prop:top_extension_stable}
    The property of an action being space-filling is stable under topological extensions.
\end{prop}

\begin{proof}
    Let $\phi \colon X \to Y$ be a topological factor map for $G$-actions on $X,Y$. Suppose $G\curvearrowright Y$ is space-filling and let $R\colon Y \to G$ be a movement function. Define $R'\colon X \to G$ by $R'(x) = R(\phi(x))$.  

    As $\phi$ is continuous, it follows that $R'$ is continuous. Furthermore, as $\phi$ is $G$-equivariant, we have that for every $g \in G$ \[R'_{x}(g) = gR'(g^{-1}x)= gR(\phi(g^{-1}x)) = gR(g^{-1}\phi(x))= R_{\phi(x)}(g) .\]

    Hence, $R'_x = R_{\phi(x)}$, as well as all iterations. Furthermore, any covering function for $R$ is a covering function for $R'$. We conclude that $G\curvearrowright X$ is space-filling.
\end{proof}

\begin{corollary}
    The property of an action being space-filling is invariant under topological conjugacy.
\end{corollary}

\subsection{Space-filling subshifts}

Of particular interest for us is the case when the space-filling action corresponds to a shift space on a group $G$, as we will later show that without loss of generality we may always assume this is the case. Given a finite set $A$, we denote by $A^G$ the space of all maps from $G$ to $A$. We endow $A^G$ with the prodiscrete topology and with the left shift action $G\curvearrowright A^G$ given by \[ (g\cdot x)(h) = x(g^{-1}h) \mbox{ for every } x \in A^G \mbox{ and } g,h \in G.   \]

\begin{defn}
    A \define{subshift} or \define{shift space} is a subset $X\subset A^G$ which is closed and invariant under the left shift action.
\end{defn}

Given a finite set $F\subseteq G$, a pattern is a map $p\colon F \to A$. The cylinder set associated to a pattern $p$ is given by \[[p]= \{ x \in A^G \mid x(g)=p(g), \mbox{ for all } g \in F\}.\]
It is well known that $X\subseteq A^G$ is a subshift if and only if there exists a set of forbidden patterns $\FF$ such that $x\in X$ if and only if for every $g\in G$ and $p \in \FF$ we have $gx\notin [p]$. We denote the subshift defined by the forbidden patterns $\FF$ by $\X_{\FF}$. This equivalence gives us a meaningful way to introduce some interesting classes of subshifts depending on the structure of their defining set of forbidden patterns.

\begin{defn}
    A subshift $X\subseteq A^G$ is a \define{subshift of finite type} (SFT) if there exists a finite set of forbidden patterns $\FF$ such that $X = \X_{\FF}$.  
\end{defn}

A useful subclass of SFTs is those defined by patterns with very simple supports. An SFT $X$ defined by forbidden patterns $\FF$ is a \define{nearest neighbor SFT} if for every pattern $p\in \FF$ there exists a generator $s\in S$ such that $p$ has support $\{1_G, s\}$. A standard argument shows that every SFT is topologically conjugate to a nearest neighbor SFT.

\begin{defn}
    A subshift $Y$ is \define{sofic} if there exists an SFT $X$ and a factor map $\pi\colon X\to Y$.
\end{defn}

When working with subshifts, topological morphisms have a very rigid structure given by the Curtis--Hedlund--Lyndon theorem.
\begin{theorem}[\cite{hedlund1969endomorphisms}]
    Consider two subshifts $X\subseteq A^G$ and $Y\subseteq B^G$, and a map $\phi\colon X\to Y$. Then, $\phi$ is continuous and $G$-equivariant if and only if there exists $F\Subset G$ and a local map $\Phi\colon A^F\to B$ such that for all $x\in X$, $g\in G$ we have $\phi(x)(g) = \Phi((g^{-1}\cdot x)|_{F})$.
\end{theorem}

In the case where $X=Y$ the Curtis--Hedlund--Lyndon theorem states that topological endomorphisms correspond to \define{cellular automata}, that is, to maps $\phi\colon X \to X$ that are obtained by applying a local rule uniformly on the space.

A standard argument shows that for a sofic shift $Y\subseteq B^G$ we can always assume that there exists a nearest neighbor SFT $X\subseteq A^G$ and a local map $\Psi\colon A\to B$ such that its induced topological morphism is a factor map.

The following corollary is a direct consequence of~\Cref{prop:top_extension_stable}.
\begin{corollary}
\label{cor:conjugacy_invariant}
    If $G$ admits a sofic space-filling subshift, then it admits a space-filling SFT.
\end{corollary}

Although it seems subshifts will merely provide particular examples of space-filling actions, we will later show that in fact whenever the group is space-filling, it admits a space-filling subshift (\Cref{thm:todo_es_subshift}) and thus we can restrict our study to this case.

Finally, we give a name to the groups that admit a space-filling SFT.

\begin{defn}\label{def:finitary_SF}
    A group $G$ is \define{finitary space-filling} if it admits a subshift of finite type which is space-filling.
\end{defn}

Below we will provide examples of finitary space-filling groups that are either in the literature or very easy to prove. Later on we will show that every finitely generated and virtually abelian group is finitary space-filling.

\begin{prop}\label{prop:finite_groups}
    Finite groups are finitary space-filling.
\end{prop}

\begin{proof}
    As every subshift on a finite group is of finite type, it suffices to construct a space-filling subshift.
    
    Take any enumeration $g_0,\dots,g_{|G|-1}$ of $G$ and let $x\in G^G$ be given by $x(g_i)=g_i^{-1}g_{i+1 \bmod{|G|}}$. Take the subshift $X = \{ gx : g \in G\}$ (which is closed because $G$ is finite) and remark that $R\colon X \to G$ given by $R(y) = y(1_G)$ is a movement function which makes $X$ space-filling.
\end{proof}

\begin{prop}\label{prop:Z_and_Z2_locallyspacefilling}
    $\ZZ^d$ is finitary space-filling for $d \in \{1,2\}$.
\end{prop}

\begin{proof}
    For $d=1$ it suffices to consider the singleton $X = \{0^{\ZZ} \}$ and the map $R\colon X \to \ZZ$ given by $R(0^{\ZZ}) = 1$. For $d=2$ the construction of Kari~\cite{Kari_reversibility_defofpeanocurve} is a space-filling $\ZZ^2$-SFT.
\end{proof}

\subsection{From actions to subshifts}
Our goal in this subsection is to show that the existence of a space-filling action implies the existence of a space-filling subshift. To do this, we introduce a subshift that will be useful throughout the rest of the article.

Fix $k\geq 1$, $f\colon\N\to\N$ and a finite symmetric generating set $S$ of $G$. Given $W\Subset G$, $p \in B(k)^{W}$ and $g \in W$ we define $R_p(g) = gp(g)$, and inductively write $R^n_p(g) = R_p(R^{n-1}_p(g))$ whenever $R^{n-1}_p(g)\in W$. We say a pattern $p\in B(k)^{B(kf(n))}$ is \define{$(k,f)$-admissible} for $n$ if there exists $g\in G$ such that $gB(n) \subseteq R_p^{\llbracket 0, f(n)\rrbracket}(1_G)$.

\begin{defn}
     Consider $k\geq 1$ and $f\colon\N\to\N$. We define the \define{$(k,f)$-covering} subshift by
    \[\mathfrak{X}(k,f) = \{x\in B(k)^G : \forall n\geq 0, \ \forall g\in G, \ g^{-1}x|_{B(kf(n))} \textnormal{ is $(k,f)$-admissible for $n$}\}.\]
\end{defn}

We note that $\mathfrak{X}(k,f)$ implicitly also depends upon the choice of generators $S$. By its definition, this subshift witnesses the space-filling property whenever it is nonempty. We note this in the following lemma.
\begin{lemma}
\label{lem:subshift}
    For any $k\geq 1$ and $f\colon \N\to\N$, if the subshift $\mathfrak{X}(k,f)$ is nonempty, then it is space-filling with movement function $x\mapsto x(1_G)$ and covering function $f$.
\end{lemma}
\begin{proof}
    Denote $X = \mathfrak{X}(k,f) \neq \varnothing$ and fix $x\in X$. Set $R\colon X \to G$ by $R(x)=x(1_G)$ and note that with this definition we have that $R_x(g) = gR(g^{-1}x) = gx(g)$. Take $n\in\N$. Then, $x|_{B(kf(n))}$ is $(k,f)$-admissible for $n$, that is, there exists $g\in G$ such that \[gB(n)\subseteq R^{\llbracket 0 ,f(n) \rrbracket}_{x|_{B(kf(n))}}(1_G) = R^{\llbracket 0 ,f(n) \rrbracket}_{x}(1_G).\] Therefore, $X$ is space-filling.
\end{proof}

We can now prove the following result which implies~\Cref{mainthm:los_subshifts_son_magnificos}.
\begin{theorem}
\label{thm:todo_es_subshift}
    Let $G\curvearrowright Y$ be a space-filling action with movement function $R$. Let $S$ be a finite symmetric set of generators of $G$ which contains $R(Y)$ and let $f$ be a covering function for $R$ with respect to $S$. Then the $(1,f)$-covering subshift is nonempty. In particular, every space-filling group admits a space-filling subshift.
\end{theorem}

\begin{proof}
      Let $X = \mathfrak{X}(1,f)\subseteq B(1)^G$. Fix $y\in Y$ and define $x\in B(1)^G$ by $x(h) = R(h^{-1}y)$ for every $h\in G$. Then $R_x(h) = hx(h) = R_y(h)$ for every $h \in G$. Let $n\in\N$ and $g\in G$. By~\Cref{prop:faciles} we have $R^j_{g^{-1}x}(1_G) = g^{-1}R^j_x(g) = g^{-1}R^j_y(g) = R^j_{g^{-1}y}(1_G)$ for every $j\geq 0$. As $f$ is a covering function, there exists $t\in G$ such that
      \[tB(n) \subseteq R^{\llbracket 0, f(n)\rrbracket}_{g^{-1}y}(1_G) = R_{g^{-1}x}^{\llbracket 0, f(n)\rrbracket}(1_G).\]
      Therefore $(g^{-1}x)|_{B(f(n))}$ is $(1,f)$-admissible for $n$, and thus $x\in X$. Finally, by~\Cref{lem:subshift} $X$ is space-filling.
  \end{proof}

\subsection{Two-sided space-filling actions}

A reader trying to introduce the notion of a space-filling action would perhaps not take the same route as us. A natural assumption is to ask for the maps $R_x$ to be invertible. In this subsection, we study whether asking for this property changes anything of significance.

Let $G$ be a discrete group and $X$ be a compact metrizable space endowed with an action $G\curvearrowright X$ by homeomorphisms. A continuous map $R\colon X \to G$ is called \define{$G$-invertible} if there is a continuous map $L\colon X \to G$ such that for all $x \in X$ we have that $R_x^{-1} = L_x$.

In particular, a $G$-invertible map $R\colon X \to G$ induces, for every $x\in X$, a $\ZZ$-action on $G$ given by $n\cdot g = R^n_x(g)$.

\begin{defn}
    An action $G\act X$ is \define{two-sided space-filling} if it is space-filling through a $G$-invertible map $R\colon X \to G$.
\end{defn}

\begin{lemma}
\label{lem:transitivo}
    Consider an infinite group $G$, and let $G\act X$ be a space-filling action with movement function $R$. For all $x\in X$ there exists $z\in\overline{\orb(x)}$ such that $R_z$ is invertible.
\end{lemma}
\begin{proof}
    By~\Cref{prop:spacefilling_isFG}, $G$ is finitely generated. Consider a word metric on $G$ such that $R(X)\subseteq B(1)$, and let $f$ be a covering function for $R$. Take $x \in X$. For each $n \in \NN$, by~\Cref{rem:covering_function} there is $g_n \in G$ such that $g_nB(n) \subset R_x^{\llbracket1,f(n)+1\rrbracket}(1_G)$. By compactness of $X$, if we consider the sequence $(x_n)_{n\in \NN}$ given by $x_n = g_n^{-1}x$, we may extract a limit point $z$.

    We want to show that for every $g\in G$, there is a unique $h\in R(X)$ such that $R(hg^{-1}z) = h$. Let $g \in G$ be arbitrary and take $m > |g|+1$ such that $R(\gamma z) = R(\gamma x_m)$ for all $\gamma\in B(|g|+1)$, which exists as $z$ is a limit point of $(x_n)_{n\in\NN}$. By definition, there is $t \in \llbracket 1,f(m)+1\rrbracket$ such that $R^t_x(1_G) = g_mg$. Hence if we let $h = (R^{t-1}_x(1_G))^{-1} g_m g \in R(X)$, we get that 
    \begin{align*}
        h^{-1}R(hg^{-1}x_m) & = (g_m g)^{-1}g_mgh^{-1}R(hg^{-1}g_m^{-1}x)\\ & = (g_mg)^{-1} R_x(g_mgh^{-1})\\ & = (g_mg)^{-1} R_x(R_x^{t-1}(1_G))\\& = 1_G.
    \end{align*}

  From where we get that $R(hg^{-1}z) = R(hg^{-1}x_m)=h$. Now take $h'\in R(X)$ such that $R(h'g^{-1}z)=h'$ as well. Take $t'\in \llbracket 1,f(m)+1\rrbracket$ such that $R^{t'}_x(1_G) = g_mg(h')^{-1}$, which exists as $g(h')^{-1}\in B(m)$, similarly as above, we deduce that $R^{t'+1}_x(1_G) = g_mg$. If $t \neq t'+1$, the sequence $R^i_x(1_G)$ goes through the element $g_mg$ twice and thus it is ultimately periodic, that is, $\{R^{k}_x(1_G) \mid k \in \NN\}$ is finite, which contradicts the fact that $G\act X$ is space-filling. We deduce that $t=t'+1$, so $g_mgh^{-1} = R_x^{t-1}(1_G) = R^{t'}_x(1_G) = g_mgh'^{-1}$, and thus $h = h'$. 

  Let $L_z\colon G \to G$ be given by $L_z(g)=gh^{-1}$ where $h \in R(X)$ is the only element such that $R(hg^{-1}z) = h$. By definition, $u\in G$ satisfies $R_z(u) = g$ if and only if $u = gh^{-1}$ for some $h\in R(X)$ with $R(hg^{-1}z) = h$. Hence the above shows that every $g\in G$ has exactly one preimage under $R_z$, namely $L_z(g)$, and thus $R_z$ is invertible.
\end{proof}

\begin{prop}
\label{prop:jaleadoble}
For every nonempty space-filling action $G\curvearrowright X$ there exists a closed and $G$-invariant subset $X_{\leftrightarrow}\subset X$ such that $X_{\leftrightarrow}\cap \overline{\orb(x)}\neq \varnothing$ for every $x\in X$, and such that $G\curvearrowright X_{\leftrightarrow}$ is a two-sided space-filling action. Furthermore, if $X$ is a $G$-SFT, then so is $X_{\leftrightarrow}$.
\end{prop}

\begin{proof}
    If $G$ is finite, then for every $x \in X$ every forward orbit of $R_x$ is all of $G$, so $R_x$ is a bijection and we may take $X_{\leftrightarrow}=X$. Suppose that $G$ is infinite. By~\Cref{prop:spacefilling_isFG}, $G$ is finitely generated.

    Let $R\colon X \to G$ be a movement function. Consider 
    \begin{align*}
        X_{\leftrightarrow} &= \{x \in X : R_x \mbox{ is invertible}\}\\ 
        & = \{ x \in X : \forall g \in G, \exists ! h \in R(X), R(hg^{-1}x)=h\}.
    \end{align*}
    By~\Cref{lem:transitivo} we have that $X_{\leftrightarrow}\cap \overline{\orb(x)}$ is nonempty for every $x\in X$. Also notice that $X_{\leftrightarrow}$ is closed and $G$-invariant, and thus $G\curvearrowright X_{\leftrightarrow}$ is a well-defined $G$-action. Finally, notice that if we let $L\colon X_{\leftrightarrow} \to G$ be the map given for $x\in X_{\leftrightarrow}$ by $L(x)=h^{-1}$ where $h\in G$ is the unique element such that $R(hx)=h$, then $L$ is continuous and $R_x^{-1}=L_x$. By~\Cref{rem:detalpalotalastilla}, $G\act X_{\leftrightarrow}$ is space-filling. 

    In the case where $X$ is an SFT, we can obtain $X_{\leftrightarrow}$ as a sub-SFT by forbidding all patterns $p$ with support $B(2r)$, where $r$ is such that $x|_{B(r)} = y|_{B(r)}$ implies $R(x) = R(y)$ and $R(X)\subseteq B(r)$, which do not satisfy the condition on the second line of the equation for $g=1_G$.   \end{proof}

The conclusion of the last proposition is that we lose nothing by assuming that our space-filling actions are two-sided. Using this, let us move on to see what happens if we ask for the action induced by the cocycle to be transitive.

\begin{lemma}
\label{lem:gamma}
    Let $G\act X$ be a two-sided space-filling action. Then, there exists $\gamma\colon\N\to\N$ such that for every $x\in X$, and $n\in\N$, there exists $g\in G$ such that for $z = g\cdot x$ we have
    \[B(i)\subset R_z^{\llbracket-\gamma(i), \gamma(i)\rrbracket}(1_G),\] 
    for all $i\leq n$.
\end{lemma}

\begin{proof}
    We may assume $G$ is infinite. Furthermore, choose a word metric on $G$ such that $R(X)\subseteq B(1)$. Let $f$ be a covering function for $R$. Note that $|R^k_x(1_G)|\leq k$ for every $k \in \NN$, and that $f(m)\geq m$ for every $m \in \NN$, as a ball of radius $m$ in an infinite group has at least $m+1$ elements. We define $\gamma$ inductively: $\gamma(0)=0$, $\gamma(1) = f(1) +1$ and $\gamma(n) = 2f(\gamma(n-1)+1)$. We shall prove by induction that for every $x \in X$, there exists $g_n \in G$ with $|g_n|\leq \gamma(n) - n$ such that 
    \[ g_nB(i)\subseteq R^{\llbracket -\gamma(i),\gamma(i)\rrbracket}_x(g_n),\]
for every $i\in\{1,\dots,n\}$.  By~\Cref{prop:faciles}, the previous equation is equivalent to $B(i)\subseteq R^{\llbracket -\gamma(i),\gamma(i)\rrbracket}_{z}(1_G)$ for $z=g_n^{-1}\cdot x$. This proves the lemma with $g=g_n^{-1}$.

    For the base case, take $g_1$ such that $g_1B(1)\subseteq R_x^{\llbracket 0,f(1)\rrbracket}(1_G)$. Then $g_1 = R^k_x(1_G)$ for some $k\leq f(1)$, so $|g_1|\leq f(1) = \gamma(1)-1$, and the inclusion follows from~\Cref{prop:faciles}. Now, let $n \geq 2$ and suppose the hypothesis is true for $n-1$. Because $f$ is a covering function, there exists $h_n\in G$ such that \[h_nB(\gamma(n-1)+1)\subseteq R_x^{\llbracket-f(\gamma(n-1)+1), f(\gamma(n-1)+1)\rrbracket}(1_G).\] Then, by our induction hypothesis applied to $h_n^{-1}\cdot x$, there exists $g_{n-1}\in G$ with $|g_{n-1}|\leq \gamma(n-1)-n+1$ such that $g_{n-1}B(i)\subseteq R^{\llbracket -\gamma(i),\gamma(i)\rrbracket}_{h_n^{-1}x}(g_{n-1})$ for every $i\in\{1,\dots,n-1\}$. Take $g_n = h_ng_{n-1}$ and notice that $g_n = R^k_x(1_G)$ for some $k\leq f(\gamma(n-1)+1)=\gamma(n)/2$, so $|g_n|\leq\gamma(n)/2\leq\gamma(n)-n$, because $\gamma(n)/2\geq \gamma(n-1)+1\geq n$. Furthermore, if we consider $z = g_{n}^{-1}\cdot x$, we have for all $i \leq n-1$ that \[   B(i)\subset R_{z}^{\llbracket-\gamma(i), \gamma(i)\rrbracket}(1_G).\]

     It suffices to check that this condition also holds for $i = n$. Indeed, we have that 
    \begin{align*}
        h_ng_{n-1}B(n)& \subset h_n B(n+|g_{n-1}|) \\ & \subset h_n B(n+\gamma(n-1)-(n-1))\\ &= h_n B(\gamma(n-1)+1)\\ & \subset R_x^{\llbracket-f(\gamma(n-1)+1), f(\gamma(n-1)+1)\rrbracket}(1_G)\\
        &\subset R_x^{\llbracket-\gamma(n),\gamma(n)\rrbracket}(h_ng_{n-1}).
    \end{align*}

    From where we deduce that $B(n) \subset R_{z}^{\llbracket-\gamma(n), \gamma(n)\rrbracket}(1_G)$.
\end{proof}

    \begin{prop}\label{prop:puntito_transitivo}
        Let $G\act X$ be a space-filling action with movement function $R\colon X\to G$. For every $x \in X$, there exists $y \in \overline{\orb(x)}$ such that $R_{y}$ is invertible and induces a transitive $\ZZ$-action on $G$. We call such a $y$ a \define{transitive point}.
    \end{prop}

    \begin{proof}
        Let $x\in X$. By~\Cref{prop:jaleadoble}, there exists a closed $G$-invariant set $X_{\leftrightarrow}\subset X$ such that $G\act X_{\leftrightarrow}$ is two-sided space-filling and such that there is $z \in X_{\leftrightarrow} \cap \overline{\orb(x)}$.

        By~\Cref{lem:gamma}, there exists a function $\gamma\colon\N\to\N$ such that for each $n \in \NN$, there is $g_n \in G$ such that if we let $y_n = g_n\cdot z$ then for $k \leq n$ we have $B(k) \subset R^{\llbracket-\gamma(k),\gamma(k)\rrbracket}_{y_n}(1_G)$.
        
        Consider an accumulation point $y$ of the sequence $(y_n)_{n\in \NN}$. By definition we have $y \in \overline{\orb(x)}$, and $y\in X_{\leftrightarrow}$ as this set is closed, so $R_y$ is invertible. We claim $R_y^{\Z}(1_G) = G$. Indeed, consider $g\in G$ and $n = |g|$. Then, for infinitely many $k\geq n$ we have that \[  R^{\llbracket-\gamma(n),\gamma(n)\rrbracket}_{y_k}(1_G) = R^{\llbracket-\gamma(n),\gamma(n)\rrbracket}_{y}(1_G).  \]
        In particular, we have $g \in R^{\llbracket-\gamma(n),\gamma(n)\rrbracket}_{y_k}(1_G) \subset R^{\ZZ}_y(1_G)$.\end{proof}

\subsection{Group theoretical restrictions on space-filling actions}

Let us now look at what necessary conditions we obtain for groups admitting space-filling actions. The first of these comes as a consequence of~\Cref{prop:puntito_transitivo}. The existence of a transitive $\Z$-action by bounded displacements is equivalent to the group having a Cayley graph which admits a bi-infinite Hamiltonian path.

\begin{prop}
\label{prop:inf_ends}
    If $G$ has infinitely many ends, it does not admit a space-filling action.
\end{prop}

\begin{proof}

    Suppose $G$ admits a space-filling action $G\curvearrowright X$ and let $R$ be the associated movement function. Fix some symmetric set of generators $S$ on $G$ such that $R(X)\subset S$, and let $K\Subset G$ such that the Cayley graph $\operatorname{Cay}(G\setminus K,S)$ has at least three distinct infinite connected components.
    
    By~\Cref{prop:puntito_transitivo} it follows that there exists $x\in X$ such that $R_x$ induces a transitive $\ZZ$-action. Therefore, there exists $N \in \NN$ such that $K\subset R_x^{\llbracket-N,N\rrbracket}(1_G)$. As $R_x$ is a bijection, $R^k_x(1_G)\notin K$ whenever $|k|>N$. As consecutive elements of the sequence $(R^k_x(1_G))_{k\in\ZZ}$ are equal or adjacent in the Cayley graph, all the elements with $k>N$ lie in a single connected component of $G\setminus K$, and so do all the elements with $k<-N$. Hence some infinite connected component is visited at most $2N+1$ times, which contradicts the fact that $R_x$ induces a transitive action. \end{proof}

The next result shows that, unless $G$ is very boring (virtually cyclic), every space-filling action must contain some points $x$ for which the map $R_x$ is noninvertible.

\begin{prop}\label{prop:transitive_implies_virt_cyclic}
    Let $G$ be a group and suppose there exists a space-filling action $G\curvearrowright X$ such that for every $x\in X$ the map $R_x$ is invertible and induces a transitive $\ZZ$-action. Then, $G$ is virtually cyclic.
\end{prop}

\begin{proof}
    By~\Cref{prop:spacefilling_isFG} $G$ is finitely generated, thus we fix a word metric on $G$. We claim that there exists $\gamma\colon \NN \to \NN$ such that for all $x \in X$ \[ B(n)\subset R_x^{\llbracket-\gamma(n),\gamma(n)\rrbracket}(1_G).    \]

    Suppose by contradiction this does not hold, that is, there exists $n \in \NN$ such that for each $m\in \NN$ there is $x_m\in X$ and $t\in B(n)$ such that $t\notin R^{\llbracket -m, m\rrbracket}_{x_m}(1_G)$. 

    As $B(n)$ is finite, there is $t_0\in B(n)$ which satisfies the above property for infinitely many $m$. Let $(m_k)_{k \in \NN}$ be those indices. By compactness of $X$, we may extract a limit point $\bar{x}$ of the sequence $(x_{m_k})_{k \in \NN}$. We claim that the $\ZZ$-action induced by $R_{\bar{x}}$ is not transitive. Indeed, if $1_G$ and $t_0$ were in the same orbit, there would exist $m \in \ZZ$ such that $t_0 = R_{\bar{x}}^m(1_G)$. But there exists $k\in \NN$ such that $m_k\geq |m|$ and $R^{\llbracket -|m|, |m|\rrbracket}_{x_{m_k}}(1_G) = R^{\llbracket-|m|,|m|\rrbracket}_{\bar{x}}(1_G)$ which leads to a contradiction.

    Fix $\gamma\colon \NN \to \NN$ the smallest function which satisfies that $B(n)\subset R^{\llbracket-\gamma(n),\gamma(n)\rrbracket}_x(1_G)$ for all $x \in X$ and $n\in\NN$. Note that for any $n \geq 1$ and $h \in B(n+1)$ there is $g \in B(n)$ such that $h \in gB(1)$. By definition of $\gamma$, we have that there is $k\in \llbracket-\gamma(n),\gamma(n)\rrbracket$ such that $R_x^k(1_G) = g$, and, as $g^{-1}h\in B(1)$, there is $k'\in \llbracket-\gamma(1),\gamma(1)\rrbracket$ such that $R^{k'}_{g^{-1}x}(1_G) = g^{-1}h$, that is, $R^{k'}_{x}(g)= gR^{k'}_{g^{-1}x}(1_G) = h$. From where we deduce that $R^{k+k'}_x(1_G)=h$. As $h$ is arbitrary, we obtain that $\gamma(n+1)\leq \gamma(n)+\gamma(1)$ and thus $\gamma(n) \leq n\gamma(1)$ for all $n \in \NN$.

    On the other hand, it is clear that we have \[ |B(n)| \leq 2\gamma(n)+1 \leq 2n \gamma(1)+1.   \]

    From where we deduce that $|B(n)|$ grows at most linearly and thus $G$ is virtually cyclic.
    \end{proof}

\begin{prop}\label{prop:transitive_points_are_free}
    Unless $G$ is virtually cyclic, every transitive point in a space-filling action has trivial stabilizer.
\end{prop}

\begin{proof}
    Let $y$ be a transitive point and suppose $g\cdot y = y$ for some $g\in G\setminus \{1_G\}$. As $y$ is transitive and $g\neq 1_G$, there exists $n\in\ZZ\setminus\{0\}$ such that $R^n_y(1_G)= g$. Note that we have $R_y^{n}(g) = gR^{n}_{g^{-1}y}(1_G) = gR^n_{y}(1_G) = g^2$. An induction argument shows that for $\ell \in \ZZ$ we have $R^{\ell n}_y(1_G) = g^{\ell}$. It follows that if we let $m \in \ZZ$ and we let $m = \ell n + r$ for $r \in \llbracket 0,|n|-1\rrbracket$, then \[ R^m_y(1_G) = R^r_y(R^{\ell n}_y(1_G)) = R^r_y(g^{\ell}) = g^{\ell}R^r_{g^{-\ell}\cdot y}(1_G) = g^{\ell}R^r_{y}(1_G).   \]
    From where it follows that \[ G = R^{\ZZ}_y(1_G) \subset \langle g \rangle R_y^{\llbracket 0,|n|-1\rrbracket}(1_G).  \]
    Hence, $G$ is virtually cyclic.\end{proof}

    In what follows, we write $\operatorname{FC}(G)=\{g \in G : |\{h^{-1}gh : h \in G\}|<\infty\}$ for the subset of group elements with finite conjugacy class.

    \begin{corollary}
        Suppose $G$ is not virtually cyclic, let $G\curvearrowright X$ be a space-filling action and let $x \in X$. The following statements hold:
        \begin{enumerate}
            \item $\operatorname{stab}(x)\cap \operatorname{FC}(G) = \{1_G\}$. In particular, if $\operatorname{FC}(G) = G$, the action must be free.
            \item $\operatorname{orb}(x)$ must be infinite. 
        \end{enumerate}

    \end{corollary}

    \begin{proof}
        By~\Cref{prop:puntito_transitivo} there is a sequence $(g_n)_{n \geq 1}$ and a transitive point $y \in X$ such that $g_n\cdot x$ converges to $y$.
        For the first statement, suppose there is $h\in G\setminus \{1_G\}$ in $\operatorname{stab}(x)\cap \operatorname{FC}(G)$. As $g_nhg_{n}^{-1}$ fixes $g_nx$ and $h$ has finitely many conjugates, it follows that one of them fixes infinitely many elements in the sequence. Passing to a subsequence shows that $y$ does not have a trivial stabilizer, which contradicts~\Cref{prop:transitive_points_are_free}.

        For the second statement, suppose $x$ has finite orbit, then $H=\operatorname{stab}(x)$ is a finite index subgroup of $G$, hence there is $N\lhd G$ of finite index with $N\subset H$. As $N$ is normal, it follows that $N$ fixes $g_n \cdot x$ for all $n$, and thus fixes $y$. As $G$ is infinite, $N$ is nontrivial, which contradicts~\Cref{prop:transitive_points_are_free}.
    \end{proof}

In our quest to find space-filling SFTs, this last corollary provides additional constraints. Item 2 implies that, at the very least, the subshift must be weakly aperiodic, that is, it has no finite orbits. It is an open problem to determine which groups admit weakly aperiodic SFTs. Carroll and Penland conjecture~\cite{carroll2015periodic} that a finitely generated group admits a weakly aperiodic SFT if and only if the group is not virtually cyclic. Although the conjecture remains open, many classes of groups have been shown to satisfy it (see~\cite{bitar2024realizability} for a recent survey on this problem).

\subsection{Stability under direct products}

\begin{defn}
    Let $\Gamma,G$ be groups, $\phi\colon \Gamma \to G$ an epimorphism and $X\subset A^G$ a subshift. The \define{pullback} subshift $\widehat{X}\subset A^{\Gamma}$ is given by \[ \widehat{X} = \{ \widehat{x}\in A^{\Gamma} : \exists x \in X \mbox{ such that } \forall g \in \Gamma, \widehat{x}(g) = x(\phi(g))\}.    \]
\end{defn}

It is straightforward to prove that if a subshift $X$ is an SFT and $\ker(\phi)$ is finitely generated, then its pullback $\widehat{X}$ is also an SFT.

\begin{prop}
\label{prop:producto}

    Consider two infinite finitely generated groups $G$ and $H$, and $\Gamma = G\times H$. If $G$ and $H$ admit space-filling subshifts, such that their pullbacks to $\Gamma$ are sofic, then $\Gamma$ admits a space-filling SFT.
\end{prop}

\begin{proof}
    The idea is to use the directions given by each configuration of the pullbacks of the two space-filling subshifts to define a grid, and then use Kari's $\Z^2$-space-filling SFT over the grid.
    
    Let $X_G$ and $X_H$ be space-filling subshifts of $G$ and $H$ whose pullbacks to $\Gamma$ are sofic. By~\Cref{prop:jaleadoble} and its proof, we may replace them by the subshifts obtained by additionally imposing the local condition which makes their movement functions invertible. The pullbacks of these conditions are local conditions on $\Gamma$, so the new pullbacks are still sofic. Denote their movement functions $R^G$ and $R^H$. Consider the Kari space-filling $\Z^2$-SFT $Z\subseteq A^{\Z^2}$ defined by nearest neighbor horizontal and vertical forbidden patterns, $\FF_{h}\cup \FF_v$, with movement function $R^Z\colon Z\to \{(\pm1,0),(0,\pm 1)\}$ which for $z\in Z$ depends only on $z\bigl( (0,0)\bigr)$. Our goal is to define $Y\subseteq\widehat{X}_G\times\widehat{X}_H\times A^{\Gamma}$, which will be a sofic space-filling subshift for $\Gamma$.

    The forbidden patterns for $Y$ are as follows: for every $(g,h)\in\Gamma$ and $(x,y,z)\in \widehat{X}_G\times\widehat{X}_H\times A^{\Gamma}$, if $R^G(g^{-1}x) = s$ and $R^H(h^{-1}y) = t$, then \[\bigl(z(g,h), z(gs,h)\bigr)\notin\FF_h, \ \bigl(z(g,h), z(g,ht)\bigr)\notin\FF_v.\] As both pullbacks are sofic and the remaining conditions are local, the subshift $Y$ is sofic. We shall show that $Y$ is nonempty and space-filling, which by~\Cref{cor:conjugacy_invariant} implies that $\Gamma$ admits a space-filling SFT. Next, consider the continuous function $R\colon Y\to\Gamma$ given by
    \[R(x,y,z) = \begin{cases}\bigl((R^G_x)^{\pm1}(1_G),1_H\bigr) & \textnormal{if } R^Z(z) = (\pm1,0),\\ \bigl(1_G,(R^H_y)^{\pm1}(1_H)\bigr)  & \textnormal{if }R^Z(z)= (0,\pm1).\end{cases}\]
    Due to the definition of $Y$, for every $ \omega = (x,y,z)\in Y$ we can define an injective function $p_\omega\colon \Z^2 \to \Gamma$ by $p_\omega(i,j) = ((R^G_x)^i(1_G), (R^H_y)^j(1_H))$, such that $z' = z\circ p_{\omega}\colon \Z^2\to A$ is a configuration from $Z$.
    
    Now, consider $\omega = (x,y,z)\in Y$ and $n\in\N$. Take the functions $f_G,f_H, f_Z\colon \N\to\N$ given by~\Cref{rem:covering_function} for the subshifts $X_G$, $X_H$, and $Z$ respectively. Define $k = \max\{f_G(n), f_H(n)\}+1$ and $M = f_Z(k)$. Then, there exists $(v_1, v_2)\in\Z^2$ such that 
    \[(v_1, v_2) + \llbracket 0, k-1\rrbracket^2 \subseteq (R^Z_{z'})^{\llbracket 0, M\rrbracket}\bigl( (0,0)\bigr).\] Let $F = p_\omega((v_1, v_2) + \llbracket 0, k-1\rrbracket^2)$. Then, $F\subseteq R^{\llbracket 0, M\rrbracket}_\omega(1_\Gamma)$. By~\Cref{rem:covering_function}, there exist $g\in G$ and $h\in H$ such that $gB_G(n)\times hB_H(n)\subseteq F$. As $B_G(n)\times B_H(n)$ contains the ball of radius $n$ of $\Gamma$ for the generating set $(S_G\times\{1_H\})\cup(\{1_G\}\times S_H)$, this implies that $Y$ is space-filling.

    It remains to see that $Y$ is nonempty. Consider the configurations $x'\in X_G$ and $y'\in X_H$ and $z'\in Z$ given by~\Cref{prop:puntito_transitivo}. Additionally, take $x = x'\circ\pi_G\in\widehat{X}_G$ and $y = y'\circ\pi_H\in\widehat{X}_H$ where $\pi_G$ and $\pi_H$ are the projections from $\Gamma$ onto $G$ and $H$ respectively. Now, for every $(g,h)\in\Gamma$ there exists unique $i,j\in\Z$ such that $(R^G_{x})^i(1_G) = g$ and $(R^H_{y})^j(1_H) = h$. Then, we define the configuration $z\in A^{\Gamma}$ by $z(g,h) = z'(i,j)$. We claim $(x,y,z)$ is in $Y$. Indeed, if we denote $s = R^G(g^{-1}x)$ and $t =R^H(h^{-1}y)$, then $gs=(R_x^G)^{i+1}(1_G)$ and $ht = (R^{H}_y)^{j+1}(1_H)$. This implies \begin{align*}
        \bigl(z(g,h), z(gs,h)\bigr) & = \bigl(z'(i,j),z'(i+1,j)\bigr) \notin\FF_{h},\\
        \bigl(z(g,h), z(g,ht)\bigr) & = \bigl(z'(i,j),z'(i,j+1)\bigr) \notin\FF_{v}.
    \end{align*}
    Therefore $(x,y,z)\in Y$, and thus $Y$ is a nonempty sofic space-filling subshift on $\Gamma$.\end{proof}

    We obtain the following direct consequence.

\begin{corollary}\label{cor:productos}
      Let $G$ and $H$ be infinite finitely generated groups. If $G$ and $H$ are space-filling (respectively computably space-filling, finitary space-filling), then so is $G\times H$.
  \end{corollary}

\begin{proof}
    The proof of~\Cref{prop:producto} shows that the subshift $Y$ is nonempty and space-filling without using that the pullbacks are sofic. If $X_G$ and $X_H$ are effective, then so are their pullbacks, and thus $Y$ is effective. If $X_G$ and $X_H$ are SFTs, then so are their pullbacks, as the kernels of the projections are finitely generated, and thus $Y$ is an SFT.
\end{proof}

As a direct consequence of the previous result and Kari's construction in $\Z^2$, we obtain the following.

\begin{corollary}\label{cor:Zd_is_space_filling_locally}
    For every $d\geq 1$, $\Z^d$ is finitary space-filling.
\end{corollary}

This corollary allows us to extend the previous result to an arbitrary number of groups, using a space-filling $\Z^n$-SFT instead of Kari's. The proof is analogous.
\begin{theorem}
\label{thm:producto}
    Let $n \geq 1$ and consider infinite and finitely generated groups $(G_{i})_{i=1}^{n}$ and let $\Gamma = \prod_{i=1}^n G_i$ be their direct product. If every $G_i$ admits a space-filling subshift such that its pullback to $\Gamma$ is sofic, then $\Gamma$ admits a space-filling SFT.
\end{theorem}

\subsection{Computable space-filling}

In this subsection, we will study a computable variant of space-filling actions. We will at first define this notion for arbitrary finitely generated and recursively presented groups, but after an initial result, we will restrict to groups with decidable word problem. We shall assume the reader is familiar with elementary notions of computability, such as the notion of decidable and recursively enumerable set. We refer the reader to~\cite{soare2016turing,sipser1996introduction} for an introduction to computability theory.

Let $G$ be a finitely generated and recursively presented group and fix a finite symmetric set of generators $S$ and a finite alphabet $A$. A \define{pattern coding} is a finite collection $c\subset S^*\times A$. We say that a subshift $X\subset A^G$ is \define{effective} (by patterns) if $X=A^G\setminus\bigcup_{g\in G,c\in \mathcal{C}}g[c]$ for a recursively enumerable set $\mathcal{C}$ of pattern codings, where $[c] = \{x \in A^G : x(\underline{w})=a \mbox{ for every } (w,a)\in c\}$ and $\underline{w}$ denotes the element of $G$ obtained by multiplying the generators of $w$.

We note that every SFT is effective by patterns. Furthermore, for recursively presented groups, a subshift is effective by patterns if and only if it is topologically conjugate to a computable action on a $\Pi_1^0$ set (see~\cite{BarCarRoj2024}) and thus we just call it \define{effective}. In particular, being an effective subshift is an invariant of topological conjugacy.

\begin{defn}\label{def:comp_SF}
    A finitely generated and recursively presented group $G$ is \define{computably space-filling} if it admits a nonempty effective space-filling subshift.
\end{defn}

The main motivation behind our interest in the computable variant of space-filling is that by using simulation theorems from the literature such as~\cite{BarSab2019_semidirect,Barbieri2019,BarSabSal2026,BarSabSal_free_2025} we may turn our effective subshift into a sofic one in a possibly larger group, and thus use this to give examples of groups which are finitary space-filling. We note that here we give the definition only in recursively presented groups for convenience (to avoid going into the subtleties of effective by pattern subshifts when the word problem is not recursively enumerable), but in the future exploring that notion beyond this class might be interesting.

First, we will show that groups which are both recursively presented and computably space-filling  necessarily have decidable word problem, and thus we may restrict our attention to that setting. This proof is very similar to the argument of Jeandel~\cite{Jeandel2015Subshifts} that shows that recursively presented groups which admit strongly aperiodic effective subshifts necessarily have decidable word problem.

\begin{prop}\label{prop:rp_decidable}
    Let $G$ be a finitely generated and recursively presented group. If $G$ is computably space-filling, then $G$ has decidable word problem.
\end{prop}

\begin{proof}
    The case of a finite group $G$ is immediate, as they all have decidable word problem. Suppose that $G$ is infinite and fix a finite set of generators $S$. As $G$ is recursively presented, in order to show that it has decidable word problem it suffices to find an algorithm that given $w\in S^*$ halts precisely when $\underline{w}\neq 1_G$.
    
    Let $X\subset A^G$ be a nonempty effective space-filling subshift with movement function $R$ and covering function $f$. Without loss of generality (by taking a topologically conjugate version of $X$) we may assume that $R(x)$ depends only on the value $x(1_G)$. Thus, for every $a \in A$ we may fix a word $r_a\in S^*$ such that $R(x) = \underline{r_{x(1_G)}}$ for every $x\in X$. Let $\ell = \max_{a\in A}|r_a| \geq 1$. For $j \geq 0 $ and $x\in X$, write $g_j = R^j_x(1_G)$ and notice that as $G$ is infinite, the map $j\mapsto g_j$ is injective.

    Let $w\in S^*$. By definition, there is $h \in G$ such that $hB_S(|w|)\subset\{g_0,\dots,g_{f(|w|)}\}$, so there exist $0 \leq i,j \leq f(|w|)$ such that $h = g_i$ and $h\underline{w}=g_j$. Therefore:
    \begin{itemize}
        \item if $\underline{w}\neq 1_G$, then for every $x\in X$ there are $i\neq j$ in $\llbracket 0,f(|w|)\rrbracket$ such that $g_i\underline{w} = g_j$;
        \item if $\underline{w} = 1_G$, then for no $x\in X$ there are $i\neq j$ such that $g_i\underline{w}=g_j$, by injectivity.
    \end{itemize}

    For $n \in \NN$ write $S^{\leq n} = \bigcup_{k\leq n}S^k$. Consider the algorithm which on input $w\in S^*$ does the following for $n=1,2,\dots$. First, it runs for $n$ steps both the algorithm that enumerates the words which represent $1_G$ and the algorithm that enumerates the forbidden pattern codings of $X$. Then it constructs the set $U_n$ of all maps $\xi \colon S^{\leq n}\to A$ such that $\xi(u)=\xi(u')$ whenever $u^{-1}u'$ has been detected to represent the identity, and such that there are no enumerated pattern coding $c$ and $u\in S^*$ with $\xi(uv)=a$ for every $(v,a)\in c$. For every $\xi\in U_n$ it constructs the words $v_0,\dots,v_{\lfloor n/\ell\rfloor}$ given by $v_0 = \varepsilon$ and $v_i = v_{i-1}r_{\xi(v_{i-1})}$, which is possible as $|v_{i-1}|\leq n$, and checks whether there are $i\neq j$ such that $v_iw v_j^{-1}$ has been enumerated. The algorithm halts if this is the case for every $\xi\in U_n$.

    Suppose first that $\underline{w}=1_G$ and take $x\in X$. For every $n$ the map $\xi(u) = x(\underline{u})$ belongs to $U_n$, and its associated words satisfy $\underline{v_i} = g_i$. By the second item above, there are no $i\neq j$ with $\underline{v_i}\,\underline{w} = \underline{v_j}$, and thus the algorithm does not halt.

    Suppose now that $\underline{w}\neq 1_G$ and let $N = \ell f(|w|)$. A map $\widehat{x}\colon S^*\to A$ whose restriction to $S^{\leq n}$ belongs to $U_n$ for every $n$ respects every relation of $G$ and contains no forbidden pattern coding, so it is of the form $\widehat{x}(u) = x(\underline{u})$ for some $x\in X$. Hence, by compactness, there is $n_0\geq N$ such that for every $n\geq n_0$ and $\xi\in U_n$ there is $x\in X$ with $\xi(u) = x(\underline{u})$ for every $u\in S^{\leq N}$. For such $\xi$ and $x$ we have $\underline{v_i} = g_i$ for every $i\leq f(|w|)$, and by the first item above there are $i\neq j$ in $\llbracket 0,f(|w|)\rrbracket$ with $\underline{v_iwv_j^{-1}} = 1_G$. As there are finitely many maps $S^{\leq N}\to A$, there is $n\geq n_0$ by which all the corresponding words $v_iwv_j^{-1}$ have been enumerated, and the algorithm halts at stage $n$ at the latest.
\end{proof}

We obtain the following natural consequence.

\begin{corollary}
      Let $G$ be a finitely generated and recursively presented group. If $G$ is finitary space-filling, then $G$ has decidable word problem.
  \end{corollary}

Indeed, as we mentioned above, every SFT is effective. In particular, if $G$ is a finitely presented group with undecidable word problem, then $G\times\Z$ is a one-ended finitely presented group which is not finitary space-filling, and not computably space-filling. We do not know whether it is space-filling.

The arguments above justify looking exclusively at groups with decidable word problem. From this point forward, we assume that $G$ is infinite, finitely generated and that it has decidable word problem. 

Fix a finite symmetric set of generators $S$ for $G$. We represent elements of $G$ by words over $S$. As the word problem is decidable, for every $n\in\NN$ we can compute a list of words representing each element of $B_S(n)$ exactly once, and given finite sets $F,F'\Subset G$ described by lists of words, we can decide whether $F\subset F'$ and compute the translates $gF$ for $g \in G$. A pattern $p\colon F\to A$ is described by the finite list of pairs $(u,p(u))$, where $u$ ranges over words representing the elements of $F$. In this setting a subshift $X\subset A^G$ is effective if and only if $X = \X_{\FF}$ for a recursively enumerable set of patterns $\FF$.

Next we will show that groups which are computably space-filling are precisely those that admit a space-filling action whose covering function is computable. We will need a folklore lemma, whose proof is a standard compactness argument and can be found in~\cite[Lemma 4.3]{Bar2021_entropy_fails}.

Fix an alphabet $A$ and let $\FF$ be a set of forbidden patterns on alphabet $A$. For $m \geq 0$ denote by $\mathcal{L}^{\mathrm{loc}}_{m}(\FF)$ the set of patterns $p \in A^{B(m)}$ which are \define{locally admissible} with respect to $\FF$, that is, those that do not contain any translate of a pattern from $\FF$ (for every $q\in \FF$ and $g\in G$, we have that $[p]$ is not contained in $g[q]$).

\begin{lemma}\label{lem:locally_admissible}
    Let $\FF$ be as above. For every $k \geq 0$ there is $m \geq 0$ such that the restriction of every $p \in \mathcal{L}^{\mathrm{loc}}_{m}(\FF)$ to $B(k)$ is \define{globally admissible} with respect to $\FF$, that is, $[p|_{B(k)}]\cap \X_{\FF}\neq \varnothing$.
\end{lemma}

\begin{prop}\label{prop:computable_cov_function}
    A group $G$ is computably space-filling if and only if it admits a space-filling action with a computable covering function.
\end{prop}

\begin{proof}
    Suppose $G$ admits a space-filling action with a computable covering function $f$ with respect to some finite symmetric set of generators which contains the image of the movement function. The subshift $\mathfrak{X}(1,f)$ of~\Cref{thm:todo_es_subshift} is nonempty and space-filling by~\Cref{lem:subshift}, and effective because $f$ is computable and $G$ has decidable word problem.

    Conversely, suppose $G$ admits an effective space-filling subshift $X\subset A^G$ and let $f\colon \NN \to \NN$ be a covering function. We will show that there exists a computable $f'\colon \NN\to \NN$ which is also a covering function. Note that the movement function $R\colon X \to G$ is computable, as it depends on a finite window $B(r)\Subset G$. Also write $\ell = \max_{x \in X}|R(x)|_S$. For a pattern $p$ with support $B(n)$ we write $R^t_p(1_G)$ for the sequence obtained by applying the local rule of $R$ to $p$, as long as it is defined.

    As $X$ is effective, there is a Turing machine which on input $n$, outputs a finite set of forbidden patterns $\FF_n$ (which generate a nested sequence) such that if we let $\FF = \bigcup_{n \geq 0}\FF_n$, then $X=\X_{\FF}$.

    Let $k \geq 1$ be an input. Iteratively on $n$ and starting on $n=1$, construct the set $\FF_n$ and the collection of patterns $\mathcal{L}^{\mathrm{loc}}_{n}(\FF_n)$. Note that it is possible because $G$ has decidable word problem. For every $p \in \mathcal{L}^{\mathrm{loc}}_{n}(\FF_n)$, compute $R^t_p(1_G)$ for all possible $t$ until it either cycles, covers a translated copy of $B(k)$ or it becomes undefined.

    If for some $p \in \mathcal{L}^{\mathrm{loc}}_{n}(\FF_n)$ the map cycles or becomes undefined, update the value of $n$ to $n+1$ and continue the algorithm. Otherwise, on every pattern $p$ the procedure covers a translated copy of $B(k)$, and we may define $f'(k)$ as the maximum value of $t$ for which this is achieved among all such $p$.

    We only need to argue that the algorithm above stops and yields a value $f'(k)$ in finitely many steps. By~\Cref{lem:locally_admissible} there is $m$ such that the restriction of every locally admissible $p \in \mathcal{L}^{\mathrm{loc}}_{m}(\FF)$ to $B(\ell f(k)+r)$ yields a globally admissible pattern. Also, as there are finitely many patterns with support $B(m)$, there is $n_0\geq m$ such that every pattern with support $B(m)$ which is not in $\mathcal{L}^{\mathrm{loc}}_{m}(\FF)$ contains a translate of a pattern of $\FF_{n_0}$. It follows that the restriction of every $q \in \mathcal{L}^{\mathrm{loc}}_{n_0}(\FF_{n_0})$ to $B(m)$ belongs to $\mathcal{L}^{\mathrm{loc}}_{m}(\FF)$, and thus its restriction to $B(\ell f(k)+r)$ yields a globally admissible pattern which occurs in $\X_{\FF}$, and thus we are guaranteed that it will cover a translate of $B(k)$ within $f(k)$ steps. Hence the algorithm stops at $n_0$ at the latest. Finally, $f'$ is a covering function: for $x\in X$ and the final value of $n$, the pattern $x|_{B(n)}$ is one of the patterns which were tested, so $R^{\llbracket 0,f'(k)\rrbracket}_x(1_G)$ contains a translate of $B(k)$. \end{proof}

\subsubsection{From effectiveness to soficity}
There are several results in the literature that allow us to go from an effective subshift to a sofic subshift. These are collectively known as \define{simulation theorems}. Their general form states that given groups $G$ and $H$, and an epimorphism $\pi\colon G\to H$, the pullback of every effective subshift on $H$ is sofic in $G$. Our goal is to use these theorems to prove certain groups are finitary space-filling, by constructing effective space-filling subshifts. These results are the following

\begin{theorem}
    Let $(G_i)_{i=1}^{3}$ be infinite and finitely generated groups with decidable word problem. The following hold.
    \begin{enumerate}
        \item \cite{Barbieri2019} For every $i\in\{1,2,3\}$, the pullback of every effective $G_i$-subshift to $G_1\times G_2\times G_3$ is sofic.
        \item \cite[Corollary 6.1 and Theorem 1.3]{BarSabSal_free_2025} If either $G_1$ or $G_2$ is non-amenable, the pullback of every effective $G_i$-subshift to $G_1\times G_2$ is sofic for $i \in \{1,2\}$.
    \end{enumerate}
\end{theorem}

Furthermore, for groups $G$ which are \define{self-simulable}, every effective $G$-subshift is sofic~\cite{BarSabSal2026}. Combining these results with our notion of space-filling, \Cref{thm:producto}, and \Cref{cor:conjugacy_invariant} we obtain the following.

\begin{corollary}
\label{cor:simulation_y_SF}
    Consider infinite, finitely generated and recursively presented groups $(G_i)_{i=1}^{3}$ which are computably space-filling. The following hold.
    \begin{enumerate}
        \item $G_1\times G_2\times G_3$ is finitary space-filling.
        \item If $G_1$ is non-amenable, then $G_1\times G_2$ is finitary space-filling.
        \item If $G_1$ is self-simulable, then $G_1$ is finitary space-filling.
    \end{enumerate}
\end{corollary}

\section{Invariance under quasi-isometries}
\label{sec:QI}

In this section we show that the property of admitting a space-filling action is a quasi-isometry invariant for finitely generated groups. Furthermore, if two quasi-isometric groups $G,H$ have decidable word problem, then $G$ is computably space-filling if and only if $H$ is, and if they are finitely presented, then $G$ is finitary space-filling if and only if $H$ is. We remark that this will not require any computability assumption on the underlying quasi-isometry.

To do this, we first provide an adequate generalization of the space-filling property in the context of graphs. Throughout this section, all graphs are simple, connected and have bounded degree. We identify a graph $\Gamma$ with its set of vertices and denote by $d=d_\Gamma$ its path metric. For $v\in \Gamma$ and $r\geq 0$ we write $B_\Gamma(v,r) = \{w \in \Gamma : d(v,w)\leq r\}$ for the ball of radius $r$ centered at $v$ and denote by $\beta_\Gamma(r) = \sup_{v\in \Gamma}|B_\Gamma(v,r)|$ the maximal cardinality of a ball of radius $r$.

Consider a graph $\Gamma$ and $k\geq 1$. A finite sequence $(w_j)_{j=0}^M$ of pairwise distinct vertices of $\Gamma$ is a \define{$k$-path} if $d(w_j,w_{j+1})\leq k$ for all $0\leq j<M$. Given $i\leq i'$ in $\llbracket 0,M\rrbracket$, we call $\{w_i,\dots,w_{i'}\}$ a \define{window} of the path.

\begin{defn}\label{def:SF_graph} Let $k\geq1$ and $f\colon\N\to\N$. We say that a graph $\Gamma$ is \define{$(k,f)$-space-filling} if for every $M\in \NN$ there is a $k$-path $(w_j)_{j=0}^M$ such that for all $r\in\N$ and all $i$ with $i+f(r)\le M$, the \define{window} $\{w_i,\dots,w_{i+f(r)}\}$ contains a ball $B_\Gamma(v,r)$ for some $v \in \Gamma$. 

We say that $\Gamma$ is \define{space-filling} if it is $(k,f)$-space-filling for some $k\geq1$ and $f\colon\N\to\N$.
\end{defn}

Notice that to verify the space-filling condition, we only have to check infinitely many $M\in\N$. Indeed, if $(w_j)_{j=0}^{M}$ is a $k$-path and $M'\leq M$, then $(w_j)_{j=0}^{M'}$ satisfies the condition of~\Cref{def:SF_graph} for $M'$, as all of its windows are windows of the longer path.
 
\begin{remark}\label{rem:SF_graph_basic}
   As was the case with our notion of space-filling actions, we can change the space-filling parameters for any upper bound. Indeed, suppose $\Gamma$ is $(k,f)$-space-filling. If $k'\geq k$ and $f'\geq f$, then $\Gamma$ is also $(k',f')$-space-filling, as every window of length $f'(r)$ contains a window of length $f(r)$. In particular, we may assume that $f$ is nondecreasing, replacing $f(r)$ by $\max_{r'\leq r}f(r')$.
\end{remark}

\subsection{Cayley graphs}

Given a finitely generated group $G$ and a finite symmetric set of generators $S$, the Cayley graph $\operatorname{Cay}(G,S)$ is the graph with vertex set $G$ in which $g$ and $gs$ are adjacent for every $g\in G$ and $s \in S\setminus\{1_G\}$. Its path metric is $d(g,h)=|g^{-1}h|_S$, its balls are $B(g,r) = gB_S(r)$ and left translations act on it by graph automorphisms.
 
Our first goal is to show that for Cayley graphs of infinite groups, \Cref{def:SF_graph} recovers the notion of space-filling action. We recall the following elementary observation: if $G$ is infinite and $G\curvearrowright X$ is a space-filling action with movement function $R$, then for every $x \in X$ the map $k\mapsto R^k_x(1_G)$ is injective.
 
\begin{prop}\label{prop:cayley_SF}
    Let $G$ be an infinite finitely generated group and $S$ a finite symmetric set of generators. Then $G$ is space-filling if and only if $\operatorname{Cay}(G,S)$ is space-filling.
\end{prop}
 
\begin{proof}
Suppose first that $G\curvearrowright X$ is a space-filling action with movement function $R$, and let $f$ be a covering function for $R$ with respect to $|\cdot|_S$ as in~\Cref{rem:covering_function}. Let $k = \max_{s \in R(X)}|s|_S$, which is finite as $R(X)$ is finite. Note that $k\geq 1$: as $X$ is nonempty, if we had $R(X) = \{1_G\}$, then $R^n_x(1_G) = 1_G$ for every $x\in X$ and $n\in\NN$, which contradicts the space-filling property since $G$ has more than one element. Fix $x\in X$ and let $w_j = R^j_x(1_G)$ for $j\in \NN$. As $w_{j+1} = w_jR(w_j^{-1}x)$, we have $d(w_j,w_{j+1}) = |R(w_j^{-1}x)|_S\leq k$, and the $w_j$ are pairwise distinct by the observation above. Now let $i,r \in \NN$ and $y = w_i^{-1}x$. By~\Cref{prop:faciles} we have for every $t \geq 0$ that
\[ w_{i+t} = R^t_x(w_i) = w_iR^t_{y}(1_G).  \]
By definition of the covering function, there is $g\in G$ such that $gB_S(r)\subset R_y^{\llbracket 0,f(r)\rrbracket}(1_G)$. It follows that \[B(w_ig,r) = w_igB_S(r)\subset \{w_i,\dots,w_{i+f(r)}\}.\] Hence $\operatorname{Cay}(G,S)$ is $(k,f)$-space-filling.
 
Conversely, suppose $\operatorname{Cay}(G,S)$ is $(k,f)$-space-filling, where by~\Cref{rem:SF_graph_basic} we may assume that $f$ is nondecreasing. By~\Cref{lem:subshift}, it suffices to show that the $(k,f)$-covering subshift $\mathfrak{X}(k,f)$ is nonempty. Let us first make the following observation. Let $x\in B_S(k)^G$, $h\in G$ and $r \in \NN$, and let $p = (h^{-1}x)|_{B_S(kf(r))}$. By~\Cref{prop:faciles}, for $j\in\llbracket 0,f(r)\rrbracket$ we have
\[ R^j_p(1_G) = R^j_{h^{-1}x}(1_G) = h^{-1}R^j_x(h). \]
It follows that $p$ is $(k,f)$-admissible for $r$ if and only if $R_x^{\llbracket 0,f(r)\rrbracket}(h)$ contains a ball of radius $r$. In particular, whether $(h^{-1}x)|_{B_S(kf(r))}$ is $(k,f)$-admissible for $r$ only depends on the restriction of $x$ to the finite set $hB_S(kf(r))$.
 
Now let $n \in \NN$ and set $M_n = 2f(n)$. Take a $k$-path $(w_j)_{j=0}^{M_n}$ as in~\Cref{def:SF_graph} and let $c_n\in G$ be such that $B(c_n,n)\subset\{w_0,\dots,w_{f(n)}\}$. Define $y_n\in B_S(k)^G$ by $y_n(w_j) = w_j^{-1}w_{j+1}$ for $j<M_n$, which belongs to $B_S(k)$ because $d(w_j,w_{j+1})\leq k$, and $y_n(g) = 1_G$ otherwise. Then $R^t_{y_n}(w_j) = w_{j+t}$ whenever $j+t \leq M_n$. Let $x_n = c_n^{-1}y_n$.
 
We claim that for every $h\in B_S(n)$ and $r\leq n$ we have that $(h^{-1}x_n)|_{B_S(kf(r))}$ is $(k,f)$-admissible for $r$. Indeed, $c_nh \in B(c_n,n)$, so $c_nh = w_j$ for some $j\leq f(n)$. As $f$ is nondecreasing, $j+f(r)\leq 2f(n) = M_n$, and thus we obtain
\[ R^{\llbracket 0,f(r)\rrbracket}_{x_n}(h) = c_n^{-1}R^{\llbracket 0,f(r)\rrbracket}_{y_n}(c_nh) = c_n^{-1}\{w_j,\dots,w_{j+f(r)}\},   \]
which contains a ball of radius $r$ because $\{w_j,\dots,w_{j+f(r)}\}$ does. The claim follows from the observation above.
 
Let $x$ be an accumulation point of $(x_n)_{n \in \NN}$, and let $h\in G$ and $r\in\NN$. For every $n\geq \max\{|h|_S,r\}$ the pattern $(h^{-1}x_n)|_{B_S(kf(r))}$ is $(k,f)$-admissible for $r$. As this only depends on the restriction of $x_n$ to the finite set $hB_S(kf(r))$, the pattern $(h^{-1}x)|_{B_S(kf(r))}$ is also $(k,f)$-admissible for $r$. Therefore $x\in \mathfrak{X}(k,f)$, which is thus nonempty.
\end{proof}


\subsection{Quasi-isometry invariance}
 
Let $\Gamma,\Delta$ be graphs and $L\geq 1, D\geq 0$. A map $\zeta\colon \Delta\to\Gamma$ is an \define{$(L,D)$-quasi-isometry} if
\begin{itemize}
    \item For every $x,y \in \Delta$
    \[\frac{1}{L}d_{\Delta}(x,y)-D\leq d_{\Gamma}(\zeta(x),\zeta(y))\leq L\,d_{\Delta}(x,y)+D.\]
    \item For every $y\in\Gamma$ there is $x\in\Delta$ such that $d_{\Gamma}(y,\zeta(x))\leq D$.
\end{itemize}

We often omit the constants and speak plainly of a quasi-isometry. Two graphs are \define{quasi-isometric} if there is a quasi-isometry between them. The proof that being space-filling is a quasi-isometry invariant of graphs proceeds in three steps: we first show that the property can be pulled back through surjective quasi-isometries, then that we can always produce such a surjection after replacing the graph by a thickened version, and finally that the property passes from the thickened graph back to the original.
 
\begin{lemma}\label{lem:pullback}
    Let $\Gamma,\Delta$ be graphs and let $\pi\colon \Delta\to\Gamma$ be a surjective quasi-isometry. If $\Gamma$ is space-filling, then so is $\Delta$.
\end{lemma}
 
\begin{proof}
    Let $L\geq 1, D\geq 0$ be integers such that $\pi$ is an $(L,D)$-quasi-isometry. 
    By definition, we have $d_\Delta(x,y)\leq L\bigl(d_\Gamma(\pi(x),\pi(y))+D\bigr)$ for every $x,y\in\Delta$. In particular, every fiber $\pi^{-1}(w)$ has diameter at most $LD$. Hence fibers have at most $b=\beta_\Delta(LD)$ elements, and they are nonempty since $\pi$ is surjective. Suppose $\Gamma$ is $(k,f)$-space-filling and set for $r\geq 0$ \[ k' =  L(k+D), \quad \rho(r) =  Lr+D \quad \mbox{ and }\quad f'(r) = b\bigl(f(\rho(r))+2\bigr).  \]
    
    Let $M \in \NN$ and take a $k$-path $(w_j)_{j=0}^M$ in $\Gamma$. Let $(w'_t)_{t=0}^{M'}$ be the sequence obtained by concatenating, in the order $j=0,\dots,M$, arbitrary enumerations of the fibers $\pi^{-1}(w_j)$, see~\Cref{fig:thicc_paths}. Its terms are pairwise distinct, as the fibers are pairwise disjoint and the $w_j$ are pairwise distinct, and $M'\geq M$ as the fibers are nonempty. Two consecutive terms in the same fiber are at distance at most $LD$, and two consecutive terms in the fibers of $w_j$ and $w_{j+1}$ are at distance at most $L\bigl(d(w_j,w_{j+1})+D\bigr)\leq L(k+D)$. Hence $(w'_t)_{t=0}^{M'}$ is a $k'$-path.

    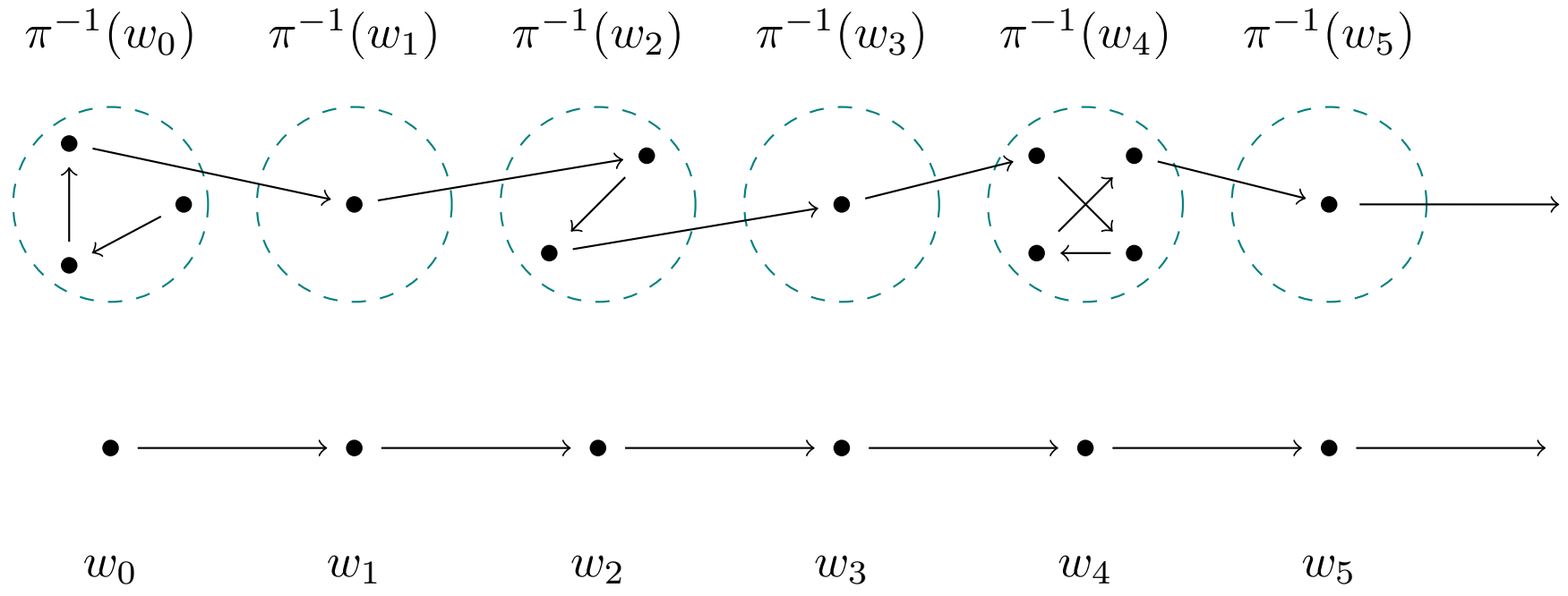
\begin{figure}
        \centering
        \begin{tikzpicture}[scale = 1.8]
            \foreach \i in {0,...,5}{
            \draw[->, shorten <=0.2cm, shorten >=0.2cm] (\i,0.5) -- (\i+1,0.5);
            \draw[fill=black] (\i,0.5) circle (0.03);
            \draw[dashed, color=teal] (\i,1.5) circle (0.4);
            \node at (\i,2.2) {$\pi^{-1}(w_{\i})$};
            \node at (\i,0) {$w_{\i}$};
            }
            \node (AA) at (0.3,1.5) {};
            \draw[fill=black] (0.3,1.5)  circle (0.03);
            \node (AB) at (-0.17,1.25) {};
            \draw[fill=black] (-0.17,1.25) circle (0.03);
            \node (AC) at (-0.17,1.75){};
            \draw[fill=black] (-0.17,1.75) circle (0.03);
            \node (BB) at (1,1.5){};
            \draw[fill=black] (1,1.5) circle (0.03);
            \node (CA) at (2.2,1.7){};
            \draw[fill=black] (2.2,1.7) circle (0.03);
            \node (CB) at (1.8,1.3){};
            \draw[fill=black] (1.8,1.3) circle (0.03);
            \node (DD) at (3,1.5){};
            \draw[fill=black] (3,1.5) circle (0.03);
            \node (EA) at (4.2,1.7){};
            \draw[fill=black] (4.2,1.7) circle (0.03);
            \node (EB) at (3.8,1.7){};
            \draw[fill=black] (3.8,1.7) circle (0.03);
            \node (EC) at (3.8,1.3){};
            \draw[fill=black] (3.8,1.3) circle (0.03);
            \node (ED) at (4.2,1.3){};
            \draw[fill=black] (4.2,1.3) circle (0.03);
            \node (FF) at (5,1.5){};
            \draw[fill=black] (5,1.5) circle (0.03);
            \draw[->, shorten <=0.05cm, shorten >=0.05cm] (AA) to (AB);
            \draw[->, shorten <=0.05cm, shorten >=0.05cm] (AB) to (AC);
            \draw[->, shorten <=0.05cm, shorten >=0.05cm] (AC) to (BB);
            \draw[->, shorten <=0.05cm, shorten >=0.05cm] (BB) to (CA);
            \draw[->, shorten <=0.05cm, shorten >=0.05cm] (CA) to (CB);
            \draw[->, shorten <=0.05cm, shorten >=0.05cm] (CB) to (DD);
            \draw[->, shorten <=0.05cm, shorten >=0.05cm] (DD) to (EB);
            \draw[->, shorten <=0.05cm, shorten >=0.05cm] (EB) to (ED);
            \draw[->, shorten <=0.05cm, shorten >=0.05cm] (ED) to (EC);
            \draw[->, shorten <=0.05cm, shorten >=0.05cm] (EC) to (EA);
            \draw[->, shorten <=0.05cm, shorten >=0.05cm] (EA) to (FF);
            \draw[->, shorten <=0.1cm, shorten >=0.1cm] (FF) to (6,1.5);
        \end{tikzpicture}
        \caption{Construction of the path $w'$ from $w$.}
        \label{fig:thicc_paths}
    \end{figure}
 
    Let $r\in \NN$ and let $W' = \{w'_i,\dots,w'_{i+f'(r)}\}$ be a window with $i+f'(r)\leq M'$. Notice that fibers have at most $b$ elements and $W'$ consists of $f'(r)+1$ consecutive terms. Then, at most two fibers meet $W'$ without being contained in it, and $W'$ contains every element of at least $\lceil (f'(r)+1)/b\rceil -2\geq f(\rho(r))+1$  consecutive fibers. That is, there is $i_0$ with $i_0+f(\rho(r))\leq M$ such that $\pi^{-1}(\{w_{i_0},\dots,w_{i_0+f(\rho(r))}\})\subset W'$. By the choice of the path in $\Gamma$, there is $v \in \Gamma$ such that $B_\Gamma(v,\rho(r))\subset\{w_{i_0},\dots,w_{i_0+f(\rho(r))}\}$. Choose $x \in \pi^{-1}(v)$. For $y\in \Delta$, if $d(x,y)\leq r$, then $d(v,\pi y)\leq Lr+D = \rho(r)$. Thus, $y\in \pi^{-1}(B_\Gamma(v,\rho(r)))\subset W'$ and $B_\Delta(x,r)\subset W'$.
 
    By the observation following~\Cref{def:SF_graph}, the $k'$-paths $(w'_t)_{t=0}^{M'}$ show that $\Delta$ is $(k',f')$-space-filling.
\end{proof}
 
Given a graph $\Delta$ and $m\geq 1$, let $\Delta\boxtimes K_m$ be its product with the complete graph on $m$ vertices. Namely, $\Delta\boxtimes K_m$ is the graph with vertex set $\Delta\times\{1,\dots,m\}$ in which $(x,i)$ and $(y,j)$ are adjacent if either $x$ and $y$ are adjacent in $\Delta$, or $x=y$ and $i\neq j$. This graph is connected and has bounded degree. Denote by $p\colon \Delta\boxtimes K_m\to\Delta$ the projection $p(x,i)=x$. As $p$ maps adjacent vertices to adjacent or equal vertices, and any path in $\Delta$ between distinct vertices can be lifted to a path of the same length with arbitrary indices at its endpoints, we have
\begin{equation*}
    d((x,i),(y,j)) = \begin{cases}
        d(x,y) & \mbox{if } x\neq y,\\
        1 & \mbox{if } x= y \mbox{ and } i\neq j.
    \end{cases}
\end{equation*}
In particular $d(p(a),p(b))\leq d(a,b)\leq d(p(a),p(b))+1$, and for every $r\geq 1$
\begin{equation*}
    B_{\Delta\boxtimes K_m}((y,i),r) = B_\Delta(y,r)\times\{1,\dots,m\}.
\end{equation*}
 
\begin{lemma}\label{lem:thicken}
Let $\Gamma,\Delta$ be graphs and $\zeta\colon\Delta\to\Gamma$ an $(L,D)$-quasi-isometry. Then there exist $m\geq 1$ and a surjective $(L,3LD+1)$-quasi-isometry $\pi\colon \Delta\boxtimes K_m\to\Gamma$.
\end{lemma}
 
\begin{proof}
    For every $w\in \Gamma$ choose $\rho(w)\in\Delta$ with $d_{\Gamma}(\zeta(\rho(w)),w)\leq D$. Then, for every $v\in \Delta$, $\rho^{-1}(v)\subset B_\Gamma(\zeta(v),D)$, so $|\rho^{-1}(v)|\leq \beta_\Gamma(D)$. Let $m=\beta_\Gamma(D)+1$. For each $v\in\Delta$ enumerate $\rho^{-1}(v)=\{w_1,\dots,w_{r_v}\}$, where $r_v<m$, and set $\pi(v,i)=w_i$ for $i\leq r_v$ and $\pi(v,i) = \zeta(v)$ for $r_v<i\leq m$. Every $w\in \Gamma$ is of the form $\pi(\rho(w),i)$ for some $i$, so $\pi$ is surjective. Moreover, $d(\pi(a),\zeta(p(a)))\leq D$ for every $a\in \Delta\boxtimes K_m$. Therefore, for every $a,b\in \Delta\boxtimes K_m$,
    \begin{align*}
        d_{\Gamma}(\pi(a),\pi(b))& \leq d_{\Gamma}(\zeta(p(a)),\zeta(p(b)))+2D\\ & \leq L\,d_{\Delta}(p(a),p(b))+3D\\ &\leq L\,d_{\Delta\boxtimes K_m}(a,b)+3LD+1,
    \end{align*}
    and, using that $d_{\Delta}(u,v)\leq L(d_{\Gamma}(\zeta(u),\zeta(v))+D)$ for $u,v\in \Delta$,
\begin{align*}
    d_{\Delta\boxtimes K_m}(a,b) & \leq d_{\Delta}(p(a),p(b))+1\\ & \leq L\bigl(d_{\Gamma}(\zeta(p(a)),\zeta(p(b)))+D\bigr)+1\\ & \leq L\,d_{\Gamma}(\pi(a),\pi(b))+3LD+1.
\end{align*}
Hence $\pi$ is an $(L,3LD+1)$-quasi-isometry.
\end{proof}

\begin{lemma}\label{lem:thin}
    Let $\Delta$ be a graph and $m\geq 1$. If $\Delta\boxtimes K_m$ is space-filling, then so is $\Delta$.
\end{lemma}
 
\begin{proof}
    Write $\Delta' = \Delta\boxtimes K_m$ and suppose that $\Delta'$ is $(k,f)$-space-filling. Let $M\in \NN$ and $N = m(M+1)-1$, and take a $k$-path $W=(w_t)_{t=0}^N$ in $\Delta'$. Let $P_W = \{p(w_t)\in \Delta : 0\leq t\leq N\}$ and for $u\in P_W$ let $T_u=\{t \in \llbracket 0,N\rrbracket : p(w_t)=u\}$ be the set of positions at which $u$ occurs. As the $w_t$ are pairwise distinct, $|T_u|\leq m$, and in particular $|P_W|\geq (N+1)/m = M+1$.
 
    Partition the set $\llbracket 0,N\rrbracket$ into the consecutive blocks $I_i = \llbracket (i-1)m, im-1\rrbracket$ for $i\in\{1,\dots, (N+1)/m\}$. Consider the bipartite graph which joins each block $I_i$ to every $u\in P_W$ with $T_u\cap I_i\neq\varnothing$. Any $q$ distinct blocks contain at least $(q-1)m+1$ positions, and each $u \in P_W$ occupies at most $m$ of them, so these $q$ blocks have at least $\lceil ((q-1)m+1)/m\rceil = q$ neighbors. By Hall's theorem, there is an injective map $\mu$ from the set of blocks to $P_W$ such that $T_{\mu(i)}\cap I_i\neq\varnothing$ for every $i$. For each $u\in P_W$ we choose $t_u \in T_u$ as follows: if $u=\mu(i)$ for some (necessarily unique) $i$, we choose $t_u\in T_u\cap I_i$, otherwise we choose $t_u \in T_u$ arbitrarily. With this choice, every block contains $t_u$ for some $u \in P_W$.
 
    Let $J = \{t_u : u \in P_W\}$ and enumerate it increasingly as $j_0<j_1<\dots<j_{|P_W|-1}$. Notice that by construction we have $|j_{s+1}-j_s|<2m$. Let $v_s = p(w_{j_s})$ which by construction are pairwise distinct and satisfy
    \[ d(v_s,v_{s+1})\leq d(w_{j_s},w_{j_{s+1}})\leq k|j_{s+1}-j_s| \leq (2m-1)k. \]
    Hence $V=(v_s)_{s=0}^{|P_W|-1}$ is a $(2m-1)k$-path in $\Delta$ with at least $M+1$ terms.
 
    Let $\{v_s,\dots,v_{s+f(r)}\}$ be a window of $V$. If $r=0$ there is nothing to prove, as every nonempty set contains a ball of radius $0$, so suppose that $r\geq 1$. As the $j$'s are distinct integers, $j_{s+f(r)}\geq j_s+f(r)$, so $\{w_{j_s},\dots,w_{j_s+f(r)}\}$ is a window of $W$. Thus, it contains a ball of radius $r$ in $\Delta'$, which is thus of the form $B_\Delta(y,r)\times\{1,\dots,m\}$. Hence, for every $u\in B_\Delta(y,r)$ all copies of $u$ occur at positions in $\llbracket j_s,j_s+f(r)\rrbracket$; in particular $t_u\in\llbracket j_s,j_{s+f(r)}\rrbracket$, that is, $u\in \{v_s,\dots,v_{s+f(r)}\}$. Therefore $B_\Delta(y,r)\subset \{v_s,\dots,v_{s+f(r)}\}$. As $M$ can be taken arbitrarily large, it follows that $\Delta$ is $((2m-1)k,f)$-space-filling.
\end{proof}

\begin{theorem}\label{thm:QI_graphs}
    Let $\Gamma$ and $\Delta$ be quasi-isometric graphs. Then $\Gamma$ is space-filling if and only if $\Delta$ is space-filling.
\end{theorem}
 
\begin{proof}
    As being quasi-isometric is a symmetric relation, it suffices to show that if $\Gamma$ is space-filling, then so is $\Delta$. Let $\zeta\colon \Delta\to\Gamma$ be a quasi-isometry. By~\Cref{lem:thicken} there are $m\geq 1$ and a surjective quasi-isometry $\pi\colon \Delta\boxtimes K_m\to\Gamma$. By~\Cref{lem:pullback}, $\Delta\boxtimes K_m$ is space-filling, and by~\Cref{lem:thin}, so is $\Delta$.
\end{proof}

\begin{corollary}\label{cor:QI_groups}
    Let $G$ and $H$ be quasi-isometric finitely generated groups. Then $G$ is space-filling if and only if $H$ is space-filling.
\end{corollary}
 
\begin{proof}
    If $G$ is finite the result is obvious by~\Cref{prop:finite_groups}. Otherwise both groups are infinite and their Cayley graphs with respect to finite symmetric sets of generators are quasi-isometric graphs. The result follows from~\Cref{prop:cayley_SF} and~\Cref{thm:QI_graphs}.
\end{proof}

\subsection{Computable space-filling on graphs}
 
In this subsection we look at an effective version of the space-filling property for graphs. In this case, the natural setting is that of highly computable graphs, which play the role that groups with decidable word problem play for groups.
 
\begin{defn}\label{def:highly_computable}
    A graph $\Gamma$ is \define{highly computable} if its vertex set is a decidable subset of $\NN$, its adjacency relation is decidable, and the map $v\mapsto \deg(v)$ which sends each vertex to its degree is computable.
\end{defn}

Notice that if $\Gamma$ is a highly computable graph, there are algorithms which, given $v,w\in \Gamma$ and $r\in\NN$, compute the finite list of elements of $B_\Gamma(v,r)$ and the distance $d(v,w)$.
 
 
\begin{defn}\label{def:effective_SF_graph}
    A highly computable graph $\Gamma$ is \define{computably space-filling} if there are $k\geq 1$ and a computable function $f\colon\NN\to\NN$ such that $\Gamma$ is $(k,f)$-space-filling.
\end{defn}
 
The following lemma shows that in a computably space-filling graph the $k$-paths required by~\Cref{def:SF_graph} can be computed.
 
\begin{lemma}\label{lem:effective_paths}
    Let $\Gamma$ be a highly computable graph which is $(k,f)$-space-filling with $f$ computable. There is an algorithm which, given $M\in\NN$, computes a $k$-path $(w_j)_{j=0}^M$ which satisfies the space-filling condition.
\end{lemma}
 
\begin{proof}
    As $\Gamma$ contains $k$-paths with arbitrarily many pairwise distinct vertices, it is infinite. As $\Gamma$ is connected, every ball of radius $r$ contains at least $r+1$ vertices, and thus the space-filling condition forces $f(r)\geq r$ for every $r$. Hence, for a given $M$, the condition only concerns the finitely many pairs $(r,i)$ with $r\leq M$, $f(r)\leq M$ and $i\leq M-f(r)$. For each such pair, a window $W$ contains a ball of radius $r$ if and only if there is $v\in W$ with $B_\Gamma(v,r)\subset W$. Given a finite sequence of vertices $(w_j)_{j=0}^M$ we can therefore decide whether it is a $k$-path satisfying the space-filling condition. The algorithm enumerates all finite sequences of $M+1$ vertices and outputs the first one which passes this test. It halts because such a sequence exists.
\end{proof}
 
Next we observe that the space-filling parameters produced in the proof of~\Cref{thm:QI_graphs} are obtained from the original ones by explicit formulas.
 
\begin{theorem}\label{thm:QI_effective_graphs}
    Let $\Gamma$ and $\Delta$ be quasi-isometric highly computable graphs. Then $\Gamma$ is computably space-filling if and only if $\Delta$ is computably space-filling.
\end{theorem}

\begin{proof}
    As being quasi-isometric is a symmetric relation, it suffices to prove one implication. Suppose $\Gamma$ is $(k,f)$-space-filling with $f$ computable and let $\zeta\colon\Delta\to\Gamma$ be an $(L,D)$-quasi-isometry for some integers $L$ and $D$.
    
   Let $m=\beta_\Gamma(D)+1$ and $\Delta'=\Delta\boxtimes K_m$. By~\Cref{lem:thicken}, there is a surjective quasi-isometry $\Delta'\to\Gamma$ with constants $L$ and $D'=3LD+1$. By the proof of~\Cref{lem:pullback}, $\Delta'$ is $(k',f')$-space-filling, where $k'= L(k+D')$, $f'(r) = b(f( Lr+D')+2)$ and $b = \beta_{\Delta'}(LD')$. By the proof of~\Cref{lem:thin}, $\Delta$ is $((2m-1)k',f')$-space-filling.
    
    As $f$ is computable, the function $f'$ is also computable. As $\Delta$ is highly computable, we conclude that it is computably space-filling.
\end{proof}

We note that the quasi-isometry and the auxiliary maps used in the proofs of~\Cref{lem:pullback,lem:thicken,lem:thin} (the choice function $\rho$, the surjection onto $\Gamma$ and the matchings given by Hall's theorem) do not need to be computable, as they only appear in the proof of existence of the $k$-paths and not in the space-filling parameters.

\begin{remark}
    The proof of~\Cref{thm:QI_effective_graphs} does not use that $\Delta$ is highly computable: every graph quasi-isometric to a computably space-filling graph is $(k,f)$-space-filling for some $k\geq 1$ and some computable $f$. The assumption is included in~\Cref{def:effective_SF_graph} because computable space-filling parameters are only useful when the graph itself can be handled algorithmically, as in~\Cref{lem:effective_paths}.
\end{remark}
 
We now turn to finitely generated groups with decidable word problem. First, notice that for every symmetric finite generating set $S$, the Cayley graph $\operatorname{Cay}(G,S)$ is highly computable. Also note that the first part of the proof of~\Cref{prop:cayley_SF} shows the following: if $G$ is infinite, $G\curvearrowright X$ is a space-filling action with movement function $R$ and covering function $f$, and $k = \max_{g \in R(X)}|g|_S$, then $\operatorname{Cay}(G,S)$ is $(k,f)$-space-filling.

\begin{prop}\label{prop:effective_SF}
    Let $G$ be an infinite finitely generated group with decidable word problem and $S$ a finite symmetric set of generators. Then, $G$ is computably space-filling if and only if $\operatorname{Cay}(G,S)$ is computably space-filling.
\end{prop}
 
\begin{proof}
    Suppose $G$ is computably space-filling. Then $G$ admits a space-filling action with a computable covering function $f$ by~\Cref{prop:computable_cov_function}. As $\operatorname{Cay}(G,S)$ is highly computable, the observation above shows that it is $(k,f)$-space-filling for some $k\geq 1$, and thus computably space-filling.
 
    Next, suppose $\operatorname{Cay}(G,S)$ is $(k,f)$-space-filling with $f$ computable, where by~\Cref{rem:SF_graph_basic} we may assume that $f$ is nondecreasing (which preserves computability). By the proof of~\Cref{prop:cayley_SF}, the $(k,f)$-covering subshift $\mathfrak{X}(k,f)$ is nonempty, and thus space-filling by~\Cref{lem:subshift}. As $f$ is computable and $G$ has decidable word problem, it is effective.
\end{proof}
 
\begin{corollary}\label{cor:QI_effective}
    Let $G$ and $H$ be quasi-isometric finitely generated groups with decidable word problem. Then $G$ is computably space-filling if and only if $H$ is computably space-filling.
\end{corollary}
 
\begin{proof}
    If $G$ is finite, then so is $H$, and the space-filling subshift for $H$ given by~\Cref{prop:finite_groups} is finite and thus effective. Suppose that $G$ and $H$ are infinite and fix finite symmetric sets of generators $S_G$ and $S_H$ and consider the highly computable Cayley graphs $\operatorname{Cay}(G,S_G)$ and $\operatorname{Cay}(H,S_H)$, which are quasi-isometric. By~\Cref{prop:effective_SF} the first one is computably space-filling and by~\Cref{thm:QI_effective_graphs}, so is $\operatorname{Cay}(H,S_H)$, and  hence by~\Cref{prop:effective_SF} again, $H$ admits an effective space-filling subshift.
\end{proof}

Finally, we remark that if instead we only care about the existence of a computable covering function, then the result holds for all finitely generated groups.

\begin{corollary}
      Let $G$ and $H$ be quasi-isometric finitely generated groups. If $G$ admits a space-filling action with a computable covering function, then so does $H$.
\end{corollary}

\begin{proof}
    If $G$ is finite the claim is clear. Otherwise, by the first part of the proof of~\Cref{prop:cayley_SF}, $\operatorname{Cay}(G,S_G)$ is $(k,f)$-space-filling with $f$ computable. By the proof of~\Cref{thm:QI_effective_graphs}, which does not use that the graphs are highly computable, $\operatorname{Cay}(H,S_H)$ is $(k',f')$-space-filling with $f'$ computable and nondecreasing. By the second part of the proof of~\Cref{prop:cayley_SF}, $\mathfrak{X}(k',f')$ is nonempty, and by~\Cref{lem:subshift} it is space-filling with covering function $f'$.
\end{proof}
 
\subsection{Finitary space-filling quasi-isometry invariance}\label{subsec:QI_local}

In this subsection we show that being finitary space-filling is a quasi-isometry invariant among finitely presented groups. In broad strokes, the proof follows the same three steps as the proof of~\Cref{thm:QI_graphs}, the main difference being that now we need to encode a family of quasi-isometries in a closed and shift-invariant manner by local rules. This is achieved through a construction originally due to Cohen~\cite{Cohen_2017} (see also~\cite{BarBit2026} for a generalization beyond the context of groups), and it is the only point where finite presentability is used (the relations of the group need to be encoded through local rules). We also replace the use of Hall's theorem in~\Cref{lem:thin} by a first return map.

Throughout this subsection, $G$ and $H$ are finitely generated groups with finite symmetric sets of generators $S_G$ and $S_H$. We write $d_G$ and $d_H$ for the associated word metrics, and say that a map $\zeta\colon G\to H$ is an $(L,D)$-quasi-isometry if it is an $(L,D)$-quasi-isometry from $\operatorname{Cay}(G,S_G)$ to $\operatorname{Cay}(H,S_H)$.

\begin{prop}\label{prop:cohen_encoding}
    Let $G$ and $H$ be finitely presented quasi-isometric groups and fix a surjective quasi-isometry $\zeta\colon G\to H$ such that $\zeta(1_G)=1_H$. There exist a finite alphabet $A$, a $G$-SFT $Q\subset A^G$ and, for every $q \in Q$, a map $\varphi_q\colon G\to H$ such that the following hold.
    \begin{enumerate}
        \item For every $q\in Q$ we have $\varphi_q(1_G) = 1_H$, and for every $g,h\in G$
        \[ \varphi_q(gh) = \varphi_q(g)\varphi_{g^{-1}q}(h). \]
        Moreover, for every $n\in\NN$ the restriction of $\varphi_q$ to $B_{S_G}(n)$ only depends on $q|_{B_{S_G}(n)}$.
        \item There exist constants $L\geq 1$, $D\geq 0$ such that for every $q \in Q$ the map $\varphi_q$ is a surjective $(L,D)$-quasi-isometry.
        \item There exists $q\in Q$ with $\varphi_q = \zeta$.
    \end{enumerate}
\end{prop}

\Cref{prop:cohen_encoding} follows from a minor modification of the construction of Cohen~\cite{Cohen_2017}. Roughly speaking, in the original construction $Q$ encodes a space of quasi-isometries with uniformly bounded constants as above, but without any constraint on whether they are surjective or not. However, it is quite simple to impose through local rules that for every $g\in G$ and $t\in S_H$ there is $g'\in gB_{S_G}(L(1+D))$ such that $\varphi_q(g') = \varphi_q(g)t$. As $\varphi_q(1_G)=1_H$, this forces $\varphi_q$ to be surjective, and it holds for $\zeta$ as soon as $\zeta$ is a surjective $(L,D)$-quasi-isometry. We will use this construction as a black box in what follows.

The next result is the analogue of~\Cref{lem:pullback}. The fibers of the encoded quasi-isometry are traversed as in~\Cref{fig:thicc_paths}. However, since there is in general no way to choose an order on each fiber which is invariant under the shift, the order in which a fiber is traversed is encoded in the configuration itself.

\begin{prop}\label{prop:local_pullback}
    Let $G$ and $H$ be finitely presented groups and suppose there is a surjective quasi-isometry $\zeta\colon G\to H$. If $H$ is finitary space-filling, then so is $G$.
\end{prop}

\begin{proof}
    Let $\zeta$ be a quasi-isometry which satisfies (without loss of generality) that $\zeta(1_G) = 1_H$. Consider the $G$-SFT $Q\subset A^G$ and the maps $(\varphi_q)_{q\in Q}$ given by~\Cref{prop:cohen_encoding}. In particular, there is $\bar{q}\in Q$ such that $\zeta = \varphi_{\bar{q}}$. Fix integers $L\geq 1$ and $D\geq 0$ large enough such that $\varphi_q$ is a surjective $(L,D)$-quasi-isometry for every $q\in Q$. 
    
    Let $E = LD$ and $b = |B_{S_G}(E)|$. For $q\in Q$ and $g\in G$, the \define{fiber} of $g$ is the set $\mathcal{C}_q(g) = \varphi_q^{-1}(\varphi_q(g))$. As $\varphi_q$ is an $(L,D)$-quasi-isometry, we have $d_G(g,g')\leq L(d_H(\varphi_q(g),\varphi_q(g'))+D)$ for every $g,g'\in G$. In particular $g^{-1}\mathcal{C}_q(g)\subset B_{S_G}(E)$ and thus $|\mathcal{C}_q(g)|\leq b$. 

    Let $Y\subset B^H$ be a space-filling SFT with movement function $T$ and covering function $f_Y$. As $T$ is continuous and takes finitely many values, there are a finite set of forbidden patterns $\FF$ and $r_0\in \NN$ such that $Y=\X_{\FF}$, every pattern in $\FF$ has support contained in $B_{S_H}(r_0)$, and $T(y)$ only depends on $y|_{B_{S_H}(r_0)}$. Let $c= \max_{y\in Y}|T(y)|_{S_H}$.

    Now we are ready to define a space-filling $G$-SFT. To that end, we consider the space $Z$ which consists of all configurations
    \[z = (q,\beta,\nu,\iota) \in \bigl(A\times B\times (B_{S_G}(E)\cup\{\bot\})\times \{0,1\}\bigr)^G,\] where $q \in Q$ and furthermore, $z$ satisfies two conditions:
    \begin{enumerate}
        \item The map $\beta$ is constant on $\mathcal{C}_q(g)$ for every $g\in G$, and the configuration $y_z\in B^H$ given by $y_z(\varphi_q(g)) = \beta(g)$ belongs to $Y$.
        \item For every $g\in G$, there is an enumeration $g_1,\dots,g_\ell$ of $\mathcal{C}_q(g)$ such that $\iota(g_1) = 1$, $\iota(g_i)=0$ for $i\geq 2$, $g_{i+1} = g_i\nu(g_i)$ for $i<\ell$, and $\nu(g_\ell) = \bot$.
    \end{enumerate}
    
    In other words, for each $q \in Q$ the coordinate $\beta$ encodes a configuration of $Y$ through $\varphi_q$, and $\nu,\iota$ encode a linear order on each fiber: $\iota$ marks the first element, $\nu$ points to the next element, and $\bot$ marks the last one. Note that $y_z$ is well defined on all of $H$, as $\varphi_q$ is surjective. We claim that $Z$ is a nonempty space-filling $G$-SFT and thus that $G$ is finitary space-filling.

    Let us first argue that $Z$ is a $G$-SFT. As $Q$ is a $G$-SFT, the condition that $q\in Q$ can obviously be implemented with finitely many forbidden patterns. For $g,k\in G$ we have $\varphi_q(g)^{-1}\varphi_q(gk) = \varphi_{g^{-1}q}(k)$, which only depends on $q|_{gB_{S_G}(|k|_{S_G})}$ by~\Cref{prop:cohen_encoding}. Since $\mathcal{C}_q(g)\subset gB_{S_G}(E)$, both the condition that $\beta$ is constant and the enumeration condition only depend on the restriction of $z$ to $gB_{S_G}(2E)$ and thus can be implemented locally. Finally, $\varphi_q(g)^{-1}y_z$ contains a pattern of $\FF$ at $1_H$ if and only if $y_z$ contains a translate of it at $\varphi_q(g)$. As every element of $\varphi_q(g)B_{S_H}(r_0)$ is the image of some element of $gB_{S_G}(L(r_0+D))$, this only depends on the restriction of $z$ to $gB_{S_G}(L(r_0+D))$. As $\varphi_q$ is surjective, we have that $y_z \in Y$ holds if and only if no pattern of $\FF$ appears at $\varphi_q(g)$ for any $g\in G$. This shows that $Z$ is a $G$-SFT.

    Next, we define the movement function $R\colon Z\to G$. Let $z = (q,\beta,\nu,\iota)\in Z$. If $\nu(1_G)\neq\bot$, we let $R(z) = \nu(1_G)$. Otherwise, we let $R(z)$ be the unique element $g$ of $\varphi_q^{-1}(T(y_z))$ with $\iota(g) = 1$. As $\varphi_q(1_G)=1_H$, we have $g\in B_{S_G}(L(c+D))$, and thus $R$ takes finitely many values. Moreover, $R(z)$ only depends on the restriction of $z$ to a finite set by the arguments above, so it is continuous.

    Note that for $z=(q,\beta,\nu,\iota)\in Z$ and $g\in G$ we have $\mathcal{C}_{g^{-1}q}(k) = g^{-1}\mathcal{C}_q(gk)$ and $y_{g^{-1}z} = \varphi_q(g)^{-1}y_z$. Recalling that $R_{z}(g) = gR(g^{-1}z)$ and $T_{y}(g) = gT(g^{-1}y)$, we deduce that for every $g\in G$:
    \begin{itemize}
        \item if $\nu(g)\neq \bot$, then $R_{z}(g) = g\nu(g)$ is the next element of $\mathcal{C}_q(g)$;
        \item if $\nu(g) = \bot$, then $R_{z}(g)$ is the first element of $\varphi_q^{-1}\bigl(T_{y_z}(\varphi_q(g))\bigr)$.
    \end{itemize}
    Let $h_j = T^j_{y_z}(1_H)$ for $j\geq 0$. It follows that the sequence $(R^t_{z}(1_G))_{t\geq 0}$ first enumerates the elements of $\varphi_q^{-1}(h_0)$ from $1_G$ up to the last element of this fiber, and then, for $j=1,2,\dots$, enumerates all elements of $\varphi_q^{-1}(h_j)$ from the first to the last one. In particular, every fiber is traversed in at most $b$ consecutive steps.

    Let us show that $R$ is a movement function which makes $Z$ space-filling. Fix $z\in Z$ and let $r\in \NN$. Set $\rho(r) =  Lr+D$ and define \[f_Z(r) = b\bigl(f_Y(\rho(r))+2\bigr).\] We claim that $f_Z$ is a covering function for $R$. By~\Cref{prop:faciles} we have $h_{1+j} = h_1T^j_{h_1^{-1}y_z}(1_H)$ for every $j\geq 0$. Applying the definition of the covering function to $h_1^{-1}y_z$, we obtain $v\in H$ such that \[vB_{S_H}(\rho(r))\subset\{h_1,\dots,h_{1+f_Y(\rho(r))}\}.\] By the description above, the fibers of $h_0,h_1,\dots,h_{1+f_Y(\rho(r))}$ are entirely visited by $R_{z}^{\llbracket 0,f_Z(r)\rrbracket}(1_G)$, except possibly for the elements of the fiber of $h_0$ which precede $1_G$. Choose $w\in\varphi_q^{-1}(v)$. If $g\in wB_{S_G}(r)$, then $d_H(\varphi_q(g),v)\leq Lr+D=\rho(r)$, so $\varphi_q(g) = h_j$ for some $j\in\llbracket 1,1+f_Y(\rho(r))\rrbracket$, and thus $g\in R_{z}^{\llbracket 0,f_Z(r)\rrbracket}(1_G)$. Hence $wB_{S_G}(r)\subset R_{z}^{\llbracket 0,f_Z(r)\rrbracket}(1_G)$, and $f_Z$ is a covering function for $R$.

    Finally we show that $Z$ is nonempty. To that end, recall that there exists $\bar{q}\in Q$ such that $\varphi_{\bar{q}}=\zeta$. Take $y\in Y$, let $\beta(g) = y(\zeta(g))$ for $g\in G$, and choose an arbitrary linear order on each fiber $\zeta^{-1}(h)$, $h\in H$, which we encode through $\nu$ and $\iota$ as above. The resulting configuration $z = (\bar{q},\beta,\nu,\iota)$ satisfies the constraints with $y_z=y$, and thus $z\in Z$.\end{proof}

The next result is the analogue of~\Cref{lem:thin}. Recall that if $S_G$ is a finite symmetric set of generators of $G$, then
\[S' = \bigl((S_G\cup\{1_G\})\times \ZZ/m\ZZ\bigr)\setminus\{(1_G,0)\}\]
is a finite symmetric set of generators of $G\times\ZZ/m\ZZ$. Since $(g,i)(s,c) = (gs,i+c)$, the Cayley graph $\operatorname{Cay}(G\times\ZZ/m\ZZ,S')$ is precisely $\operatorname{Cay}(G,S_G)\boxtimes K_m$. In particular, for every $r\geq 1$ the ball of radius $r$ in $G\times\ZZ/m\ZZ$ with respect to $S'$ centered at $(g,i)$ is $gB_{S_G}(r)\times\ZZ/m\ZZ$.

\begin{prop}\label{prop:local_thinning}
    Let $G$ be a finitely generated group and $m\geq 1$. If $G\times\ZZ/m\ZZ$ is finitary space-filling, then so is $G$.
\end{prop}

\begin{proof}
    Let $\Gamma = G\times\ZZ/m\ZZ$ with the set of generators $S'$ defined above, and let $X\subset A^\Gamma$ be a space-filling SFT with movement function $R$ and covering function $f$. Consider the map $\Psi\colon A^\Gamma\to(A^{\ZZ/m\ZZ})^G$ given by $\Psi(x)(g) = (x(g,i))_{i\in\ZZ/m\ZZ}$, and let $Y = \Psi(X)$. The map $\Psi$ is a homeomorphism and satisfies $\Psi((g,0)x) = g\Psi(x)$, so it is a topological conjugacy between the action of $G\times\{0\}$ on $X$ and the shift action of $G$ on $Y$.

    The subgroup $G\times\{0\}$ has finite index in $\Gamma$, with $\{1_G\}\times\ZZ/m\ZZ$ as a set of coset representatives, and $Y$ is the subshift induced by $X$ on this subgroup. As the subaction of an SFT by a subgroup of finite index is an SFT (see for instance~\cite{carroll2015periodic}), $Y$ is a $G$-SFT.

    Next we define the movement function. For $x\in X$, let $\tau(x)=\min\{t\geq 1: R^t_x(1_\Gamma)\in G\times\{0\}\}$. Let $w = R_x(1_\Gamma)$. By~\Cref{prop:faciles} we have $R^{1+j}_x(1_\Gamma) = wR^j_{w^{-1}x}(1_\Gamma)$ for every $j\geq 0$, so applying the definition of the covering function to $w^{-1}x$ shows that the set $R_x^{\llbracket 1,1+f(1)\rrbracket}(1_\Gamma)$ contains a ball of radius $1$ of $\Gamma$, which is of the form $gB_{S_G}(1)\times\ZZ/m\ZZ$ and thus contains $(g,0)$. Hence $\tau(x)\leq f(1)+1$. Let $T\colon Y\to G$ be given by $T(\Psi(x)) = g$, where $R^{\tau(x)}_x(1_\Gamma) = (g,0)$. As $\tau$ is bounded, $T$ takes finitely many values and only depends on the restriction of $x$ to a finite set, so it is continuous.

    Fix $x\in X$ and let $y = \Psi(x)$. Let $0=t_0<t_1<t_2<\cdots$ be the times $t\geq 0$ such that $R^t_x(1_\Gamma)\in G\times\{0\}$, and let $g_n\in G$ be given by $R^{t_n}_x(1_\Gamma) = (g_n,0)$. We claim that $T^n_y(1_G) = g_n$ for every $n\geq 0$. By induction, it suffices to show that $T_y(g_n) = g_{n+1}$. By~\Cref{prop:faciles}, for every $t\geq 0$
    \[ R^t_{(g_n,0)^{-1}x}(1_\Gamma) = (g_n,0)^{-1}R^t_x\bigl((g_n,0)\bigr) = (g_n,0)^{-1}R^{t_n+t}_x(1_\Gamma). \]
    Hence $\tau((g_n,0)^{-1}x) = t_{n+1}-t_n$ and $T(g_n^{-1}y) = T(\Psi((g_n,0)^{-1}x)) = g_n^{-1}g_{n+1}$. We conclude that $T_y(g_n) = g_nT(g_n^{-1}y) = g_{n+1}$.

    Finally, let $r\in\NN$, $r' = \max\{r,1\}$ and $f_G(r) = f(r')$. By definition of the covering function there is $g\in G$ such that $gB_{S_G}(r')\times\ZZ/m\ZZ\subset R_x^{\llbracket 0,f(r')\rrbracket}(1_\Gamma)$. Hence for every $h\in gB_{S_G}(r)$ there is $n$ with $h = g_n$ and $t_n\leq f(r')$. As the $t_n$ are strictly increasing integers, $n\leq t_n\leq f(r')$, and thus $gB_{S_G}(r)\subset T_y^{\llbracket 0,f_G(r)\rrbracket}(1_G)$. As $X$ is nonempty, so is $Y$, and therefore $Y$ is a space-filling $G$-SFT with movement function $T$ and covering function $f_G$.
\end{proof}

\begin{theorem}\label{thm:QI_local}
    Let $G$ and $H$ be quasi-isometric finitely generated groups and suppose that $G$ is finitely presented. Then $G$ is finitary space-filling if and only if $H$ is finitary space-filling.
\end{theorem}

\begin{proof}
    Finite presentability is a quasi-isometry invariant (see for instance~\cite[Chapter I.8, Proposition 8.24]{MR1744486}), so $H$ is finitely presented as well. As being quasi-isometric is a symmetric relation, it suffices to show that if $H$ is finitary space-filling, then so is $G$. Let $\zeta\colon G\to H$ be a quasi-isometry. By~\Cref{lem:thicken} applied to the Cayley graphs $\operatorname{Cay}(G,S_G)$ and $\operatorname{Cay}(H,S_H)$, there are $m\geq 1$ and a surjective quasi-isometry from $\operatorname{Cay}(G,S_G)\boxtimes K_m = \operatorname{Cay}(G\times\ZZ/m\ZZ,S')$ to $\operatorname{Cay}(H,S_H)$, that is, a surjective quasi-isometry $\pi\colon G\times\ZZ/m\ZZ\to H$. As $G\times\ZZ/m\ZZ$ is finitely presented and $H$ is finitary space-filling,~\Cref{prop:local_pullback} shows that $G\times\ZZ/m\ZZ$ is finitary space-filling, and by~\Cref{prop:local_thinning} so is $G$.
\end{proof}

We remark that, as in~\Cref{thm:QI_effective_graphs}, the quasi-isometry is not required to be computable, and the covering functions obtained in~\Cref{prop:local_pullback} and~\Cref{prop:local_thinning} are given by explicit formulas from the original ones. We end this section with a direct corollary.

\begin{corollary}\label{cor:virtab_is_finitary_SF}
      Every finitely generated virtually abelian group is finitary space-filling.
\end{corollary}

\begin{proof}
    On the one hand, for finite groups the result is given by~\Cref{prop:finite_groups}. On the other hand, infinite, finitely generated and virtually abelian groups are finitely presented and quasi-isometric to $\Z^d$ for some $d\geq 1$. By~\Cref{cor:Zd_is_space_filling_locally} we have that $\Z^d$ is finitary space-filling  and thus the result follows from~\Cref{thm:QI_local}.
\end{proof}

\section{Groups quasi-isometric to hyperbolic spaces are space-filling}\label{sec:hyp-group}

In this section we prove that every finitely generated group quasi-isometric to $\HH^{d+1}$ for $d \geq 1$ admits an effective space-filling subshift. To prove this result, it suffices to exhibit, for each $d\geq 1$, a single graph which is quasi-isometric to $\HH^{d+1}$ and show that it is computably space-filling. We use the natural generalization of the dual graph of the binary tiling of the hyperbolic plane to higher dimensions, see~\Cref{fig:grafo2dmegabolico}.

\subsection{The binary \texorpdfstring{$(d+1)$}{(d+1)}-hyperbolic graph}

Fix $d\geq 1$ and write $e_1,\dots,e_d$ for the canonical basis of $\ZZ^d$ and $\mathbf{1}=(1,\dots,1)\in \ZZ^d$. For $(i,j) \in \ZZ\times\ZZ^d$ we call $i \in \ZZ$ a \define{row} and $j\in \ZZ^d$ a \define{horizontal position}.

\begin{defn}
    The \define{binary $(d+1)$-hyperbolic graph} $\mathcal{B}_{d+1}$ is the graph with vertex set $\ZZ\times\ZZ^d$ whose edges are given by
\[ (i,j)\sim (i,j\pm e_k) \mbox{ for }1\leq k\leq d \quad \mbox{and} \quad (i,j)\sim (i+1,2j-\mathbf{1}+\varepsilon)\mbox{ for } \varepsilon\in\{0,1\}^d. \]
\end{defn}

\begin{figure}[ht!]
    \centering
    \begin{tikzpicture}[scale = 0.4]
        \begin{scope}[shift = {(0,-6)},cm={1,0,0.6,0.3,(0,0)}]
            \draw[fill = blue!10,opacity = 0.8] (0,0) rectangle (16,16);
            \draw[step = 1cm] (0,0) grid (16,16);
        \end{scope}
        \begin{scope}[shift = {(0,-6)}]
            \foreach \i in {0,...,8}{
                \foreach \j in {0,...,8}{
                \draw (2*\i+1.2*\j,\j*0.6) to (2*\i+1.2*\j,\j*0.6+6);
            }
        }
        \end{scope} 
        \begin{scope}[cm={1,0,0.6,0.3,(0,0)}]
            \draw[fill = blue!10,opacity = 0.8] (0,0) rectangle (16,16);
            \draw[step = 2cm] (0,0) grid (16,16);
        \end{scope}
        \foreach \i in {0,...,4}{
            \foreach \j in {0,...,4}{
                \draw (4*\i+2.4*\j,\j*1.2) to (4*\i+2.4*\j,\j*1.2+6);
            }
        }
        \begin{scope}[shift = {(0,6)},cm={1,0,0.6,0.3,(0,0)}]
            \draw[fill = blue!10,opacity = 0.8] (0,0) rectangle (16,16);
            \draw[step=4.0cm] (0,0) grid (16,16);
        \end{scope}
        \begin{scope}[shift = {(0,6)}]
            \foreach \i in {0,...,2}{
                \foreach \j in {0,...,2}{
                \draw (8*\i+4.8*\j,\j*2.4) to (8*\i+4.8*\j,\j*2.4+6);
            }
        }
        \end{scope} 
        \begin{scope}[shift = {(0,12)},cm={1,0,0.6,0.3,(0,0)}]
            \draw[fill = blue!10,opacity = 0.8] (0,0) rectangle (16,16);
            \draw[step=8cm] (0,0) grid (16,16);
        \end{scope}
        
    \end{tikzpicture}
    \caption{A portion of the tiles that define $\mathcal{B}_3$. Vertices correspond to cubes, and edges correspond to their common faces. We omit drawing the graph to avoid visual clutter.}
    \label{fig:grafo3dmegabolico}
\end{figure}
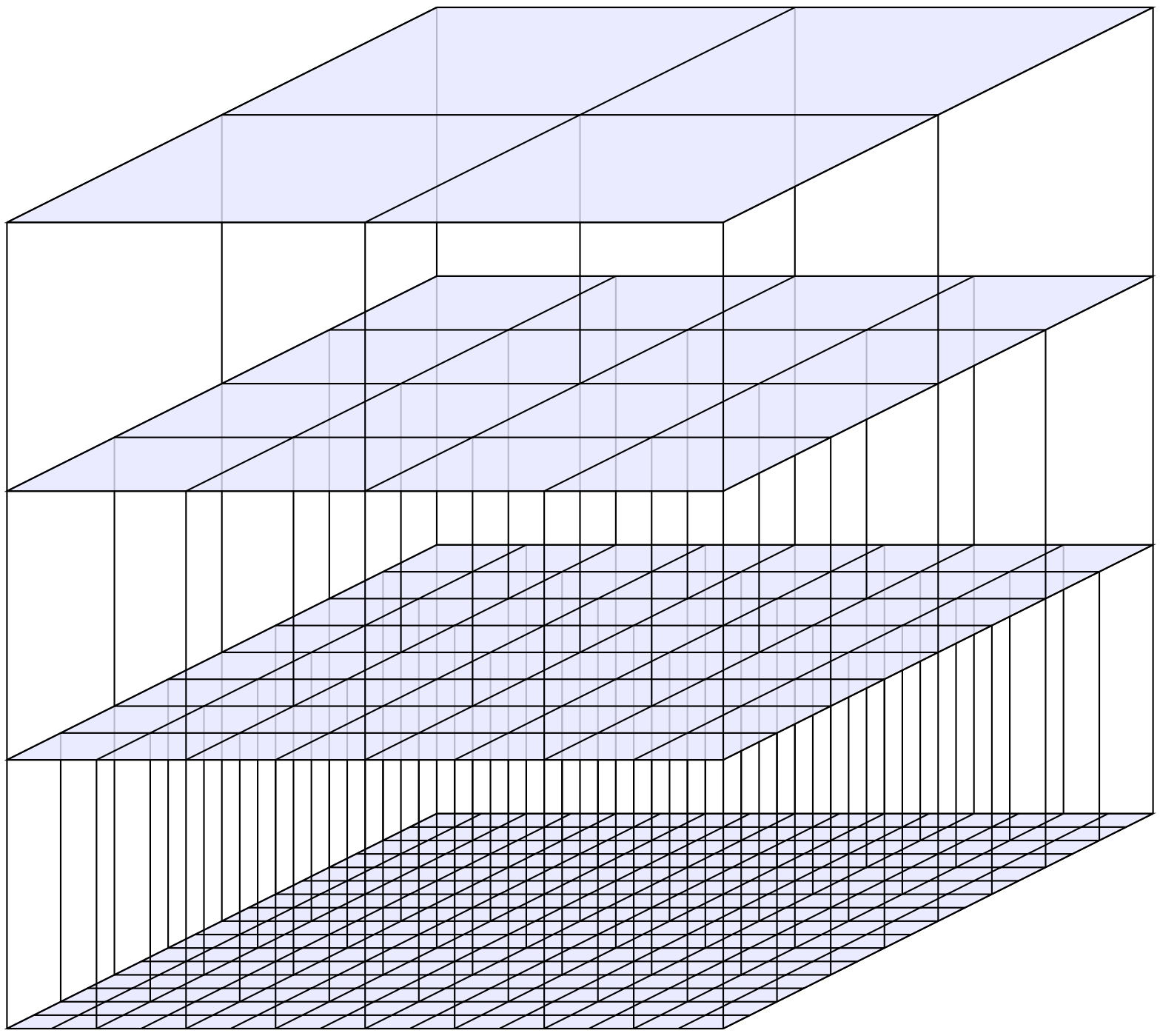

Thus every vertex has $2d$ \define{horizontal neighbors}, $2^d$ \define{children} and one \define{parent}. The cases $d=1$ and $d=2$ are shown in \Cref{fig:grafo2dmegabolico} and \Cref{fig:grafo3dmegabolico} respectively. We endow this graph with the path metric $d_{\mathcal{B}_{d+1}}$. Identifying $\ZZ\times\ZZ^d$ with $\NN$ through a computable bijection, the graphs $\mathcal{B}_{d+1}$ are highly computable, since adjacency is decidable and every vertex has degree $2d+2^d+1$.

\begin{figure}[ht!]
    \centering
    \begin{tikzpicture}[scale =0.8]

	\begin{scope}[scale=0.3, shift={(0,0)}, opacity=0.5]
    \foreach \xx in {-12,0,12,24}{
	\begin{scope}[shift={(\xx,3)}, scale = 2]
	\draw (0,0) -- (0,3);
	\draw (0,3) to [controls=+(45:2) and +(135:2)] (6,3);
	\draw (6,3) -- (6,0);
	\draw (0,0) to [controls=+(45:1) and +(135:1)] (3,0);
	\draw (3,0) to [controls=+(45:1) and +(135:1)] (6,0);
	\end{scope}
	}
	\foreach \xx in {-12,-6,...,30}{
	\begin{scope}[shift={(\xx,0)}]
	\draw (0,0) -- (0,3);
	\draw (0,3) to [controls=+(45:2) and +(135:2)] (6,3);
	\draw (6,3) -- (6,0);
	\draw (0,0) to [controls=+(45:1) and +(135:1)] (3,0);
	\draw (3,0) to [controls=+(45:1) and +(135:1)] (6,0);
	\end{scope}
	}
	\foreach \xx in {-12,-9,...,33}{
	\begin{scope}[shift={(\xx,-3)},scale=0.5]
	\draw (0,0) -- (0,6);
	\draw (0,6) to [controls=+(45:2) and +(135:2)] (6,6);
	\draw (6,6) -- (6,0);
	\draw (0,0) to [controls=+(45:1) and +(135:1)] (3,0);
	\draw (3,0) to [controls=+(45:1) and +(135:1)] (6,0);
	\end{scope}
	}
	\foreach \xx in {-12,-10.5,...,34.5}{
	\begin{scope}[shift={(\xx,-6)},scale=0.25]
	\draw (0,0) -- (0,12);
	\draw (0,12) to [controls=+(45:2) and +(135:2)] (6,12);
	\draw (6,12) -- (6,0);
	\draw (0,0) to [controls=+(45:1) and +(135:1)] (3,0);
	\draw (3,0) to [controls=+(45:1) and +(135:1)] (6,0);
	\end{scope}
	}
    \foreach \xx in {-12,-11.25,...,35.25}{
	\begin{scope}[shift={(\xx,-7.5)},scale=0.125]
	\draw (0,0) -- (0,12);
	\draw (0,12) to [controls=+(45:2) and +(135:2)] (6,12);
	\draw (6,12) -- (6,0);
	\draw (0,0) to [controls=+(45:1) and +(135:1)] (3,0);
	\draw (3,0) to [controls=+(45:1) and +(135:1)] (6,0);
	\end{scope}
	}
	\end{scope}
     \draw[blue] (-1.8,2) -- (9,2);
    \foreach \i in {0,3.6,7.2,10.8}{
    \draw[blue]  (-1.8+\i,2) -- (-1.8+\i-0.9,0.5);
    \draw[blue]  (-1.8+\i,2) -- (-1.8+\i+0.9,0.5);
    \filldraw (-1.8+\i,2) circle (0.1);
    }
    \draw[blue] (-2.7,0.5) -- (9.9,0.5);
    \foreach \i in {0,1.8,...,13}{
     \draw[blue]  (-2.7+\i,0.5) -- (-2.7+\i+0.45,-0.5);
     \draw[blue]  (-2.7+\i,0.5) -- (-2.7+\i-0.45,-0.5);
    \filldraw (-2.7+\i,0.5) circle (0.06);
    }
    \draw[blue] (-3.15,-0.5) -- (10.35,-0.5);
    \foreach \i in {0,0.9,...,14}{
     \draw[blue]  (-3.15+\i,-0.5) -- (-3.15+\i+0.225,-1.5);
     \draw[blue]  (-3.15+\i,-0.5) -- (-3.15+\i-0.225,-1.5);
    \filldraw (-3.15+\i,-0.5) circle (0.06);
    }
    \draw[blue] (-3.375,-1.5) -- (10.575,-1.5);
    \foreach \i in {0,0.45,...,14}{
     \draw[blue]  (-3.375+\i,-1.5) -- (-3.375+\i+0.1125,-2);
     \draw[blue]  (-3.375+\i,-1.5) -- (-3.375+\i-0.1125,-2);
    \filldraw (-3.375+\i,-1.5) circle (0.04);
    }
    \draw[blue] (-3.4875,-2) -- (10.575+0.1125,-2);
     \foreach \i in {0,0.225,...,14.2}{
     \filldraw (-3.4875+\i,-2) circle (0.03);
    }

	\end{tikzpicture}

    \caption{A portion of the graph $\mathcal{B}_2$. The tiles correspond to vertices and edges are between adjacent tiles.}
    \label{fig:grafo2dmegabolico}
\end{figure}
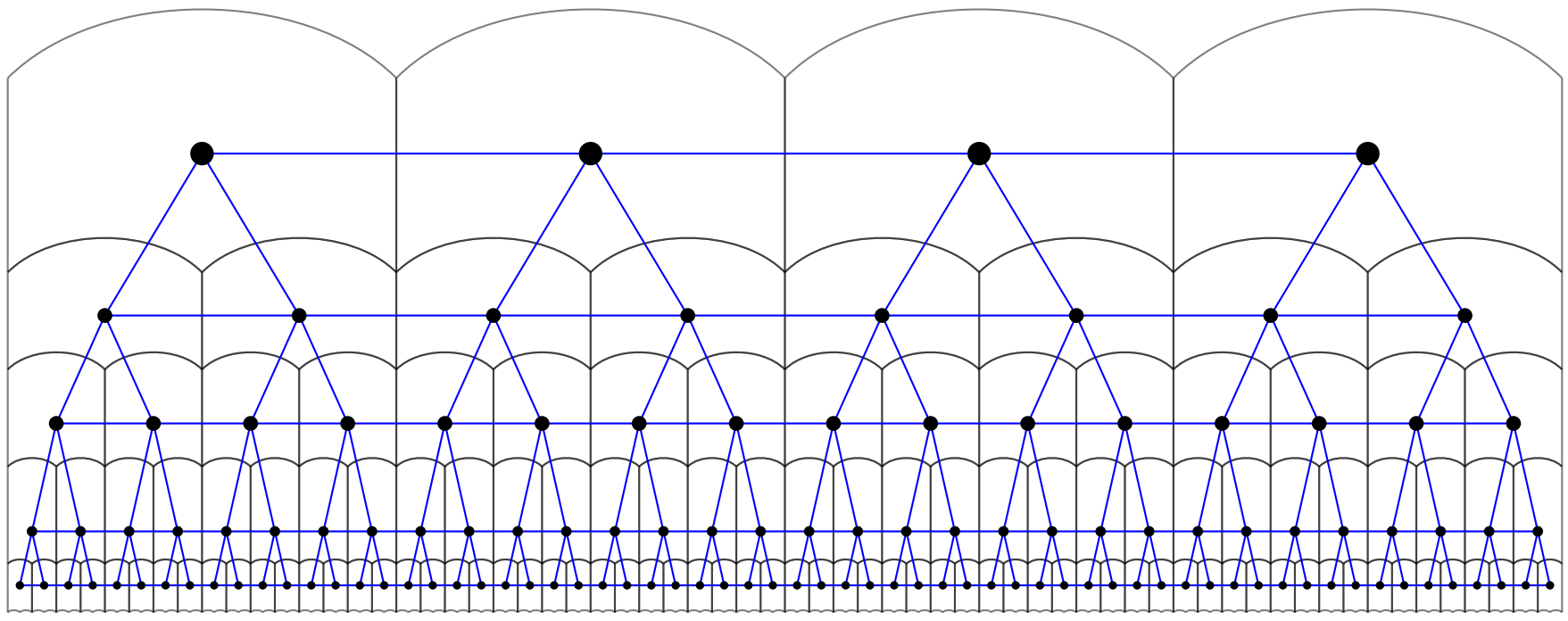

These graphs are naturally quasi-isometric to the $(d+1)$-dimensional hyperbolic space $\HH^{d+1}$ seen through its upper half-space model \[\HH^{d+1} = \{ (x_1,\dots,x_d,x_{d+1}) \in \RR^{d+1} : x_{d+1}>0\}.\]  Indeed, the vertex $(i,j)$ can be thought of as representing the \textbf{tile}
\[t(i,j)=t(i,j_1,\dots,j_d) = \left(\prod_{k=1}^d[(j_k-1)2^{-i},j_k2^{-i}] \right)\times[2^{-i},2^{1-i}],\] 
and two vertices are adjacent if and only if their tiles share a face of positive $d$-dimensional measure. These tiles are congruent, as $z\mapsto 2^{-i}\bigl(z+(j-\mathbf{1},0)\bigr)$ maps the tile of $(0,\mathbf{1})$ onto the tile of $(i,j)$, and intersecting tiles correspond to vertices at distance at most $d+1$ in $\mathcal{B}_{d+1}$: they lie in the same or in consecutive rows, and in the latter case the horizontal position of one of them differs by at most one in each coordinate from that of a child of the other.
 
\begin{lemma}\label{lem:Bd_QI}
    For every $d\geq 1$, the graph $\mathcal{B}_{d+1}$ is quasi-isometric to $\HH^{d+1}$.
\end{lemma}

\begin{proof}
    Fix any point $z_0\in t(0,\mathbf{1})$. Consider the function $\zeta:\mathcal{B}_{d+1}\to \HH^{d+1}$ given by
    \[\zeta(i,j) = 2^{-i}\bigl(z_0 +(j-\mathbf{1},0)\bigr).\]
    Let $\delta>0$ be the diameter of the tiles $t(i,j)$ with $(i,j)\in \ZZ \times \ZZ^d$, and consider $\rho>0$ the radius of the largest ball contained in $t(0,\mathbf{1})$. Then, as the tiles $t(i,j)$ cover $\HH^{d+1}$, every point in $\HH^{d+1}$ is at distance at most $\delta$ from an element of $\zeta(\ZZ\times\ZZ^d)$.
    
    If two vertices $u,v$ are at distance $n = d_{\mathcal{B}_{d+1}}(u,v)$, there is a sequence of $n+1$ tiles from $t(u)$ to $t(v)$ in which consecutive tiles share a face. As $\zeta(w)\in t(w)$ for every vertex $w$, it follows that $d_{\HH}(\zeta(u),\zeta(v))\leq 2\delta n$. Now, notice that the tiles which intersect a given ball of radius $1$ are all contained in the concentric ball of radius $1+\delta$. Furthermore, each of those tiles contains a ball of radius $\rho$, disjoint from the rest. Therefore, there are at most $K = \operatorname{vol}(B(1+\delta))/\operatorname{vol}(B(\rho))$ of them.
     Finally, take $u,v$ vertices and consider a geodesic $\gamma$ from $\zeta(u)$ to $\zeta(v)$ of length $\lambda$. Then $\gamma$ can be covered by at most $\frac{\lambda}{2}+1$ balls of radius $1$, so at most $K\left(\frac{\lambda}{2} + 1\right)$ tiles intersect $\gamma$. These tiles contain $t(u)$ and $t(v)$, and as $\gamma$ is connected and covered by finitely many closed tiles, any two of them are joined by a chain of tiles meeting $\gamma$ in which consecutive tiles intersect. As intersecting tiles correspond to vertices at distance at most $d+1$, we obtain \[d_{\mathcal{B}_{d+1}}(u,v)\leq (d+1)K\left(\frac{1}{2}d_{\HH}(\zeta(u),\zeta(v)) + 1\right).\] 
     Thus, $\zeta$ is a quasi-isometry.
\end{proof}

By~\Cref{thm:QI_effective_graphs} and~\Cref{prop:effective_SF}, and since groups quasi-isometric to $\HH^{d+1}$ have decidable word problem (for instance by~\cite[Chapter III.$\Gamma$, Theorem 2.6]{MR1744486}), in order to show that every finitely generated group which is quasi-isometric to $\HH^{d+1}$ is computably space-filling, it will suffice to show that the graphs $\mathcal{B}_{d+1}$ are computably space-filling.

\subsection{A hierarchical structure on \texorpdfstring{$\mathcal{B}_{d+1}$}{B(d+1)}}

In each of the graphs $\mathcal{B}_{d+1}$ we will define a sequence of finite subsets of vertices $(F_{\ell})_{\ell \geq 0}$ that we shall refer to as \define{fans} and that will satisfy three fundamental properties (along with some minor technical ones): first, each $F_{\ell+1}$ can be decomposed into disjoint copies of $F_{\ell}$; second, each $F_{\ell}$ can be traversed by a Hamiltonian path in such a way that the lower level sets it is composed of are traversed completely in succession; third, for every radius $r\geq 0$ there is some level $\ell(r) \geq 0$ for which $F_{\ell(r)}$ contains a ball of radius $r$.

With this sequence in mind, a path for which every window sees a ball of radius $r$ can be obtained by just following the Hamiltonian path in $F_{\ell}$ for $\ell \geq \ell(r)+1$ and taking a window of size $2|F_{\ell(r)}|$, which always contains a complete copy of $F_{\ell(r)}$, as these copies are traversed in consecutive blocks of $|F_{\ell(r)}|$ vertices.

\begin{defn}
    The \define{fan} of level $\ell\geq 0$ is
\[ F_\ell = \{(i,j_1,\dots,j_d): 1\leq i\leq 2^\ell,\ 1\leq j_k\leq 2^{i-1} \textnormal{ for } 1\leq k\leq d\}. \]
\end{defn}

\begin{figure}
    \centering
    \begin{tikzpicture}[scale =0.8]

	\begin{scope}[scale=0.3, shift={(0,0)}, opacity = 0.5]
    \foreach \xx in {-12}{
	\begin{scope}[shift={(\xx,3)}, scale = 2]
	\draw (0,0) -- (0,3);
	\draw (0,3) to [controls=+(45:2) and +(135:2)] (6,3);
	\draw (6,3) -- (6,0);
	\draw (0,0) to [controls=+(45:1) and +(135:1)] (3,0);
	\draw (3,0) to [controls=+(45:1) and +(135:1)] (6,0);
	\end{scope}
	}
	\foreach \xx in {-12,-6}{
	\begin{scope}[shift={(\xx,0)}]
	\draw (0,0) -- (0,3);
	\draw (0,3) to [controls=+(45:2) and +(135:2)] (6,3);
	\draw (6,3) -- (6,0);
	\draw (0,0) to [controls=+(45:1) and +(135:1)] (3,0);
	\draw (3,0) to [controls=+(45:1) and +(135:1)] (6,0);
	\end{scope}
	}
	\foreach \xx in {-12,-9,...,-3}{
	\begin{scope}[shift={(\xx,-3)},scale=0.5]
	\draw (0,0) -- (0,6);
	\draw (0,6) to [controls=+(45:2) and +(135:2)] (6,6);
	\draw (6,6) -- (6,0);
	\draw (0,0) to [controls=+(45:1) and +(135:1)] (3,0);
	\draw (3,0) to [controls=+(45:1) and +(135:1)] (6,0);
	\end{scope}
	}
	\foreach \xx in {-12,-10.5,...,-1.5}{
	\begin{scope}[shift={(\xx,-6)},scale=0.25]
	\draw (0,0) -- (0,12);
	\draw (0,12) to [controls=+(45:2) and +(135:2)] (6,12);
	\draw (6,12) -- (6,0);
	\draw (0,0) to [controls=+(45:1) and +(135:1)] (3,0);
	\draw (3,0) to [controls=+(45:1) and +(135:1)] (6,0);
	\end{scope}
	}
	\end{scope}
    \filldraw[green!50!black, opacity = 0.2] (-1.8,2.4) -- (-3,0.3) -- (-0.6,0.3);
    \filldraw (-1.8,2) circle (0.1);
    \draw[blue] (-2.7,0.5) -- (-1.8,2);
    \draw[blue] (-0.9,0.5) -- (-1.8,2);
    \draw[blue] (-0.9,0.5) -- (-2.7,0.5);
    \foreach \i in {0,1.8}{
    \draw[blue] (-0.9-\i,0.5) -- (-1.35-\i,-0.5);
    \draw[blue] (-0.9-\i,0.5) -- (-0.45-\i,-0.5);
    \filldraw (-0.9-\i,0.5) circle (0.06);
    
    }
    \foreach \i in {0,0.9,1.8,2.7}{
    \filldraw[green!50!black, opacity = 0.2] (-0.45-\i,-0.2) -- (-0.1-\i,-1.6) -- (-0.8-\i,-1.6);
    }
    \draw[blue] (-0.45,-0.5) -- (-3.15,-0.5);
    \foreach \i in {0,0.9,1.8,2.7}{
    \draw[blue] (-0.45-\i,-0.5) -- (-0.225-\i,-1.5);
    \draw[blue] (-0.45-\i,-0.5) -- (-0.675-\i,-1.5);
    \filldraw (-0.45-\i,-0.5) circle (0.06);
    }
    \draw[blue] (-0.225,-1.5) -- (-3.375,-1.5);
    \foreach \i in {0,0.45,...,3.15}{
    \filldraw (-0.225-\i,-1.5) circle (0.04);
    }
	\end{tikzpicture}

    \caption{A fan $F_{\ell}$ of level $\ell =2$ on $\mathcal{B}_2$. It is composed of $5$ fans of level $1$ marked in green.}
    \label{fig:onlyfans}
\end{figure}
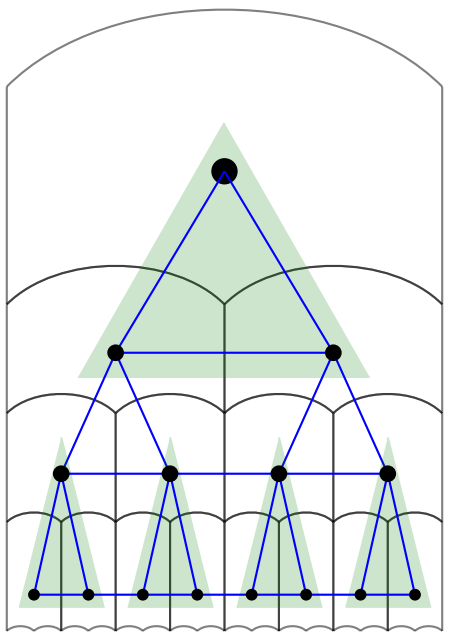

The name \define{fans} comes from what these objects look like for $d=1$, see~\Cref{fig:onlyfans}. The name for $d>1$ is less appropriate but we shall use it anyway.

Fans of level $\ell+1$ admit a natural decomposition into copies of fans of level $\ell$. Namely, for $(a,b)\in\ZZ\times\Z^d$, consider the map $\varphi_{a,b}$ defined on the vertices $(i,j)\in \ZZ\times \ZZ^d$ with $i\geq 1$ by
\[ \varphi_{a,b}(i,j) = (a+i,\ b2^{i-1}+j). \]
This map is injective and maps adjacent vertices to adjacent vertices. Moreover, if $i\geq2$ then $\varphi_{a,b}$ maps the set of neighbors of $(i,j)$ bijectively onto the set of neighbors of $\varphi_{a,b}(i,j)$. For horizontal neighbors and children this is clear, and the parent $(i-1,\lceil j/2\rceil)$ is mapped to $(a+i-1,b2^{i-2}+\lceil j/2\rceil)$, which is the parent of $\varphi_{a,b}(i,j)$ because $b2^{i-1}$ is even. In particular, if every vertex of $B_{\mathcal{B}_{d+1}}(v,r-1)$ lies in a row $\geq 2$, then $\varphi_{a,b}(B_{\mathcal{B}_{d+1}}(v,r)) = B_{\mathcal{B}_{d+1}}(\varphi_{a,b}(v),r)$. Note also that $\varphi_{a,b}\circ\varphi_{a',b'} = \varphi_{a+a',2^{a'}b+b'}$, so the copies of $F_{\ell'}$ appearing at any level of the hierarchy of $F_\ell$ are all of the form $\varphi_{a,b}(F_{\ell'})$.

Let $N_\ell = 2^{2^\ell}$ and note that
\begin{equation*}
    F_{\ell+1} = F_\ell\ \sqcup\ \bigsqcup_{m\in\{1,\dots,N_{\ell}\}^d}\varphi_{2^\ell,m-\mathbf{1}}(F_\ell).
\end{equation*}

Hence each $F_{\ell+1}$ is naturally decomposed into $1+N_{\ell}^d$ copies of $F_{\ell}$. The following lemma shows that fans contain arbitrarily large balls. While this is quite natural, the proof requires some care.

\begin{lemma}\label{lem:balls_fans}
    For every $r\geq 1$ there exist $\ell(r)\geq 0$ and $v\in \ZZ\times\ZZ^d$ such that $B_{\mathcal{B}_{d+1}}(v,r)\subseteq F_{\ell(r)}$ and every vertex of $B_{\mathcal{B}_{d+1}}(v,r)$ lies in a row $i\geq 3$. Moreover, the value $\ell(r)$ is computable from $r$, and for every $a,b$ the copy $\varphi_{a,b}(F_{\ell(r)})$ contains the ball $B_{\mathcal{B}_{d+1}}(\varphi_{a,b}(v),r)$.
\end{lemma}

\begin{proof}
    To every vertex $(i,j)\in F=\bigcup_{\ell \geq 0}F_{\ell}$ we assign a $d$-tuple of words $(w_1,\dots,w_d)\in (\{0,1\}^*)^d$ which correspond to their symbolic coordinates in the product of $d$ copies of the binary tree, where $i-1$ is the depth. Namely, for $(i,j_1,\dots,j_d)$ we assign the tuple where for each $k \in \{1,\dots,d\}$, $|w_k|=i-1$ and $w_k$ is the binary expansion of $j_k-1$, padded with leading zeros to length $i-1$.

    Now consider $h \geq r+1$ and take the vertex $v=(h+2,(2^h+1)\cdot \mathbf{1})$. This vertex clearly belongs to $F$ and its symbolic coordinates are $(w_1,\dots,w_d)$ where all words are given by $w_k = 10^{h}$. Now we note that moving towards the parent corresponds to removing the last symbol in each word in $(w_1,\dots,w_d)$, moving horizontally corresponds to adding or subtracting $1$ as an odometer in a single coordinate, and moving toward a child corresponds to appending a symbol to each word. 
    
    Let us show that no combination of $r$ operations can move beyond $F_{\ell(r)}$ starting from $v$, where $\ell(r)$ is such that $h+2 \leq 2^{\ell(r)-1}$. First, the row stays in $\llbracket h+2-r, h+2+r\rrbracket\subseteq\llbracket 3, 2^{\ell(r)}\rrbracket$, as $r<h+2\leq 2^{\ell(r)-1}$. Next, fix $k$ and let $a\in\llbracket 0,2^L-1\rrbracket$ be the integer with binary expansion $w_k$, where $L=|w_k|$. A path leaves $F$ only if some horizontal move turns $a=0$ into $-1$ or $a=2^L-1$ into $2^L$. Both $a$ and $2^L-1-a$ are replaced by their integer half when moving to the parent, by twice their value plus $0$ or $1$ when moving to a child, and change by at most one under a horizontal move. Hence the number of binary digits of each of them decreases by at most one at each step. Initially $a=2^h$ and $2^L-1-a = 2^h-1$ have $h+1$ and $h$ binary digits respectively, so at least $h\geq r+1$ steps are needed before either of them vanishes, and therefore $r$ steps never leave $F$. Finally, as every vertex of $B_{\mathcal{B}_{d+1}}(v,r)$ lies in a row $\geq 3$, the discussion following the definition of $\varphi_{a,b}$ shows that $\varphi_{a,b}(F_{\ell(r)})$ contains $B_{\mathcal{B}_{d+1}}(\varphi_{a,b}(v),r)$ for every $a,b$. The fact that $\ell(r)$ is computable from $r$ is direct.
    \end{proof}

Given the previous lemma, the only remaining step is to find Hamiltonian paths on each fan that are consistent with the hierarchical structure of the $(F_{\ell})_{\ell \geq 0}$. Recall that a \define{Hamiltonian path} on a set $V$ is a sequence of distinct vertices $v_1,\dots, v_{|V|}$ such that $v_i$ is adjacent to $v_{i+1}$ for all $i\in\{1,\dots,|V|-1\}$. 

\begin{defn}
    A Hamiltonian path $v_1,\dots, v_{|F_\ell|}$ on $F_\ell$ is \textbf{hereditary Hamiltonian} if either $\ell = 1$, or $\ell \geq 2$ and, writing $n = |F_{\ell-1}|$, for every $k\in\{0,\dots, |F_\ell|/n-1\}$ the block $v_{kn+1}, \dots, v_{(k+1)n}$ is the image under some $\varphi_{a,b}$ of a hereditary Hamiltonian path on $F_{\ell-1}$. 
    
    The collection $(F_{\ell})_{\ell \geq 0}$ is \define{hierarchically Hamiltonian} if every $F_{\ell}$ for $\ell \geq 1$ admits a hereditary Hamiltonian path. 
\end{defn}

Our goal will be to show that $(F_{\ell})_{\ell \geq 0}$ is hierarchically Hamiltonian. To that end, we define a \textbf{port} of $F_\ell$ as an ordered pair of vertices $(u,v)$ such that there exists a hereditary Hamiltonian path $v_1,\dots, v_{|F_\ell|}$ with $u = v_1$ and $v = v_{|F_\ell|}$. Notice that, by definition, proving $(F_{\ell})_{\ell \geq 0}$ is hierarchically Hamiltonian is equivalent to the fact that the set of ports of $F_{\ell}$ is nonempty for every $\ell \geq 1$. We will do this by inductively using ports of level $\ell$ to construct a set of ports in level $\ell+1$. The strategy to do this is very different for the $d=1$ case and the $d\geq 2$ case, so we treat them separately in the next two subsections.

Before going into the proofs of both cases, we remark that the set of ports is symmetric, that is, if $(u,v)$ is a port of $F_{\ell}$, then so is $(v,u)$, since the reverse of a hereditary Hamiltonian path is hereditary Hamiltonian.

\subsection{The case \texorpdfstring{$d=1$}{d=1}}

We begin by giving names to the following three vertices of $F_\ell$, which are pairwise distinct for $\ell\geq 1$:
\[ \uparrow = (1,1), \qquad \swarrow = (2^\ell,1), \qquad \searrow = (2^\ell,2^{2^\ell-1}). \]
Recall that $N_\ell = 2^{2^\ell}$ is the number of vertices of $F_{\ell+1}$ in row $2^\ell+1$. Also recall that we decompose each level $\ell\geq 0$ into the previous levels as
\begin{equation*}
    F_{\ell+1} = F_\ell \ \sqcup\ \bigsqcup_{m=1}^{N_\ell}\varphi_{2^\ell,m-1}(F_\ell).
\end{equation*}

We call $T=F_\ell$ the \define{top piece} and $B_m = \varphi_{2^\ell,m-1}(F_\ell)$ the \define{$m$-th bottom piece} of $F_{\ell+1}$. The bottom piece $B_m$ is the fan of height $2^\ell$ hanging from its \define{root} $\varphi_{2^\ell,m-1}(\uparrow)=(2^\ell+1,m)$, and it is the image of $F_\ell$ under a map which preserves adjacency. We denote the images of the distinguished vertices of $F_\ell$ in $B_m$ by $\uparrow_{B_m}$, $\swarrow_{B_m}$ and $\searrow_{B_m}$, and similarly for $T$.

Finally, we also need an additional pair which we call the \define{seam} of $F_\ell$. It is given by the pair
\[ \Pi_{\ell} = \left((2^\ell, 2^{2^\ell -2}),\ (2^\ell,2^{2^\ell -2}+1)\right) \mbox{ for } \ell\geq 1,\]
which consists of the vertex at the center of the bottom row of $F_{\ell}$ and its right neighbor. Note that for $\ell=1$ we have $(\swarrow,\searrow) = \Pi_1$.
 
\begin{prop}\label{prop:fan_ports}
    For every $\ell\geq 1$, the pairs $(\uparrow,\swarrow)$, $(\uparrow,\searrow)$, $(\swarrow,\searrow)$ and $\Pi_{\ell}$ are ports of $F_\ell$.
\end{prop}
 
\begin{proof}
    We proceed by induction on $\ell$. For $\ell=1$ we have $F_1 = \{(1,1),(2,1),(2,2)\}$. The Hamiltonian paths
    \[ (1,1),(2,2),(2,1); \qquad (1,1),(2,1),(2,2); \qquad (2,1),(1,1),(2,2) \]
    realize $(\uparrow,\swarrow)$, $(\uparrow,\searrow)$ and $(\swarrow,\searrow) = \Pi_1$.
 
    Now suppose the statement holds for some $\ell\geq 1$, and consider the decomposition of $F_{\ell+1}$ into $T=F_\ell$ and the bottom pieces $B_1,\dots,B_{N_{\ell}}$. We consider two ways of connecting ports from these pieces, which we call \define{junctions}. Note that the bottom row of each piece has $t_{\ell} = N_{\ell}/2$ vertices. The two junctions are
    
    \begin{enumerate}
        \item For $1\leq m<N_\ell$, the vertices $\searrow_{B_m}=(2^{\ell+1},mt_{\ell})$ and $\swarrow_{B_{m+1}}=(2^{\ell+1},mt_{\ell}+1)$ are horizontal neighbors.
        \item For $1\leq m\leq N_\ell$, the root $\uparrow_{B_m}=(2^\ell+1,m)$ is a child of the vertex $(2^\ell,\lceil m/2\rceil)$ in the bottom row of $T$. In particular $\swarrow_T = (2^\ell,1)$ is the parent of $\uparrow_{B_1}$, and $\searrow_T=(2^\ell,t_{\ell})$ is the parent of $\uparrow_{B_{N_\ell}}$.
    \end{enumerate}
    Moreover, the seam of $T$ is $\Pi_{\ell} = \left((2^\ell, t_{\ell}/2),\ (2^\ell,1+t_{\ell}/2)\right)$. Note that using the second junction we see that the vertex $(2^\ell, t_{\ell}/2)$ is the parent of $\uparrow_{B_{t_{\ell}}}$ and $(2^\ell,1+t_{\ell}/2)$ is the parent of $\uparrow_{B_{t_{\ell}+1}}$.

    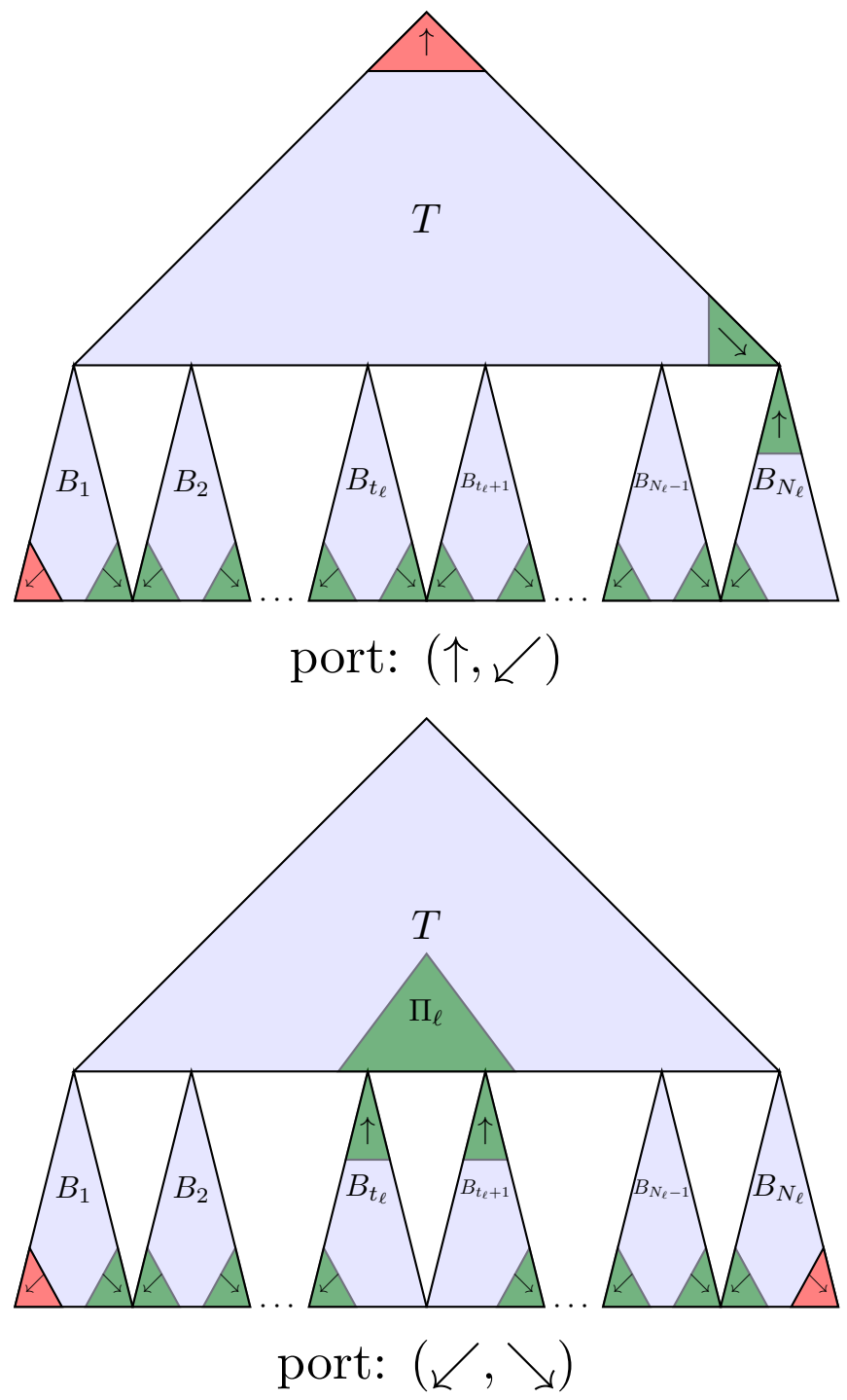
\begin{figure}[ht!]
        \centering
        \begin{tikzpicture}[scale = 0.4]
         \begin{scope}[shift = {(0,0)}]
                \draw[fill = blue!10] (-2,0) -- (10,0) -- (4,6) -- cycle;
            \draw[fill = blue!10] (-3,-4) -- (-1,-4) -- (-2,0) -- cycle;
            \draw[fill = blue!10] (-1,-4) -- (1,-4) -- (0,0) -- cycle;
            \draw[fill = blue!10] (2,-4) -- (4,-4) -- (3,0) -- cycle;
            \draw[fill = blue!10] (4,-4) -- (6,-4) -- (5,0) -- cycle;
            \draw[fill = blue!10] (7,-4) -- (9,-4) -- (8,0) -- cycle;
            \draw[fill = blue!10] (9,-4) -- (11,-4) -- (10,0) -- cycle;
            \node at (1.5,-4){\scalebox{0.6}{$\cdots$}};
            \node at (6.5,-4) {\scalebox{0.6}{$\cdots$}};
            \node at (4,-5) {port: $(\uparrow,\swarrow)$};
            \draw[fill = red!50] (3,5) -- (5,5) -- (4,6) -- cycle;
            \node at (4,5.5)  {\scalebox{0.6}{$\uparrow$}};
            \draw[fill = green!50!black, opacity = 0.5] (10,0) -- (8.8,0) -- (8.8,1.2) -- cycle;
            \node at (9.2,0.4)  {\scalebox{0.6}{$\searrow$}};
            \draw[fill = green!50!black, opacity = 0.5] (10-0.375,-1.5) -- (10+0.375,-1.5) -- (10,0) -- cycle;
            \node at (10,-1)  {\scalebox{0.6}{$\uparrow$}};
            \foreach \i in {-1,1,4,6,9}{
             \draw[fill = green!50!black, opacity = 0.5] (\i-0.25,-3) -- (\i-0.8,-4) -- (\i,-4) -- cycle;
            \node at (\i-0.35,-3.6)  {\scalebox{0.4}{$\searrow$}};
            }
            \foreach \j in {-1,2,4,7,9}{
             \draw[fill = green!50!black, opacity = 0.5] (\j+0.25,-3) -- (\j+0.8,-4) -- (\j,-4) -- cycle;
            \node at (\j+0.35,-3.6)  {\scalebox{0.4}{$\swarrow$}};
            }
            \draw[fill = red!50] (-3+0.25,-3) -- (-3+0.8,-4) -- (-3,-4) -- cycle;
            \node at (-3+0.35,-3.6)  {\scalebox{0.4}{$\swarrow$}};
            \node at (4,2.5) {\scalebox{0.8}{$T$}};
            \node at (-2,-2) {\scalebox{0.6}{$B_1$}};
            \node at (0,-2) {\scalebox{0.6}{$B_2$}};
            \node at (3,-2) {\scalebox{0.6}{$B_{t_{\ell}}$}};
            \node at (5,-2) {\scalebox{0.4}{$B_{t_{\ell}+1}$}};
            \node at (8,-2) {\scalebox{0.4}{$B_{N_{\ell}-1}$}};
            \node at (10,-2) {\scalebox{0.6}{$B_{N_{\ell}}$}};

            \end{scope}

             \begin{scope}[shift = {(16,0)}]
            \draw[fill = blue!10] (-2,0) -- (10,0) -- (4,6) -- cycle;
            \draw[fill = blue!10] (-3,-4) -- (-1,-4) -- (-2,0) -- cycle;
            \draw[fill = blue!10] (-1,-4) -- (1,-4) -- (0,0) -- cycle;
            \draw[fill = blue!10] (2,-4) -- (4,-4) -- (3,0) -- cycle;
            \draw[fill = blue!10] (4,-4) -- (6,-4) -- (5,0) -- cycle;
            \draw[fill = blue!10] (7,-4) -- (9,-4) -- (8,0) -- cycle;
            \draw[fill = blue!10] (9,-4) -- (11,-4) -- (10,0) -- cycle;
            \node at (1.5,-4) {\scalebox{0.6}{$\cdots$}};
            \node at (6.5,-4) {\scalebox{0.6}{$\cdots$}};
            \node at (4,-5) {port: $(\uparrow,\searrow)$}; 
            \draw[fill = red!50] (3,5) -- (5,5) -- (4,6) -- cycle;
            \node at (4,5.5)  {\scalebox{0.6}{$\uparrow$}};
            \draw[fill = green!50!black, opacity = 0.5] (-2,0) -- (-0.8,0) -- (-0.8,1.2) -- cycle;
            \node at (-1.2,0.4)  {\scalebox{0.6}{$\swarrow$}};
             \draw[fill = green!50!black, opacity = 0.5] (-2-0.375,-1.5) -- (-2+0.375,-1.5) -- (-2,0) -- cycle;
            \node at (-2,-1)  {\scalebox{0.6}{$\uparrow$}};

            \foreach \i in {-1,1,4,6,9}{
             \draw[fill = green!50!black, opacity = 0.5] (\i-0.25,-3) -- (\i-0.8,-4) -- (\i,-4) -- cycle;
            \node at (\i-0.35,-3.6)  {\scalebox{0.4}{$\searrow$}};
            }
            \foreach \j in {-1,2,4,7,9}{
             \draw[fill = green!50!black, opacity = 0.5] (\j+0.25,-3) -- (\j+0.8,-4) -- (\j,-4) -- cycle;
            \node at (\j+0.35,-3.6)  {\scalebox{0.4}{$\swarrow$}};
            }

            \draw[fill = red!50] (11-0.25,-3) -- (11-0.8,-4) -- (11,-4) -- cycle;
            \node at (11-0.35,-3.6)  {\scalebox{0.4}{$\searrow$}};

            \node at (4,2.5) {\scalebox{0.8}{$T$}};
            \node at (-2,-2) {\scalebox{0.6}{$B_1$}};
            \node at (0,-2) {\scalebox{0.6}{$B_2$}};
            \node at (3,-2) {\scalebox{0.6}{$B_{t_{\ell}}$}};
            \node at (5,-2) {\scalebox{0.4}{$B_{t_{\ell}+1}$}};
            \node at (8,-2) {\scalebox{0.4}{$B_{N_{\ell}-1}$}};
            \node at (10,-2) {\scalebox{0.6}{$B_{N_{\ell}}$}};
            
            \end{scope}

            \begin{scope}[shift = {(0,-12)}]
                \draw[fill = blue!10] (-2,0) -- (10,0) -- (4,6) -- cycle;
            \draw[fill = blue!10] (-3,-4) -- (-1,-4) -- (-2,0) -- cycle;
            \draw[fill = blue!10] (-1,-4) -- (1,-4) -- (0,0) -- cycle;
            \draw[fill = blue!10] (2,-4) -- (4,-4) -- (3,0) -- cycle;
            \draw[fill = blue!10] (4,-4) -- (6,-4) -- (5,0) -- cycle;
            \draw[fill = blue!10] (7,-4) -- (9,-4) -- (8,0) -- cycle;
            \draw[fill = blue!10] (9,-4) -- (11,-4) -- (10,0) -- cycle;
            \node at (1.5,-4){\scalebox{0.6}{$\cdots$}};
            \node at (6.5,-4) {\scalebox{0.6}{$\cdots$}};
            \node at (4,-5) {port: $(\swarrow,\searrow)$};

            \foreach \i in {-1,1,6,9}{
             \draw[fill = green!50!black, opacity = 0.5] (\i-0.25,-3) -- (\i-0.8,-4) -- (\i,-4) -- cycle;
            \node at (\i-0.35,-3.6)  {\scalebox{0.4}{$\searrow$}};
            }
            \foreach \i in {11}{
             \draw[fill = red!50] (\i-0.25,-3) -- (\i-0.8,-4) -- (\i,-4) -- cycle;
            \node at (\i-0.35,-3.6)  {\scalebox{0.4}{$\searrow$}};
            }
            \foreach \j in {-1,2 ,7,9}{
             \draw[fill = green!50!black, opacity = 0.5] (\j+0.25,-3) -- (\j+0.8,-4) -- (\j,-4) -- cycle;
            \node at (\j+0.35,-3.6)  {\scalebox{0.4}{$\swarrow$}};
            }
            \foreach \j in {-3}{
             \draw[fill = red!50] (\j+0.25,-3) -- (\j+0.8,-4) -- (\j,-4) -- cycle;
            \node at (\j+0.35,-3.6)  {\scalebox{0.4}{$\swarrow$}};
            }
            \draw[fill = green!50!black, opacity = 0.5] (3-0.375,-1.5) -- (3+0.375,-1.5) -- (3,0) -- cycle;
            \node at (3,-1)  {\scalebox{0.6}{$\uparrow$}};
            \draw[fill = green!50!black, opacity = 0.5] (5-0.375,-1.5) -- (5+0.375,-1.5) -- (5,0) -- cycle;
            \node at (5,-1)  {\scalebox{0.6}{$\uparrow$}};

            \draw[fill = green!50!black, opacity = 0.5] (2.5,0) -- (5.5,0) -- (4,2) -- cycle;
            \node at (4,1)  {\scalebox{0.6}{$\Pi_{\ell}$}}; 

            \node at (4,2.5) {\scalebox{0.8}{$T$}};
            \node at (-2,-2) {\scalebox{0.6}{$B_1$}};
            \node at (0,-2) {\scalebox{0.6}{$B_2$}};
            \node at (3,-2) {\scalebox{0.6}{$B_{t_{\ell}}$}};
            \node at (5,-2) {\scalebox{0.4}{$B_{t_{\ell}+1}$}};
            \node at (8,-2) {\scalebox{0.4}{$B_{N_{\ell}-1}$}};
            \node at (10,-2) {\scalebox{0.6}{$B_{N_{\ell}}$}};
            \end{scope}

            \begin{scope}[shift = {(16,-12)}]
                \draw[fill = blue!10] (-2,0) -- (10,0) -- (4,6) -- cycle;
            \draw[fill = blue!10] (-3,-4) -- (-1,-4) -- (-2,0) -- cycle;
            \draw[fill = blue!10] (-1,-4) -- (1,-4) -- (0,0) -- cycle;
            \draw[fill = blue!10] (2,-4) -- (4,-4) -- (3,0) -- cycle;
            \draw[fill = blue!10] (4,-4) -- (6,-4) -- (5,0) -- cycle;
            \draw[fill = blue!10] (7,-4) -- (9,-4) -- (8,0) -- cycle;
            \draw[fill = blue!10] (9,-4) -- (11,-4) -- (10,0) -- cycle;
            \node at (1.5,-4){\scalebox{0.6}{$\cdots$}};
            \node at (6.5,-4) {\scalebox{0.6}{$\cdots$}};
            \node at (4,-5) {port: $\Pi_{\ell+1}$}; 
            \draw[fill = green!50!black, opacity = 0.5] (-2,0) -- (-0.8,0) -- (-0.8,1.2) -- cycle;
            \node at (-1.2,0.4)  {\scalebox{0.6}{$\swarrow$}};
            \draw[fill = green!50!black, opacity = 0.5] (10,0) -- (8.8,0) -- (8.8,1.2) -- cycle;
            \node at (9.2,0.4)  {\scalebox{0.6}{$\searrow$}};

            \draw[fill = green!50!black, opacity = 0.5] (-2-0.375,-1.5) -- (-2+0.375,-1.5) -- (-2,0) -- cycle;
            \node at (-2,-1)  {\scalebox{0.6}{$\uparrow$}};
            \draw[fill = green!50!black, opacity = 0.5] (10-0.375,-1.5) -- (10+0.375,-1.5) -- (10,0) -- cycle;
            \node at (10,-1)  {\scalebox{0.6}{$\uparrow$}};

            \foreach \i in {-1,1,6,9}{
             \draw[fill = green!50!black, opacity = 0.5] (\i-0.25,-3) -- (\i-0.8,-4) -- (\i,-4) -- cycle;
            \node at (\i-0.35,-3.6)  {\scalebox{0.4}{$\searrow$}};
            }
            \foreach \i in {4}{
             \draw[fill = red!50] (\i-0.25,-3) -- (\i-0.8,-4) -- (\i,-4) -- cycle;
            \node at (\i-0.35,-3.6)  {\scalebox{0.4}{$\searrow$}};
            }
            \foreach \j in {-1,2 ,7,9}{
             \draw[fill = green!50!black, opacity = 0.5] (\j+0.25,-3) -- (\j+0.8,-4) -- (\j,-4) -- cycle;
            \node at (\j+0.35,-3.6)  {\scalebox{0.4}{$\swarrow$}};
            }
            \foreach \j in {4}{
             \draw[fill = red!50] (\j+0.25,-3) -- (\j+0.8,-4) -- (\j,-4) -- cycle;
            \node at (\j+0.35,-3.6)  {\scalebox{0.4}{$\swarrow$}};
            }

            \node at (4,2.5) {\scalebox{0.8}{$T$}};
            \node at (-2,-2) {\scalebox{0.6}{$B_1$}};
            \node at (0,-2) {\scalebox{0.6}{$B_2$}};
            \node at (3,-2) {\scalebox{0.6}{$B_{t_{\ell}}$}};
            \node at (5,-2) {\scalebox{0.4}{$B_{t_{\ell}+1}$}};
            \node at (8,-2) {\scalebox{0.4}{$B_{N_{\ell}-1}$}};
            \node at (10,-2) {\scalebox{0.6}{$B_{N_{\ell}}$}};
            \end{scope}
        \end{tikzpicture}
        \caption{Construction of ports for $d=1$.}
        \label{fig:lacosaunodimensional}
    \end{figure}
 
    We now list hereditary Hamiltonian paths realizing each port. In each case, the path inside a piece is given by the inductive hypothesis, and consecutive pieces are joined by one of the two junctions mentioned above. The paths described below are also schematized in~\Cref{fig:lacosaunodimensional}.

    \begin{description}
        \item[$(\uparrow,\swarrow)$:] $(\uparrow,\searrow)_T$, then $(\uparrow,\swarrow)_{B_{N_{\ell}}}$, then $(\searrow,\swarrow)_{B_m}$ for $m=N_{\ell}-1,\dots,1$. It ends at $\swarrow_{B_1}=\swarrow$.
        \item[$(\uparrow,\searrow)$:] $(\uparrow,\swarrow)_T$, then $(\uparrow,\searrow)_{B_1}$, then $(\swarrow,\searrow)_{B_m}$ for $m=2,\dots,N_\ell$. This starts at $\uparrow_T = \uparrow$ and ends at $\searrow_{B_{N_\ell}}=(2^{\ell+1},t_{\ell+1}) = \searrow$.
        \item[$(\swarrow,\searrow)$:] $(\swarrow,\searrow)_{B_m}$ for $m=1,\dots,t_\ell-1$, then $(\swarrow,\uparrow)_{B_{t_{\ell}}}$, then the seam $\Pi_{\ell}$ of $T$, then $(\uparrow,\searrow)_{B_{t_{\ell}+1}}$, then $(\swarrow,\searrow)_{B_m}$ for $m = t_{\ell}+2,\dots,N_{\ell}$.
        \item[$\Pi_{\ell+1}$:] $(\searrow,\swarrow)_{B_m}$ for $m=t_{\ell},t_{\ell}-1,\dots,2$, then $(\searrow,\uparrow)_{B_1}$, then $(\swarrow,\searrow)_T$, then $(\uparrow,\swarrow)_{B_{N_{\ell}}}$, then $(\searrow,\swarrow)_{B_m}$ for $m=N_{\ell}-1,\dots,t_{\ell}+1$. 
    \end{description}
    With this, we have described how to create the ports within $F_{\ell+1}$ from the ports in $F_{\ell}$. It follows that the collection of fans $(F_{\ell})_{\ell \geq 0}$ is hierarchically Hamiltonian.
\end{proof}

\subsection{The case \texorpdfstring{$d\geq 2$}{d>=2}}

As we did in the previous case, we distinguish the vertices of $F_{\ell}$ which will belong to its ports. In this case we will not need to use the seam, and in fact the description of the ports is easier (although the paths are more complex). The \define{root} of $F_{\ell}$ is the vertex $\uparrow=(1,\mathbf{1})$ and, for every $\varepsilon\in\{0,1\}^d$, the \textbf{corner} $C_\varepsilon$ is the vertex $(2^\ell,c_\varepsilon)$ of its bottom row, where $(c_\varepsilon)_k = 1$ if $\varepsilon_k=0$ and $(c_\varepsilon)_k = 2^{2^\ell-1}$ if $\varepsilon_k=1$. For $\ell\geq 1$ these $2^d+1$ points are pairwise distinct. We write $\bar\varepsilon = \mathbf{1}-\varepsilon$ and call $C_{\bar\varepsilon}$ the corner \define{opposite} to $C_\varepsilon$.
 
Once again, consider $N_\ell = 2^{2^\ell}$ and recall the decomposition,
\begin{equation*}
    F_{\ell+1} = F_\ell\ \sqcup\ \bigsqcup_{m\in\{1,\dots,N_\ell\}^d}\varphi_{2^\ell,m-\mathbf{1}}(F_\ell).
\end{equation*}
We call the elements $m \in \{1,\dots,N_{\ell}\}^d$ \define{cells}, $T=F_\ell$ the \define{top piece} and $B_m=\varphi_{2^\ell,m-\mathbf{1}}(F_\ell)$ the \define{bottom piece} at the cell $m$. Two cells are \textbf{adjacent} if they differ by $\pm e_k$ for some $k$. We write $m'-m>0$ if $m'-m=e_k$ and $m'-m<0$ if $m'-m=-e_k$ for some $k$ (this will make the notations simpler in the proofs below). The root of $B_m$ is $\uparrow_{B_m} = (2^\ell+1,m)$, and we denote the images of the distinguished points of $F_\ell$ in a piece $P$ by $\uparrow_P$ and $C_{\varepsilon,P}$. 

Before going into the proof, we will need the fact that the $d$-dimensional grid $\{1,\dots,n\}^d$ admits a Hamiltonian cycle whenever $n\geq 2$ is even and $d\geq 2$. This is what makes the seams defined in the previous section unnecessary here. This elementary result can be easily shown by induction by starting with a Hamiltonian path and inductively using the product to produce a snake-like cycle. An even stronger property satisfied by these grids can be found in~\cite{Ruskey_Sawada_hamiltoniangrid}.

\begin{prop}\label{prop:Hamiltonian_grid}
    For every $\ell \geq 0$, the grid graph $\{1,\dots,N_{\ell}\}^{d}$ admits a Hamiltonian cycle.
\end{prop}

Next, we show that we can reduce the whole analysis to only two ports.

\begin{lemma}
\label{lem:port_symmetric}
    For every $\ell\geq 1$, if both $(\uparrow, C_{\mathbf{0}})$ and $(C_{\mathbf{0}}, C_{\mathbf{1}})$ are ports of $F_{\ell}$, then so are $(\uparrow, C_{\varepsilon})$ and $(C_{\varepsilon}, C_{\bar{\varepsilon}})$ for every $\varepsilon\in\{0,1\}^d$.
\end{lemma}

\begin{proof}
    Let $1\leq k\leq d$ and consider the reflection $\sigma_k\colon F_\ell\to F_{\ell}$ which fixes coordinates $j_{k'}$ with $k'\neq k$, and maps $j_k$ to $2^{i-1}+1-j_k$. It is a bijection of $F_\ell$ which preserves adjacency. Moreover, a direct computation shows that $\sigma_k\circ\varphi_{a,b} = \varphi_{a,b'}\circ\sigma_k$, where $b'$ is obtained from $b$ by replacing $b_k$ by $2^a-1-b_k$. Hence $\sigma_k$ maps the top piece of $F_{\ell+1}$ onto itself and $B_m$ onto $B_{m'}$, where $m'_k=N_\ell+1-m_k$. Finally, $\sigma_k$ fixes $\uparrow$ and maps $C_\varepsilon$ to $C_{\varepsilon+e_k}$, with the sum taken modulo $2$. Consequently, by induction on $\ell$, if $(u,v)$ is a port, so is $(\sigma_k(u),\sigma_k(v))$.

    Since ports are symmetric and invariant under the reflections $\sigma_k$, and the group generated by the reflections acts transitively on the corners and maps opposite corners to opposite corners, if both $(\uparrow,C_{\mathbf{0}})$ and $(C_{\mathbf{0}},C_{\mathbf{1}})$ are ports, then so are $(\uparrow,C_{\varepsilon})$ and $(C_{\varepsilon}, C_{\bar{\varepsilon}})$ for every $\varepsilon\in\{0,1\}^d$.
\end{proof}

\begin{prop}\label{prop:fan_ports_d}
    For every $\ell\geq 1$ and every $\varepsilon\in\{0,1\}^d$, the pairs $(\uparrow,C_\varepsilon)$ and $(C_\varepsilon,C_{\bar\varepsilon})$ are ports of $F_\ell$.
\end{prop}

\begin{figure}[ht!]
    \centering
    \fadeimage{width=0.3\textwidth}{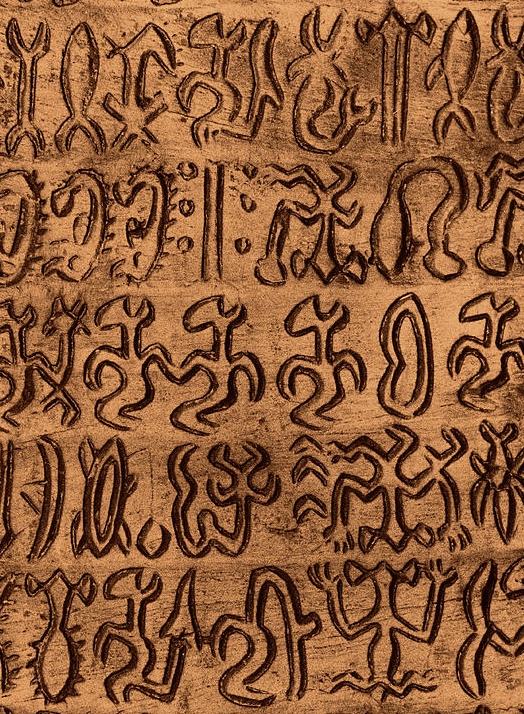}
    \caption{A portion of rongorongo script from Rapa Nui. Notice that the symbols on subsequent rows are inverted. Public domain image taken from~\cite{englert1970island}.}
    \label{fig:Rongorongo}
\end{figure}

\begin{proof}
    By~\Cref{lem:port_symmetric}, it suffices to show that for every $\ell\geq 1$ the pairs $(\uparrow, C_{\mathbf{0}})$ and $(C_{\mathbf{0}}, C_{\mathbf{1}})$ are ports.

    Let us consider the base case $\ell=1$. Here $F_1$ consists of $\uparrow=(1,\mathbf{1})$ and the $2^d$ vertices $(2,c)$ with $c\in\{1,2\}^d$. The vertices $(2,c)$ with the horizontal edges form a $d$-dimensional hypercube, which admits a Hamiltonian cycle $v_0,v_1,\dots,v_{2^d-1}$ by~\Cref{prop:Hamiltonian_grid} (for instance, the one given by a Gray code~\cite{frank1953pulse}) with $v_0 = C_{\mathbf{0}}$. Then $\uparrow,v_1,v_2,\dots,v_{2^d-1},v_0$ realizes $(\uparrow,C_{\mathbf{0}})$. If $C_{\mathbf{1}}=v_s$, then $v_0,v_1,\dots,v_{s-1},\uparrow,v_{2^d-1},v_{2^d-2},\dots,v_s$ realizes $(C_{\mathbf{0}},C_{\mathbf{1}})$.

    Now suppose the statement holds for some $\ell\geq 1$, and decompose $F_{\ell+1}$ as above. We will use the following three kinds of junctions between pieces.

    \begin{enumerate}
    \item If $m$ and $m+e_k$ are cells, then $\uparrow_{B_m}$ and $\uparrow_{B_{m+e_k}}$ are horizontal neighbors.
    \item If $m$ and $m+e_k$ are cells and $\varepsilon\in \{0,1\}^d$ is such that $\varepsilon_k=1$, then $C_{\varepsilon,B_m}$ and $C_{\varepsilon-e_k,B_{m+e_k}}$ are horizontal neighbors.
    \item For $\varepsilon\in\{0,1\}^d$, let $Q_\varepsilon\subset \{1,\dots,N_{\ell}\}^d$ be the set of cells $m$ with $m_k\in\{1,2\}$ if $\varepsilon_k=0$ and $m_k\in\{N_{\ell}-1,N_{\ell}\}$ if $\varepsilon_k=1$. Then $C_{\varepsilon,T}$ is the parent of $\uparrow_{B_m}$ for every $m\in Q_\varepsilon$.
\end{enumerate}

Note that the points of $F_{\ell+1}$ we are interested in realizing through ports are $\uparrow=\uparrow_T$, $C_{\mathbf{0}} = C_{\mathbf{0},B_{\mathbf{1}}}$ and $C_{\mathbf{1}} = C_{\mathbf{1},B_{N_{\ell}\mathbf{1}}}$.

For the remainder of the proof, we fix a Hamiltonian cycle $v_0,\dots, v_{N_{\ell}^d-1}$ on the grid graph $\{1,\dots,N_{\ell}\}^d$, which we will interpret as a cycle between the bottom pieces. We interpret the indices as elements of $\ZZ/(N_{\ell}^d)\ZZ$. We fix the cycle in such a way that $v_0 = \mathbf{1}$, so that it corresponds to the piece $B_{\mathbf{1}}$, and we denote by $s \in \{1,\dots,N_{\ell}^d-1\}$ the index for which $v_s = N_{\ell}\mathbf{1}$, so that it corresponds to $B_{N_{\ell}\mathbf{1}}$. As the grid is bipartite, with the parity of the sum of the coordinates as a coloring, and the sums of the coordinates of $\mathbf{1}$ and $N_\ell\mathbf{1}$ differ by $d(N_\ell-1)$, we have $s\equiv d \pmod 2$. Note also that $v_1$ and $v_{N_\ell^d-1}$ are adjacent to $\mathbf{1}$, and thus belong to $Q_\mathbf{0}$, while $v_{s-1}$ and $v_{s+1}$ are adjacent to $N_\ell\mathbf{1}$, and thus belong to $Q_\mathbf{1}$.

To create our hereditary Hamiltonian paths, we will group our ports into three types of moves. We denote the \define{down move} by $D$, which represents any port of the form $(\uparrow_{B_m}, C_{\varepsilon,B_m})$ for some cell $m$ and $\varepsilon\in \{0,1\}^d$. The \define{up move} is denoted by $U$ and represents the path with ports $(C_{\varepsilon,B_m}, \uparrow_{B_m})$. Finally, a \define{diagonal move} $\Delta$ is a path with ports $(C_{\varepsilon,B_m}, C_{\bar{\varepsilon},B_m})$.

A \textbf{rongorongo path} is a sequence of consecutive bottom pieces along the cycle which are traversed alternately by down and up moves. We chose this name because the alternation resembles the rongorongo script from Rapa Nui, in which every other line is written upside down (see~\Cref{fig:Rongorongo}). The idea is to use a rongorongo path in different portions of the Hamiltonian cycle to create a hereditary Hamiltonian path for each port by using the junctions. The main idea of the construction is depicted in~\Cref{fig:ports_higher_dimension} and the details are given below.

\begin{figure}[ht!]
    \centering
    \begin{tikzpicture}[scale = 0.2]
    
        \begin{scope}[shift = {(0,6)},cm={1,0,0.5,0.5,(0,0)}]
            \draw[fill = blue!10,opacity = 0.8] (0,0) rectangle (16,16);
            \begin{scope}[shift = {(2,2)}]
                \draw[blue, thick, opacity = 0.5] (0,0) -- (0,12) -- (4,12) -- (4,4) -- (8,4) -- (8,12) -- (12,12) -- (12,0) -- cycle;
                \node[rotate=0] at (0,0) {$v_0$};
                \node[rotate=0] at (0,4) {$v_1$};
                \node[rotate=0] at (0,8) {$v_2$};
                \node[rotate=0] at (0,12) {$v_3$};
                \node[rotate=0] at (4,12) {$v_4$};
                \node[rotate=0] at (4,8) {$v_5$};
                \node[rotate=0] at (4,4) {$v_6$};
                \node[rotate=0] at (8,4) {$v_7$};
                \node[rotate=0] at (8,8) {$v_8$};
                \node[rotate=0] at (8,12) {$v_9$};
                \node[rotate=0] at (12,12) {$v_{10}$};
                \node[rotate=0] at (12,8) {$v_{11}$};
                \node[rotate=0] at (12,4) {$v_{12}$};
                \node[rotate=0] at (12,0) {$v_{13}$};
                \node[rotate=0] at (8,0) {$v_{14}$};
                \node[rotate=0] at (4,0) {$v_{15}$};
            \end{scope}
            \draw[step=4.0cm] (0,0) grid (16,16);
        \end{scope}
        \begin{scope}[shift = {(0,14)},cm={1,0,0.5,0.5,(0,0)}]
            \draw[fill = blue!10,opacity = 0.8] (4,4) rectangle (12,12);
            \node[rotate=0] at (8,8) {$T$};
        \end{scope}

        \begin{scope}[shift = {(-20,0)}]
            \begin{scope}[shift = {(0,6)},cm={1,0,0.5,0.5,(0,0)}]
            \draw[fill = blue!10,opacity = 0.8] (0,0) rectangle (16,16);
            \begin{scope}[shift = {(2,2)}]
                \draw[color=green!50!black, very thick, dashed, ->] (0,4) -- (0,12) -- (4,12) -- (4,4) -- (8,4) -- (8,12) -- (12,12) -- (12,0) -- (0,0);
                \draw[ fill =red ] (-2,-2) rectangle (-1,-1);
                \begin{scope}[opacity = 0.6]
                \node at (0,0) {$\Delta$};
                \node[rotate=0] at (0,4) {$D$};
                \node[rotate=0] at (0,8) {$U$};
                \node[rotate=0] at (0,12) {$D$};
                \node[rotate=0] at (4,12) {$U$};
                \node[rotate=0] at (4,8) {$D$};
                \node[rotate=0] at (4,4) {$U$};
                \node[rotate=0] at (8,4) {$D$};
                \node[rotate=0] at (8,8) {$U$};
                \node[rotate=0] at (8,12) {$D$};
                \node[rotate=0] at (12,12) {$U$};
                \node[rotate=0] at (12,8) {$D$};
                \node[rotate=0] at (12,4) {$U$};
                \node[rotate=0] at (12,0) {$D$};
                \node[rotate=0] at (8,0) {$U$};
                \node[rotate=0] at (4,0) {$D$};
                \end{scope}
            \end{scope}
            \draw[step=4.0cm] (0,0) grid (16,16);
        \end{scope}
        \begin{scope}[shift = {(0,14)},cm={1,0,0.5,0.5,(0,0)}]
            \draw[fill = blue!10,opacity = 0.8] (4,4) rectangle (12,12);
            \draw[fill = green!50!black,opacity = 1] (4,4) rectangle (6,6);
            \draw[fill = red] (7,7) rectangle (9,9);
            
        \end{scope}
        \draw[dashed, very thick, color=green!50!black, ->] (6,16) to (4,9);

        \end{scope}

        \begin{scope}[shift = {(20,0)}]
            \begin{scope}[shift = {(0,6)},cm={1,0,0.5,0.5,(0,0)}]
            \draw[fill = blue!10,opacity = 0.8] (0,0) rectangle (16,16);
            \begin{scope}[shift = {(2,2)}]
                \draw[color=green!50!black, very thick, dashed, ->] (0,0) -- (0,4) -- (0,12) -- (4,12) -- (4,4) -- (8,4) -- (8,12);
                \draw[color=green!50!black, very thick, dashed, ->] (8,12) -- (2,24);
                \draw[color=green!50!black, very thick, dashed, ->] (-5,18) -- (4,0);
                \draw[color=green!50!black, very thick, dashed, ->] (4,0) -- (12,0) -- (12,12);
                \draw[ fill =red ] (-2,-2) rectangle (-1,-1);
                \draw[ fill =red ] (13,13) rectangle (14,14);
                \begin{scope}[opacity = 0.6]
                \node at (0,0) {$\Delta$};
                \node[rotate=0] at (0,4) {$U$};
                \node[rotate=0] at (0,8) {$D$};
                \node[rotate=0] at (0,12) {$U$};
                \node[rotate=0] at (4,12) {$D$};
                \node[rotate=0] at (4,8) {$U$};
                \node[rotate=0] at (4,4) {$D$};
                \node[rotate=0] at (8,4) {$U$};
                \node[rotate=0] at (8,8) {$D$};
                \node[rotate=0] at (8,12) {$U$};
                \node[rotate=0] at (12,12) {$\Delta$};
                \node[rotate=0] at (12,8) {$D$};
                \node[rotate=0] at (12,4) {$U$};
                \node[rotate=0] at (12,0) {$D$};
                \node[rotate=0] at (8,0) {$U$};
                \node[rotate=0] at (4,0) {$D$};
                \end{scope}
            \end{scope}
            \draw[step=4.0cm] (0,0) grid (16,16);
        \end{scope}
        \begin{scope}[shift = {(0,14)},cm={1,0,0.5,0.5,(0,0)}]
            \draw[fill = blue!10,opacity = 0.8] (4,4) rectangle (12,12);
            \draw[fill = green!50!black,opacity = 1] (4,4) rectangle (6,6);
            \draw[fill = green!50!black,opacity = 1] (10,10) rectangle (12,12);
            \draw[color=green!50!black, very thick, dashed, ->] (10,10) -- (6,6);
            
        \end{scope}

        \end{scope}
        \node at (-10,2) {port: $(\uparrow,C_{\mathbf{0}})$};
        \node at (30,2) {port: $(C_{\mathbf{0}},C_{\mathbf{1}})$};
    \end{tikzpicture}
    \caption{Diagram depicting the hereditary Hamiltonian paths in $F_{\ell+1}$ for $d= 2$. In the middle we depict the Hamiltonian cycle on the bottom pieces, here $s=10$.}
    \label{fig:ports_higher_dimension}
\end{figure}

Case $(\uparrow, C_{\mathbf{0}})$: We first describe the order in which the pieces will be traversed. Our path will traverse in order:

    \[ T \rightarrow B_{v_1}\rightarrow B_{v_2}\rightarrow \dots\rightarrow B_{v_{N_{\ell}^d-1}} \rightarrow B_{v_0}. \]

    We begin assigning the port $(\uparrow_T,C_{\mathbf{0},T})$ to $T$, and then we assign a rongorongo path to $B_{v_1} \dots B_{v_{N_{\ell}^d-1}}$ which begins with $D$. Note that this path ends with $D$, as it consists of an odd number of pieces. Finally, we assign the port $\Delta = (C_{\mathbf{1},B_{v_0}},C_{\mathbf{0},B_{v_0}})$ to $B_{v_0}$.

    It remains to specify the exact ports we use along the rongorongo path. To each cell $v_i$ with $i$ odd, we assign the port $D_i = (\uparrow_{B_{v_i}},C_{\varepsilon_i,B_{v_i}})$ where $\varepsilon_i$ is given by \[ \varepsilon_i = \begin{cases}
        \mathbf{1}+(v_{i+1}-v_i) & \mbox{ if } v_{i+1}-v_i < 0\\
        (v_{i+1}-v_i) & \mbox{ if } v_{i+1}-v_i > 0.
    \end{cases}  \]
    With this choice, we can assign to each even index $i$ the port $U_i = (C_{\delta_i,B_{v_i}},\uparrow_{B_{v_i}})$ where $\delta_i$ is given by \[ \delta_i = \begin{cases}
        \mathbf{1} & \mbox{ if } v_{i}-v_{i-1} < 0\\
        \mathbf{0} & \mbox{ if } v_{i}-v_{i-1} > 0.
    \end{cases}  \]

    Hence the complete path of ports is given by \[ (\uparrow_T,C_{\mathbf{0},T}) \rightarrow D_1 \rightarrow U_2 \rightarrow \dots \rightarrow D_{N_{\ell}^d-1} \to (C_{\mathbf{1},B_{v_0}},C_{\mathbf{0},B_{v_0}}). \]
    Note that each connection between ports corresponds to one of the junctions described above, hence $(\uparrow, C_{\mathbf{0}})$ is a port of $F_{\ell+1}$.

    Case $(C_{\mathbf{0}}, C_{\mathbf{1}})$: As in the previous case, we begin by specifying the order in which the pieces are traversed. In this case we take \[ B_{v_0} \to B_{v_1} \to \dots \to B_{v_{s-1}} \to T \to B_{v_{N_{\ell}^d-1}} \to B_{v_{N_{\ell}^d-2}} \to \dots \to B_{v_{s+1}} \to B_{v_{s}}.  \]
    
    Suppose first that $d$ is even, so that $s$ is even. We start by assigning the $\Delta$ port $\Delta_0 = (C_{\mathbf{0},B_{v_0}},C_{\mathbf{1},B_{v_0}})$ to $B_{v_0}$. Next, we assign a rongorongo path to $B_{v_1},\dots, B_{v_{s-1}}$ starting with $U$. As $s$ is even, this path ends with a $U$. Next, we assign the $\Delta$ port $\Delta_T = (C_{\mathbf{1},T},C_{\mathbf{0},T})$ to $T$, and continue by assigning a rongorongo path starting with $D$ from $B_{v_{N_{\ell}^d-1}}$ up to $B_{v_{s+1}}$, which again, by parity, ends with a $D$. Finally, we assign the $\Delta$ port $\Delta_s =(C_{\mathbf{0},B_{v_s}},C_{\mathbf{1},B_{v_s}})$ to $B_{v_s}$. This can be schematized as follows \[ \Delta_0 \rightarrow U_1 \rightarrow D_2 \to \dots \to U_{s-1} \to \Delta_T \to D_{N_{\ell}^d-1} \to U_{N_{\ell}^d-2} \to \dots \to D_{s+1} \to \Delta_s.   \]

    If $d$ is odd, then $s$ is odd, and we modify the construction at its two ends: we assign the up move $(C_{\mathbf{0},B_{v_0}},\uparrow_{B_{v_0}})$ to $B_{v_0}$ and the down move $(\uparrow_{B_{v_s}},C_{\mathbf{1},B_{v_s}})$ to $B_{v_s}$. Both rongorongo paths then start with $D$ and, having an even number of pieces, end with $U$. They are joined to $B_{v_0}$ and $B_{v_s}$ through their roots by the first junction, and to $T$ as before. The ports along the rongorongo paths are chosen exactly as in the even case, according to whether each piece is traversed by an up or a down move.

    It remains to describe the precise ports $U_i$ and $D_i$ assigned to the pieces $B_{v_i}$ along the rongorongo paths. To that end, for the first rongorongo path ($i \in \{1,\dots,s-1\}$) we take $U_i = (C_{\delta_i,B_{v_i}},\uparrow_{B_{v_i}})$ for the pieces traversed by an up move, and $D_i = (\uparrow_{B_{v_i}}, C_{\varepsilon_i,B_{v_i}})$ for those traversed by a down move, where \[ \delta_i = \begin{cases}
        \mathbf{1}-(v_{i}-v_{i-1}) & \mbox{ if } v_{i}-v_{i-1} > 0\\
        -(v_{i}-v_{i-1}) & \mbox{ if } v_{i}-v_{i-1} < 0,
    \end{cases}
    \quad\quad \varepsilon_i = \begin{cases}
        \mathbf{1} & \mbox{ if } v_{i+1}-v_{i} > 0\\
        \mathbf{0} & \mbox{ if } v_{i+1}-v_{i} < 0.
    \end{cases} \]

    For the second rongorongo path ($i \in \{N_{\ell}^d-1,\dots,s+1\}$), we take $D_i = (\uparrow_{B_{v_i}}, C_{\varepsilon'_i,B_{v_i}})$ for the pieces traversed by a down move, and $U_i = (C_{\delta'_i,B_{v_i}},\uparrow_{B_{v_i}})$ for those traversed by an up move, where in this case the positions are given by \[ \varepsilon'_i = \begin{cases}
        v_{i-1}-v_i & \mbox{ if } v_{i-1}-v_i > 0\\
        \mathbf{1}+(v_{i-1}-v_i) & \mbox{ if } v_{i-1}-v_i < 0,
    \end{cases} \quad\quad \delta'_i = \begin{cases}
        \mathbf{0} & \mbox{ if } v_{i}-v_{i+1} > 0\\
        \mathbf{1} & \mbox{ if } v_{i}-v_{i+1} < 0.
    \end{cases} \]

    With these choices, each consecutive sequence of ports corresponds to one of the junctions described above, hence we conclude that $(C_{\mathbf{0}}, C_{\mathbf{1}})$ is a port of $F_{\ell+1}$.\end{proof}

\subsection{Grand finale}

With all the tools in hand, we are ready to show the main lemma of this section.

\begin{lemma}\label{lem:BDSM}
    For every $d\geq 1$, the graph $\mathcal{B}_{d+1}$ is $(1,f)$-space-filling for some computable function $f$.
\end{lemma}

\begin{proof}
    By~\Cref{lem:balls_fans}, for every $r\geq 1$ there exists a computable map $r\mapsto \ell(r)$ with the property that $F_{\ell(r)}$ contains a ball of radius $r$. Take $f(0)=0$ and $f(r)=2|F_{\ell(r)}|$ for $r\geq 1$, and note that this map is computable. Moreover, by~\Cref{prop:fan_ports} and~\Cref{prop:fan_ports_d} the collection of fans $(F_{\ell})_{\ell \geq 0}$ is hierarchically Hamiltonian for all $d\geq 1$. Let $M \in \NN$ and consider any subpath $(w_j)_{j=0}^M$ of a hereditary Hamiltonian path on $F_{\ell}$ for some $\ell$ large enough. By our choice of $f(r)=2|F_{\ell(r)}|$ and the definition of hierarchically Hamiltonian, any window of length $f(r)$ must entirely contain one of the copies $\varphi_{a,b}(F_{\ell(r)})$ in the hierarchy of $F_\ell$, which by~\Cref{lem:balls_fans} contains a ball of radius $r$.
\end{proof}

With this lemma in hand, we are finally able to prove~\Cref{mainthm:hyperbolic}. 

\begin{proof}[Proof of~\Cref{mainthm:hyperbolic}]
    Let $G$ be a finitely generated group which is quasi-isometric to $\HH^{d+1}$ for some $d\geq 0$. Note that $G$ is hyperbolic, and thus it has decidable word problem (see for instance~\cite{MR1744486}). The case $d=0$ is trivial, as in this case $G$ is infinite and virtually cyclic, and thus finitary space-filling by~\Cref{cor:virtab_is_finitary_SF}. Hence we suppose that $d\geq 1$.
    
    By~\Cref{lem:Bd_QI}, the highly computable graph $\mathcal{B}_{d+1}$ is quasi-isometric to $\HH^{d+1}$, hence $\mathcal{B}_{d+1}$ is quasi-isometric to $G$. Moreover, by~\Cref{lem:BDSM} we have that $\mathcal{B}_{d+1}$ is computably space-filling.

    Finally, by~\Cref{thm:QI_effective_graphs} we have that every Cayley graph of $G$ is computably space-filling, and thus by~\Cref{prop:effective_SF} we have that $G$ is computably space-filling.
\end{proof}

We close this section exhibiting a family of finitary space-filling groups which are not virtually abelian.

\begin{example}\label{ex:productitos}
    Let $G$ be any finitely generated group which is quasi-isometric to $\HH^{d+1}$ for $d\geq 1$. By~\Cref{cor:simulation_y_SF} we obtain that every product of the form $\Z^n \times G$ for $n \geq 1$ is finitary space-filling. An explicit new example of a finitary space-filling group is $\pi_1(S_2) \times \ZZ$, where $\pi_1(S_2)$ is the fundamental group of the closed orientable surface of genus $2$.
\end{example}

\section{Applications}
\label{sec:applications}

Our goal in this section is to show that in groups with undecidable domino problem which admit a strongly space-filling directed tileset (a strengthening of space-filling introduced in~\Cref{def:SSF}), some decision problems related to cellular automata and tilings are undecidable. 

The domino problem was first introduced by Wang for $\Z^2$~\cite{wang1962dominoes}, and famously shown to be undecidable by Berger~\cite{berger1966undecidability}. For finitely generated groups, its current formulation is the following.
\begin{defn}
Let $G$ be a finitely generated group and $S$ a finite symmetric set of generators. The \define{domino problem} on $G$ with respect to $S$, denoted $\DP(G,S)$, is the following decision problem: given a finite alphabet $A$ and a finite set $\FF$ of nearest neighbor forbidden patterns, determine whether $\X_{\FF}$ is nonempty.
\end{defn}

It has been conjectured by Ballier and Stein~\cite{ballier2013domino} that a finitely generated group has decidable domino problem if and only if it is virtually free. Although the conjecture remains unproven, it has been established for many classes of groups (see~\cite[Chapter 2]{bitar2024subshifts} for a survey).

The other decision problems we consider are related to cellular automata. We note that a cellular automaton is called \define{reversible} if it admits an inverse which is a cellular automaton, and this notion coincides with being bijective (see~\cite[Chapter 1]{Ceecherini_Coornaert_book_2023}).

\begin{defn}
The \define{injectivity problem} for cellular automata on $G$ is the following decision problem: given a finite alphabet $A$, a finite set $F\Subset G$, and a map $\Phi\colon A^F\to A$, determine whether the cellular automaton $\phi\colon A^G\to A^G$ given by $\phi(x)(g) = \Phi((g^{-1}\cdot x)|_F)$ is injective.

Similarly, the \define{reversibility problem} asks, given the same input, if $\phi$ is reversible.
\end{defn}


\subsection{Snake tiling problems}

Our goal is to follow Kari's proof strategy~\cite{kari2002infinite} of reducing the domino problem to the infinite snake problem in $\Z^2$ to prove its undecidability. To introduce the latter problem, we need the following definitions.

We write sets of nearest neighbor forbidden patterns through triples $(a,b,s)\in A\times A\times S$ which represent a pattern of support $\{1_G,s\}$, with $a$ at the identity and $b$ at $s$. Consider a configuration $x\in A^G$ and a set of nearest neighbor forbidden patterns $\FF\subseteq A\times A\times S$. We say that the edge $\{g,gs\}$ of the Cayley graph is \define{valid in $x$ with respect to $\FF$} if $(x(g),x(gs),s)\notin\FF$ and $(x(gs),x(g),s^{-1})\notin\FF$. This way the subshift $\X_{\FF}$ corresponds to the set of configurations in which every edge is valid.
 
\begin{defn}\label{def:directed_tileset}
A \define{directed tileset} over $(G,S)$ is a triple $(A,\FF,\delta)$ where $A$ is a finite alphabet, $\FF\subseteq A\times A\times S$ is a set of nearest neighbor forbidden patterns and $\delta\colon A\to S$. Given $x\in A^G$ its \define{skeleton} is the map $\omega_x\colon\N\to G$ given by $\omega_x(n)=R^n_x(1_G)$, where $R\colon A^G\to G$ is given by $R(x) = \delta(x(1_G))$, so that $R_x(g) = g\delta(x(g))$. In particular, $\omega_x(0)=1_G$ and $\omega_x(n+1)=\omega_x(n)\,\delta(x(\omega_x(n)))$.
\end{defn}
 
Domino snake problems were first introduced by Myers in~\cite{myers1979decidability} as a variant of the standard domino problem. Among these, the infinite snake problem stands out; it was proven to be undecidable by Adleman, Kari, Kari, Reishus and Sosík~\cite{adleman2002decidability,adleman2009undecidability}. The definitions for general groups we present were introduced in~\cite{aubrun2023domino} and further used in~\cite{bartholdi2024snakes}.
 
\begin{defn}
Let $(A,\FF,\delta)$ be a directed tileset and $x\in A^G$.
\begin{itemize}
\item We say $x$ is an \define{infinite snake} if $\omega_x$ is injective, and every edge $\{g,gs\}$ with $g,gs\in\omega_x(\N)$ is valid in $x$.
\item For $N\in\N$, we say $x$ is a \define{snake of length $N$} if $\omega_x$ is injective on $\{0,\dots,N\}$ and every edge with both endpoints in $\omega_x(\{0,\dots,N\})$ is valid in $x$.
\end{itemize}
The \define{strong directed snake problem}, denoted $\SDSP(G,S)$, asks, given a directed tileset, whether it admits an infinite snake.
\end{defn}

This problem is known to be undecidable for virtually nilpotent groups which are not virtually cyclic, and decidable for virtually free groups~\cite{aubrun2023domino}. However, the groups with decidable strong directed snake problem do not coincide with those with decidable domino problem: Bartholdi and Salo constructed a non-virtually free group that has decidable strong directed snake problem~\cite{bartholdi2024snakes}.

\begin{remark}
    In the literature there are multiple versions of the infinite snake problem. One can ask for the snake to be a weak snake, where edges are required to be valid only along the skeleton, or even undirected, where the snake can advance in any direction provided the skeleton remains injective. All of these variants many-one reduce to the strong directed snake problem~\cite[Lemma 2.4]{bartholdi2024snakes}. In particular, they are decidable whenever $\SDSP(G,S)$ is.
\end{remark}

\subsection{Locally space-filling}

To achieve the goal of getting a reduction from the domino problem to the snake problem, we need to introduce a stronger version of space-filling. Indeed, if $X$ is a space-filling subshift, the skeleton of every configuration in $X$ covers translates of every finite subset of the group. However, an infinite snake only has to be correct in the support defined by its skeleton, so it need not be the restriction of a configuration of $X$, and its skeleton need not cover anything. Therefore, we must not only ask full configurations to be space-filling, but also any infinite snake over the directed tileset.

\begin{defn}\label{def:SSF}
A directed tileset $(A,\FF,\delta)$ over $(G,S)$ is \define{strongly space-filling} if
\begin{enumerate}
\item (Existence) there exists $x\in \X_{\FF}$ whose skeleton $\omega_x$ is injective, and
\item (Space-filling) for every infinite snake $x\in A^G$ and every $F\Subset G$ there is $h\in G$ with $hF\subseteq\omega_x(\N)$.
\end{enumerate}
\end{defn}

For this notion of space-filling, a standard compactness argument shows the following.
\begin{lemma}\label{lem:compact_strong_snake}
    A directed tileset $(A,\FF,\delta)$ satisfies the space-filling condition of~\Cref{def:SSF} if and only if for every $n\in\N$ there exists $N\in\N$ such that for every snake $x$ of length $N$ there exists $h\in G$ with $hB(n)\subseteq\omega_x(\llbracket 0,N\rrbracket)$.
\end{lemma}

The natural question that follows is how this stronger version of space-filling relates to our original definition. To answer it we need the following notion.
\begin{defn}
\label{def:strong_sf}
Given a directed tileset $(A,\FF,\delta)$ over $(G,S)$ we say a configuration $x\in A^G$ has a \define{valid directed cycle} if $\omega_x$ is periodic and every edge $\{g,gs\}$ with $g,gs\in\omega_x(\N)$ is valid in $x$. We say that the directed tileset is \define{cycle-free} if no configuration $x\in A^G$ has a valid directed cycle.
\end{defn}

\begin{example}
The notion of being strongly space-filling is actually already present in Kari's original construction~\cite{Kari_reversibility_defofpeanocurve,adleman2002decidability,adleman2009undecidability}. Indeed, it is shown there that any infinite snake of the tileset covers arbitrarily large squares. This is precisely the notion of strongly space-filling. Furthermore, the constructed tileset is also cycle-free.
\end{example}

\begin{prop}\label{prop:SSFtoSF}
Let $(A,\FF,\delta)$ be a strongly space-filling directed tileset over $(G,S)$. Let $X_{\FF}^{\mathrm{inj}}$ be the set of $x\in \X_{\FF}$ such that the skeleton $n\mapsto R^n_x(g)$ is injective for every $g\in G$.
\begin{itemize}
\item $X_{\FF}^{\mathrm{inj}}$ is a nonempty space-filling subshift with respect to $R(x)=\delta(x(1_G))$.
\item If moreover $(A,\FF,\delta)$ is cycle-free, then $X_{\FF}^{\mathrm{inj}}=\X_{\FF}$, so $\X_{\FF}$ is a space-filling SFT.
\end{itemize}
\end{prop}
\begin{proof}
    Clearly $X_{\FF}^{\mathrm{inj}}$ is closed and $G$-invariant. To prove it is nonempty, let $x\in \X_{\FF}$ be the infinite snake given by the strongly space-filling property. For every $k\in\N$ there exists $h_k\in G$ such that $h_kB(k)\subseteq \omega_x(\N)$. Then, take an accumulation point $y$ of the sequence $(h_k^{-1}\cdot x)_{k\in\N}$. As $\X_{\FF}$ is closed and $G$-invariant, $y\in\X_{\FF}$. Next, take $g\in G$ and $N\in\N$. For infinitely many $k\geq|g|$, $y$ and $h_k^{-1}\cdot x$ coincide on the ball of radius $N + |g|$. For such $k$ we have $h_kg\in\omega_x(\N)$, so the sequence $n\mapsto R^n_{h_k^{-1}x}(g) = h_k^{-1}R^n_x(h_kg)$ is a translate of a tail of $\omega_x$, and is thus injective. Its first $N+1$ terms coincide with those of $n\mapsto R^n_y(g)$. As $N$ is arbitrary, the map $n\mapsto R^n_y(g)$ is injective.

    Let us see it is space-filling with $R(x)=\delta(x(1_G))$. Consider a configuration $x\in X_{\FF}^{\mathrm{inj}}$ and a finite subset $F\Subset G$. By~\Cref{lem:compact_strong_snake} with $n\in\N$ such that $F\subseteq B(n)$, there exists $N\in\N$ and $h\in G$ such that \[hF\subseteq hB(n)\subseteq \omega_x(\llbracket 0,N\rrbracket). \]
    Therefore, $X_{\FF}^{\mathrm{inj}}$ is a space-filling subshift. For the second point, suppose there exist $x\in \X_{\FF}$ and $g\in G$ such that $n\mapsto R^n_x(g)$ is not injective. Replacing $x$ by $g^{-1}\cdot x$, we may assume that $g=1_G$. Then, the skeleton must be eventually periodic, that is, for sufficiently large $m$, the configuration $R^m_x(1_G)^{-1}\cdot x$ contains a valid directed cycle.
\end{proof}

The previous result motivates the following definition.

\begin{defn}\label{def:locally_SF}
    A group $G$ is called \define{locally space-filling} if there exists a finite symmetric set of generators $S$ and a cycle-free strongly space-filling directed tileset $(A,\FF,\delta)$ over $(G,S)$.
\end{defn}

\begin{remark}
    Notice that if $G$ has decidable word problem, the subshift $X_{\FF}^{\mathrm{inj}}$ is effective. Therefore, the existence of a strongly space-filling directed tileset over $(G,S)$ for some generating set $S$, implies that $G$ is computably space-filling.
\end{remark}

The natural question that remains is whether every finitary space-filling group is locally space-filling. A more modest question is whether every finitary space-filling group admits a strongly space-filling directed tileset. Although we do not yet know the answer to these questions, we can at least give a restriction on the groups which admit strongly space-filling directed tilesets.

\begin{defn}
An infinite snake $x$ has \define{eventually periodic moves} if its sequence of moves 
\[n\mapsto \delta(x(\omega_x(n)))\in S^{\N}\]
is eventually periodic. This holds in particular if $x$ is \define{eventually periodic}, that is, if its sequence of tiles $\bigl(x(\omega_x(n))\bigr)_{n\ge0}$ is eventually periodic.
\end{defn}
 
\begin{prop}\label{prop:noperiodic}
Let $G$ be a finitely generated group which is not virtually cyclic, and let $(A,\FF,\delta)$ be a strongly space-filling directed tileset. Then it admits no infinite snake with eventually periodic moves.
\end{prop}

\begin{proof}
    Suppose there exists $x$, an infinite snake with eventually periodic moves of period $p$ starting from $n_0\in\N$. Consider the infinite snake $y=\omega_x(n_0)^{-1}\cdot x$, which has periodic moves of period $p$. Consider $s_n = \delta(y(\omega_y(n)))$ and define $g_k = s_0s_1\dots s_{k-1}$ for $0\leq k\leq p$, $g = g_p$, and $D = \max_{0\leq k <p}|g_k|_S$. A short inductive argument shows $\omega_y(mp + r) = g^mg_r$ for all $m\in\N$ and $0\leq r < p$. Thus, $\omega_y(\N)\subseteq \langle g\rangle B(D)$. As the tileset is strongly space-filling, for every $n\geq D$ there exists $h\in G$ such that 
    \[hB(n)\subseteq\omega_y(\N)\subseteq \langle g\rangle B(D).\]
    In particular, there exist $i\in\Z$ and $f\in B(D)$ such that $h = g^if$. Then,
    \[fB(n)\subseteq g^{-i}\langle g\rangle B(D) = \langle g\rangle B(D).\]

    Now, take any $t\in B(n-D)$. We have that $f^{-1}t\in B(n)$, and therefore $t=f(f^{-1}t)\in fB(n)$. This implies $B(n-D)\subseteq\langle g\rangle B(D)$. As $n$ is arbitrary, $G=\langle g\rangle B(D)$, so $\langle g\rangle$ has finite index in $G$ and $G$ is virtually cyclic, a contradiction.
\end{proof}

In~\cite{bartholdi2024snakes}, Bartholdi and Salo construct a finitely generated non-virtually free group $G$ with the following property: for every generating set $S$, and every directed tileset over $(G,S)$, if it admits an infinite snake, it admits an infinite snake with periodic moves (\cite[Proposition 5.4]{bartholdi2024snakes}). Therefore, by the previous proposition, $G$ does not admit strongly space-filling directed tilesets over $(G,S)$ for any $S$. Thus, their group is not locally space-filling.

\subsection{Reducing decision problems}

The classical result by Amoroso and Patt~\cite{AmorosoPatt1972} shows that it is decidable whether a cellular automaton defined on $\ZZ$ is injective (equivalently, as $\ZZ$ is surjunctive, reversible). However, Kari showed in~\cite[Theorem 1]{Kari_reversibility_defofpeanocurve} that the analogous problem in $\ZZ^2$ is undecidable, by reducing the complement of the domino problem to it. The following is a straightforward adaptation of Kari's proof to our setting. 

\begin{theorem}
\label{thm:CA}
    Let $G$ be a finitely generated group and $S$ a finite symmetric set of generators such that there exists a cycle-free strongly space-filling directed tileset over $(G,S)$. Then the complement of $\DP(G,S)$ many-one reduces to both the injectivity problem and the reversibility problem for CA on $G$.
\end{theorem}

\begin{proof}
Let $(A,\FF,\delta)$ be a cycle-free strongly space-filling directed tileset over $(G,S)$. Recall that $R\colon A^G\to G$ is given by $R(x)=\delta(x(1_G))$, so that $R_x(g) = g\delta(x(g))$. This way, $\omega_x(n) = R^n_x(1_G)$. Consider an instance of the domino problem $(B,\G)$. For ease in the definition of the cellular automaton, we say a configuration $(x,y)\in (A\times B)^G$ is valid at $g$ if for every $s\in S$, the edge $\{g,gs\}$ is valid in $x$ for $\FF$, and valid in $y$ for $\G$. Now, let $\varphi$ be the CA on alphabet $\ZZ/2\ZZ\times A \times B$ given by the following rule. Given $(z,x,y)\in (\ZZ/2\ZZ\times A\times B)^G$ and $g\in G$, we let \[  \varphi(z,x,y)(g) = \begin{cases}
        (z(g)+z(R_x(g)),x(g),y(g)) & \mbox{if } (x,y) \textnormal{ is valid at } g\\
        (z(g),x(g),y(g)) & \mbox{otherwise. }
    \end{cases}\]
In simpler words, if $(x,y)$ is valid at $g$, we add to the bit at $g$ the bit at the next position given by the movement function. Otherwise, we do nothing. Clearly a description of the local rule of $\varphi$ is computable from $(B,\G)$.

We claim that $\varphi$ is bijective if $\X_{\G}$ is empty, and not injective otherwise.

If $\X_{\G}$ is nonempty, take $y\in\X_{\G}$ and $x\in\X_{\FF}$, which exists by the existence condition of~\Cref{def:SSF}. We have that $(x,y)$ is valid at every $g$, and get that $\varphi(0^G,x,y)=\varphi(1^G,x,y) = (0^G,x,y)$, thus $\varphi$ is not injective.

Now, suppose $\X_\G$ is empty. We first show that in this case for every $(x,y)\in(A\times B)^G$ and $g\in G$ there exists $n\in\N$ such that $(x,y)$ is not valid at $R_x^n(g)$. Suppose this is not the case, that is, there exists $(x,y)$ and $g_0\in G$ such that $(x,y)$ is valid at $R_x^n(g_0)$ for all $n\in\N$. Consider the configurations $x' = g_0^{-1}\cdot x$ and $y' = g_0^{-1}\cdot y$. As $R^n_{x'}(1_G) = g_0^{-1}R^n_x(g_0)$, the pair $(x',y')$ is valid at every element of $\omega_{x'}(\N)$. If there are two positions $i<j$ such that $\omega_{x'}(i) = \omega_{x'}(j)$, then the configuration $\omega_{x'}(i)^{-1}\cdot x'$ has a valid directed cycle. This is a contradiction to the cycle-free hypothesis. Thus, $n\mapsto R^n_{x'}(1_G)$ is injective and $x'$ is an infinite snake. Then, by the strongly space-filling property, for every $k\in\N$ there exists $h_k\in G$ such that $h_kB(k)\subseteq\omega_{x'}(\N)$. This implies that all edges in $B(k)$ are valid in $h_k^{-1}\cdot y'$ for $\G$. Taking an accumulation point of this sequence gives a point in $\X_\G$, which is a contradiction. 

Next, let us see $\varphi$ is injective. Because the automaton only changes the first coordinate, suppose there exist $z,z'$ with $z(g)\neq z'(g)$ for some $g\in G$ such that $\varphi(z,x,y) = \varphi(z',x,y)$. Notice that $(x,y)$ is valid at $g$, otherwise the images would differ. In particular, $z(R_x(g)) \neq z'(R_x(g))$. By repeating this argument we obtain that $(x,y)$ is valid at $R^n_x(g)$ for every $n\in\N$. This contradicts our previous claim. 

Finally, let us see $\varphi$ is surjective. Consider a configuration $(z,x,y)\in (\ZZ/2\ZZ\times A\times B)^G$. By the previous claim, for every $g\in G$ we may take $n_g\in\N$ as the least integer such that $(x,y)$ is not valid at $R^{n_g}_x(g)$. In particular, $n_{R_x(g)} = n_g -1$ if $n_g\geq 1$. Then, we define $z'(g)$ as $z(g)$ if $n_g = 0$, and $z(g) + z'(R_x(g))$ otherwise. A short computation then shows $\varphi(z',x,y) = (z,x,y)$.
\end{proof}

\begin{corollary}
      If $G$ has undecidable domino problem and admits a cycle-free strongly space-filling directed tileset over $(G,S)$, then the injectivity and reversibility problems for cellular automata on $G$ are undecidable.
  \end{corollary}

Next we look at the strong directed snake problem.

\begin{theorem}\label{thm:reduction}
If $G$ admits a strongly space-filling directed tileset over $(G,S)$, then $\DP(G,S)$ many-one reduces to $\SDSP(G,S)$.
\end{theorem}

\begin{proof}
Let $(A,\FF,\delta)$ be a strongly space-filling directed tileset over $(G,S)$, and consider an instance of $\DP(G,S)$ given by $(B,\G)$. We define the set of nearest neighbor forbidden patterns $\mathcal{H}$ over $A\times B$ given by
\[((a,b), (a',b'),s)\in\mathcal{H} \iff (a,a',s)\in \FF \ \vee \ (b,b',s)\in\G. \]
Furthermore, we define the directions $\delta'\colon A\times B \to S$ as $\delta'(a,b) = \delta(a)$.
Notice that the map that sends $(B,\G)$ to the instance $(A\times B, \mathcal{H}, \delta')$ of the strong directed snake problem is computable. 

Now, suppose there exists a configuration $y\in \X_{\G}\subseteq B^G$. By the strongly space-filling property, there exists an infinite snake $x\in \X_{\FF}\subseteq A^G$. Then, the configuration $(x,y)\in (A\times B)^G$ belongs to $\X_{\mathcal{H}}$ by definition. In addition, the skeleton for this configuration is $\omega_{(x,y)} = \omega_x$, therefore it is injective. Finally, every edge in $\omega_x(\N)$ is valid, and therefore, $(x,y)$ is an infinite snake.

Conversely, consider an infinite snake $(x,y)\in (A\times B)^G$ for $(A\times B, \mathcal{H},\delta')$. Because $\omega_x = \omega_{(x,y)}$, and every edge that is valid in $(x,y)$ with respect to $\mathcal{H}$ is valid in $x$ with respect to $\FF$, the configuration $x$ is an infinite snake for $(A,\FF, \delta)$. As the tileset is strongly space-filling, for every $n\in\N$ there exists $g_n\in G$ such that $g_nB(n)\subseteq \omega_{(x,y)}(\N)$. Because $(x,y)$ is an infinite snake, every edge in $\omega_{(x,y)}(\N)$ is valid. Therefore, all edges in $B(n)$ are valid for the configuration $y_n = g_n^{-1}y$. Taking a limit point of the sequence $(y_n)_{n\in\N}$, we obtain a configuration $y_*$ such that every edge is valid, that is, $\X_{\G}\neq\varnothing$.
\end{proof}

\begin{corollary}
    If $G$ has undecidable domino problem and admits a strongly space-filling directed tileset over $(G,S)$, then $\SDSP(G,S)$ is undecidable.
\end{corollary}

\section{Questions}
\label{sec:questions}

In this article we introduced the notion of space-filling actions and groups and provided some tools for their study, but there are still many natural questions that we were not able to answer.

In~\Cref{sec:space-filling} we established some basic facts about space-filling groups. In particular, we saw in~\Cref{prop:inf_ends} that groups with infinitely many ends are not space-filling. As virtually cyclic groups are finitary space-filling by~\Cref{cor:virtab_is_finitary_SF}, the remaining case is that of one-ended groups. We obtained some examples, such as finitely generated virtually abelian groups, groups quasi-isometric to hyperbolic spaces and some direct products, but the general question remains open. By a result of Seward~\cite{Seward2014}, every one-ended finitely generated group has a Cayley graph which admits a bi-infinite Hamiltonian path, but for us it is not clear if his construction can be used to construct a space-filling action.

\begin{question}
    Is there any finitely generated one-ended group which is not space-filling?
\end{question}

For the finitary and computable variants the answer is positive: by~\Cref{prop:rp_decidable}, a finitely generated and recursively presented group with undecidable word problem is neither finitary nor computably space-filling, and there are one-ended groups of this kind, such as $G\times\Z$ for a finitely presented group $G$ with undecidable word problem. We do not know whether these groups can be space-filling, nor whether there is a one-ended group with decidable word problem which is not computably space-filling.

An illustrative case is $F_2\times F_2$. We cannot use~\Cref{prop:producto}, as $F_2$ is not space-filling by~\Cref{prop:inf_ends}. Moreover, this group is self-simulable~\cite{BarSabSal2026}, so by~\Cref{cor:simulation_y_SF} it would be finitary space-filling as soon as it is computably space-filling.

\begin{question}
     Is $F_2\times F_2$ space-filling? Is it computably (and thus finitary) space-filling?
\end{question}

The effectiveness of space-filling actions also raises an important question. All our examples of space-filling groups with decidable word problem are also computably space-filling. It may be possible for such a group to be space-filling, but to only admit uncomputable covering functions for its space-filling actions. 

\begin{question}
    Is a space-filling finitely generated group with decidable word problem computably space-filling?
\end{question}

All known groups which are not finitary space-filling fail to be computably space-filling. This raises the following question.

\begin{question}
    Is there a finitely generated group with decidable word problem which is computably space-filling but not finitary space-filling?
\end{question}

A natural candidate is $F_2\times\Z$, which is not self-simulable, as shown in~\cite{BarSabSal2027_UFO}. While we expect it to be computably space-filling, we were not able to obtain a proof.

We introduced locally space-filling groups because the other notions do not suffice to establish the reductions of~\Cref{sec:applications}. In this article we did not produce any strongly space-filling directed tileset other than the one described in the work of Kari~\cite{Kari_reversibility_defofpeanocurve}, and we do not know if we can separate the groups that admit these objects from the finitary space-filling groups. This raises two questions.

\begin{question}
    If a group is finitary space-filling, is it locally space-filling as well?
\end{question}

\begin{question}
    Which groups are locally space-filling?
\end{question}

\printbibliography

\Addresses

\end{document}